# Stochastic analysis on configuration spaces: basic ideas and recent results

Michael Röckner
Fakultät für Mathematik
Universität Bielefeld
D 33615 Bielelefeld (Germany)

**Abstract**

The purpose of this paper is to provide a both comprehensive and summarizing account on recent results about analysis and geometry on configuration spaces $\Gamma_X$ over Riemannian manifolds $X$. Particular emphasis is given to a complete description of the so–called "lifting–procedure", Markov resp. strong resp. $L^1$–uniqueness results, the non–conservative case, the interpretation of the constructed diffusions as solutions of the respective classical "heuristic" stochastic differential equations, and a self–contained presentation of a general closability result for the corresponding pre–Dirichlet forms. The latter is presented in the general case of arbitrary (not necessarily pair) potentials describing the singular interactions. A support property for the diffusions, the intrinsic metric, and a Rademacher theorem on $\Gamma_X$, recently proved, are also discussed.

AMS Subject Classification
Primary: 60G60, 82C22

Key words and phrases: geometry of configuration spaces, mixed Poisson measures, canonical and grand canonical Gibbs measures, integration by parts, diffusion processes, stochastic dynamics, Dirichlet forms, irreducibility, ergodicity, extreme Gibbs measures, weak Sobolev spaces on configuration spaces, Rademacher theorem, intrinsic metric

# Contents

# 1 Introduction

The purpose of this paper is on one hand to provide sufficiently many details to complement lectures given in Anogia at the Euro–Conference on "Dirichlet forms and their Applications in Geometry and Stochastics" in June 1997 and at the Mathematical Sciences Research Institute in Berkeley in October/November 1997 within the "MSRI–Year in Stochastic Analysis". On the other hand the aim is to provide a both comprehensive and summarizing "pedagogic" account of the quite substantial material recently published on analysis and geometry on configuration spaces.

In addition, this work might serve as a guide for the reader to the following papers on the subject: [AKR97a], [AKR97b], (see also [AKR96a], [AKR96b]), [MR97], [RS97], [RSch97], [dSKR98]. We also include some new results and a number of new resp. modified proofs for crucial results in those articles, as well as describe some very recent developments, not covered by the lectures in Anogia and Berkeley.

For an overview of the single topics treated in these notes we refer to the list of contents which we hope is descriptive enough. Here we only would like to point out the following parts of this paper to which we put special emphasis resp. which are complements to the above articles:

(1) We describe the lifting of the geometry of the underlying manifold $X$ to the configuration space $\Gamma_X$ in Section 2 in more detail than it was done in [AKR97a].

(2) In Subsection 3.2 (see also Proposition 4.6) we give a particularly detailed account of the (Markov–, strong, $L^1$–) uniqueness results for the "mixed Poisson case" obtained in [AKR97a]. We complement [AKR97a] by additionally treating the *non–conservative case* (cf. Remark 3.8 (i)) in Subsection 3.2. In Subsection 3.3 (among other things also) the consequences of the uniqueness results for the corresponding diffusion processes on $\Gamma_X$ are discussed.

(3) In Subsection 8.2 in case the underlying measure $\mu$ on $\Gamma_X$ is a Gibbs (Ruelle) measure coming from a pair potential $\phi$, we include an explanation in what sense our diffusion process solves the *heuristic* stochastic differential equation

$$dX_t^i = dW_t^i + \sum_{j:j\neq i} \nabla\phi(X_t^i - X_t^j)\, dt,\ i \in \mathbb{N}.$$

(4) In Subsection 6.3 we give an essentially self–contained proof of the closability results recently obtained in [dSKR98].

(5) In Section 5 we include the joint results with Alexander Schied (who presented a part of them himself in Anogia and Berkeley) from [RSch97]. In particular, we discuss the Dirichlet forms $(\mathcal{E}_\mu^\Gamma, \mathcal{F})$ and $(\mathcal{E}_\mu^\Gamma, \mathcal{F}^{(c)})$ introduced in that reference (see Proposition 4.9 below), and additionally prove that the latter has the local property (cf. Proposition 5.6 below).

(6) The recent results in [RS97], not particularly mentioned in Anogia, are presented in Theorems 3.16 and 8.2 (iv).

(7) Most of the main results from [AKR97b] are discussed in these notes for simplicity under the additional *finite range* condition (C) (cf. Subsection 7.2). This is supposed to ease the reading for non–experts.

We emphasize, however, that the quite large number of consequences for mathematical physics is not discussed in sufficient detail in this paper. We refer instead to [AKR97c] and [KLRRSh97].

For the reader, who is only interested in the (much simpler) "free case", i.e., where the underlying measure on $\Gamma_X$ is a mixed Poisson measure, we tried to keep this part completely independent and have devoted Sections 2 and 3 entirely to this case. (The prize we had to pay in turn is however, that we have to be a bit repetitious in the more general "Gibbsian case" discussed in Sections 6 – 9). Also for simplicity, in Sections 2 and 3 we mainly consider the situation where the intensity measures $\sigma$ of the mixed Poisson measures is just the volume element $m$ of the manifold $X$ and leave the more general "$\sigma$–case" to respective remarks.

As far as prerequisites in general are concerned, we would like to point out that our framework is described in detail in Section 2. Some a-priori knowledge on Dirichlet forms is, however, advisable, e.g. Chapter I in [MR92] (though a *generator of a Dirichlet form* in this paper is $(-1)$ times the generator of a Dirichlet form in [MR92]).

W.r.t. the (anyway quite substantial) length of this paper we did not include all proofs of the results stated here. Main references for the single sections (apart from the modifications and complements mentioned above which are all completely proved below) are as follows: Sections 2 and 3: [AKR97a], [RS97]; Section 4: [AKR97b], [MR97]; Section 5: [RSch97]; Section 6: [AKR97b], [dSKR98]; Section 7–9: [AKR97b].

Finally, it is a great pleasure to thank all my colleagues and friends with whom I shared the exciting recent work on analysis and geometry on configuration spaces, in particular, Sergio Albeverio, Yuri Kondratiev, Zhi–Ming Ma, Jiagang Ren, Alexander Schied, Byron Schmuland and Jose Luis da Silva. I would also like to thank all of them for the permission to report here on very recent results, only just submitted for publication. I also thank the organizers of the above mentioned two conferences, in particular, Susanna Papadopoulou and Karl Theodor Sturm resp. Steve Evans and Ruth Williams for very nices stays in Anogia resp. Berkeley and for very stimulating meetings.

# 2 Lifting the geometry from the base manifold to the configuration space

In this paper let $X$ be a connected, oriented $C^\infty$ Riemannian manifold. For each point $x \in X$, the tangent space to $X$ at $x$ will be denoted by $T_xX$; and the tangent bundle endowed with its natural differentiable structure will be denoted $TX = \cup_{x\in X} T_xX$. The Riemannian metric on $X$ associates to each $x \in X$ an inner product on $T_xX$, which we denote by $\langle\cdot,\cdot\rangle_{TX}$. The associated norm will be denoted by $|\cdot|_{TX}$. Let $m$ denote the volume element.

$\mathcal{O}(X)$ is defined as the family of all open sets of $X$ and $\mathcal{B}(X)$ denotes the corresponding Borel $\sigma$–algebra. $\mathcal{O}_c(X)$ and $\mathcal{B}_c(X)$ denote the systems of all elements in $\mathcal{O}(X)$, $\mathcal{B}(X)$ respectively, which have compact closures.

The *configuration space* $\Gamma_X$ over the manifold $X$ is defined as the set of all locally finite subsets (configurations) in $X$:

$$\Gamma_X := \{\gamma \subset X \mid |\gamma \cap K| < \infty \text{ for any compact } K \subset X\}.$$

Here $|A|$ denotes the cardinality of a set $A$. For $\Lambda \subset X$ we sometimes use the shorthand $\gamma_\Lambda$ for $\gamma \cap \Lambda$ and define

$$\Gamma_\Lambda := \{\gamma \in \Gamma_X \mid \gamma \cap (X \setminus \Lambda) = \emptyset\}.$$

The reader should watch, however, that if $\Lambda$ is a sub–manifold not equal to $X$ there is an ambiguity in the notation (consider e.g. the case where $\Lambda$ is relatively compact). As in the standard literature we shall, however, never consider sub–manifolds of $X$ and always understand $\Gamma_\Lambda$ in the above sense.

We can identify any $\gamma \in \Gamma_X$ with a positive integer–valued Radon measure on $(X, \mathcal{B}(X))$, i.e.,

$$\gamma \equiv \sum_{x \in \gamma} \varepsilon_x \ ,$$

where $\sum_{x \in \emptyset} \varepsilon_x :=$ zero measure. The space $\Gamma_X$ can hence be endowed with the vague topology, i.e., the weakest topology on $\Gamma_X$ such that all maps

$$\Gamma_X \ni \gamma \mapsto \langle f, \gamma \rangle := \int_X f(x) \gamma(dx) = \sum_{x \in \gamma} f(x) \tag{2.1}$$

are continuous. Here $f \in C_0(X)$ (:= the set of all continuous functions on $X$ with compact support). Let $\mathcal{B}(X)$ denote the corresponding Borel $\sigma$–algebra.

For later use we recall the "localized" description of $\Gamma_X$: for $\Lambda \in \mathcal{B}_c(X)$ define

$$\Gamma_\Lambda \ := \ \{\gamma \in \Gamma_X \mid \gamma(X \setminus \Lambda) = 0\},$$

and for $n \in \mathbb{Z}_+$

$$\Gamma_\Lambda^{(n)} \ := \ \{\gamma \in \Gamma_\Lambda \mid \gamma(\Lambda) = n\}.$$

There is a bijection

$$\tilde{\Lambda}^n / S_n \ \to \Gamma_\Lambda^{(n)},$$

where

$$\tilde{\Lambda}^n \ := \ \{(x_1, \ldots, x_n) \in \Lambda^n | x_k \neq x_j, \ j \neq k\}, \tag{2.2}$$

and $S_n$ is the permutation group over $\{1,\dots,n\}$. If $\Lambda \in \mathcal{O}_c(X)$, this bijection defines a locally compact metrizably Hausdorff topology on $\Gamma_\Lambda^{(n)}$, hence a corresponding (sum) topology on

$$\Gamma_\Lambda = \bigcup_{n=0}^{\infty} \Gamma_\Lambda^{(n)}. \tag{2.3}$$

If $\Gamma_\Lambda$ is equipped with the associated Borel $\sigma-$ algebra then $(\Gamma_X, \mathcal{B}_X(\Gamma))$ is the projective limit of the measurable spaces $(\Gamma_\Lambda, \mathcal{B}(\Gamma_\Lambda))$ as $\Lambda \nearrow X$ (cf. e.g. [MR97]).

In this section we shall "lift the geometry" of $X$ onto $\Gamma_X$.

## 2.1 Test functions, flows, directional derivatives

We already know how to lift functions in the test functions space $\mathcal{D} = C_0^\infty(X)$ (:= the set of all $C^\infty$–functions on $X$ with compact support) onto $\Gamma_X$, namely according to (2.1) we obtain for each $f \in \mathcal{D}$ the function

$$\gamma \mapsto \langle f, \gamma \rangle$$

on $\Gamma_X$. We would like, however, an algebra of *test functions on* $\Gamma_X$. So, we define

$$\mathcal{F}C_b^\infty(\mathcal{D}, \Gamma_X) := \{g(\langle f_1, \cdot\rangle, \dots, \langle f_N, \cdot\rangle) \;\mid\; N \in \mathbb{N},\, g \in C_b^\infty(\mathbb{R}^N);\; f_1, \dots, f_N \in \mathcal{D}\} \tag{2.4}$$

and sometimes set for simplicity $\mathcal{F}C_b^\infty := \mathcal{F}C_b^\infty(\mathcal{D}, \Gamma_X)$. Elements in $\mathcal{F}C_b^\infty$ are called *finitely based*, since they only depend on finitely many points in $\gamma$.

After lifting the test functions, we want to lift flows. Let us first consider the group of diffeomorphisms $\mathrm{Diff}_0(X)$ on $X$ which are equal to the identity outside some compact set. Any $\psi \in \mathrm{Diff}_0(X)$ defines a transformation

$$\Gamma \ni \gamma \mapsto \psi(\gamma) = \{\psi(x) \mid x \in \gamma\} = \sum_{x \in \gamma} \varepsilon_{\psi(x)} \in \Gamma\ . \tag{2.5}$$

Correspondingly, if $V_0(X)$ denotes the set of smooth vector fields on $X$ with compact support and if $v \in V_0(X)$ with corresponding flow $\psi_t^v$, $t \in \mathbb{R}$ (i.e., $\frac{d}{dt}\psi_t^v(x) = v(\psi_t^v(x))$, $\psi_0^v(x) = x$), then $\psi_t^v$, $t \in \mathbb{R}$, lifts to a *flow on* $\Gamma_X$ via (2.3). Immediately, we now gain *directional derivatives on* $\Gamma_X$ for any $F \in \mathcal{F}C_b^\infty$ by defining for $v \in V_0(X)$

$$\left(\nabla_v^\Gamma F\right)(\gamma) := \frac{d}{dt}\, F(\psi_t^v(\gamma))|_{t=0}\ . \tag{2.6}$$

## 2.2 Tangent bundle, gradient, vector fields and divergence

The next step is to define a gradient for functions in $\mathcal{F}C_b^\infty$ which corresponds to the directional derivatives in (2.6). At the same time we shall find the appropriate "tangent bundle" over $\Gamma_X$.

So let $F = g(\langle f_1, \cdot\rangle, \ldots, \langle f_N, \cdot\rangle) \in \mathcal{F}C_b^\infty$, $v \in V_0(X)$ and $\gamma \in \Gamma_X$. By (2.6) and the chain rule we have with $\nabla_v^X f_i := \langle \nabla^X f_i, v\rangle_{TX}$

$$\begin{aligned}
\nabla_v^\Gamma F(\gamma) &= \sum_{i=1}^N \partial_i g(\langle f_1,\gamma\rangle, \ldots, \langle f_N,\gamma\rangle)\ \langle \nabla_v^X f_i, \gamma\rangle \\
&= \int \left\langle \sum_{i=1}^N \partial_i g(\langle f_1,\gamma\rangle, \ldots, \langle f_N,\gamma\rangle)\ \nabla^X f_i, v \right\rangle_{TX} d\gamma \\
&= \langle \nabla^\Gamma F(\gamma), v\rangle_{L^2(X\to TX;\gamma)},
\end{aligned} \tag{2.7}$$

where the *gradient* $\nabla^\Gamma$ *on* $\Gamma_X$ is defined by

$$\nabla^\Gamma F(\gamma) := \sum_{i=1}^N \partial_i g(\langle f_1,\gamma\rangle, \ldots, \langle f_N,\gamma\rangle)\ \nabla^X f_i \in C_0^\infty(X), \tag{2.8}$$

$\nabla^X$ denotes the gradient on $X$, $\partial_i$ directional derivative w.r.t. the $i$–th coordinate and $L^2(X \to TX;\gamma)$ the space of $\gamma$–square integrable vector fields on $X$. We note that here and henceforth we identify every $v$ in $V_0(X)$ with the corresponding class in $L^2(X \to TX;\gamma)$ though this is, of course, not a one–to–one operation. In particular, $\nabla^\Gamma F(\gamma) \in L^2(X \to TX;\gamma)$ and by (2.7) this class is independent of the representation of $F$ in (2.8).

Equation (2.7) immediately leads to the appropriate *tangent bundle on* $\Gamma_X$, namely

$$T_\gamma\Gamma_X := L^2(X \to TX;\gamma),\ \gamma \in \Gamma_X, \tag{2.9}$$

equipped with the usual $L^2$–inner product

$$\langle\ ,\ \rangle_{T_\gamma\Gamma_X} := \langle\ ,\ \rangle_{L^2(X\to TX;\gamma)},\ \gamma\in\Gamma_X. \tag{2.10}$$

Corresponding, *finitely based vector fields* on $(\Gamma_X, T\Gamma_X)$ can be defined as follows:

$$\Gamma_X \ni \gamma \mapsto \sum_{i=1}^N F_i(\gamma)\ v_i \in C_0^\infty(X), \tag{2.11}$$

where $F_1,\dots,F_N \in \mathcal{F}C_b^\infty$; $v_1,\dots,v_N \in V_0(X)$. Let $\mathcal{V}\mathcal{F}C_b^\infty := \mathcal{V}\mathcal{F}C_b^\infty(\mathcal{D},\Gamma_X)$ be the set of all such maps. We note that, of course, $\nabla^\Gamma F \in \mathcal{V}\mathcal{F}C_b^\infty$ for all $F\in\mathcal{F}C_b^\infty$ and that each $v\in V_0(X)$ is identified with the vector field $\gamma\mapsto v$ in $T\Gamma_X$ which is constant modulo taking $\gamma$–classes. Its norm function $\gamma\mapsto \|v\|_{T_\gamma\Gamma_X} = (\int \|v\|_{T_xX}\gamma(dx))^{1/2}$ is, however, not bounded.

We can now lift the divergence $\mathrm{div}^X$ on $X$. For $v\in V_0(X)$ by our rule to lift functions we must define

$$\mathrm{div}^\Gamma v(\gamma) := \langle \mathrm{div}^X v, \gamma\rangle,\ \gamma\in\Gamma_X\ . \tag{2.12}$$

Requiring the usual product rule to hold we must define for $V := \sum_{i=1}^N F_i\, v_i \in \mathcal{V}\mathcal{F}C_b^\infty$

$$\mathrm{div}^\Gamma V := \sum_{i=1}^N \left(\langle\nabla^\Gamma F_i, v_i\rangle_{T\Gamma_X} + F_i\,\mathrm{div}^\Gamma v_i\right)\ . \tag{2.13}$$

Also this will turn out to be a definition independent of the representation of $V$ (see Remark 2.3 (ii) below).

The next question that now arises is whether $\nabla^\Gamma$ and $\mathrm{div}^\Gamma$ are "connected" through a volume element on $\Gamma_X$ as $\nabla^X$ and $\mathrm{div}^X$ are through the volume element $m$ on $X$. Since $\Gamma_X$ is an infinite dimensional space, one might be sceptical whether such volume elements exist. But, in fact they do, if one uses the appropriate definition of "volume element" on $\Gamma_X$. Let us first recall the following well–known characterization of the volume element $m$ on $X$ (cf. e.g. [Ch84]):

$m$ is up to a constant the unique *Radon* measure $\sigma$ on $(X,\mathcal{B}(X))$ such that $(\nabla^X, C_0^\infty(X))$ and $(\mathrm{div}^X, V_0(X))$ are dual operators on $L^2(X;\sigma)$ w.r.t. $\langle\ ,\ \rangle_{TX}$, in the sense that

$$\int_X \langle v, \nabla^X f\rangle_{TX}\ d\sigma = -\int \mathrm{div}^X v\ f\ d\sigma \text{ for all } f\in C_0^\infty(X),\ v\in V_0(X).$$

**Definition 2.1** *A probability measure $\mu$ on $(\Gamma_X, \mathcal{B}(\Gamma_X))$ with Radon mean (i.e., $\int_{\Gamma_X} \gamma(K)\mu(d\gamma) < \infty$ for all compact $K \subset X$) is called a* volume element *on $\Gamma_X$, if*

$$\int \langle V, \nabla^\Gamma F\rangle_{T\Gamma_X}\, d\mu = -\int div^\Gamma V \;\, F\, d\mu \text{ for all } F \in \mathcal{F}C_b^\infty, V \in \mathcal{V}\mathcal{F}C_b^\infty. \tag{2.14}$$

## 2.3 Characterization of the volume elements as the mixed Poisson measures

We have a complete characterization of the volume elements on $\Gamma_X$. Let us first recall some standard definitions.

For any Radon measure $\sigma$ on $X$ the (*pure*) *Poisson measure with intensity* $\sigma$ is the unique measure $\pi_\sigma$ on $(\Gamma_X, \mathcal{B}(\Gamma_X))$ with Laplace transform given by

$$\int_{\Gamma_X} e^{\langle f,\gamma\rangle}\; \pi_\sigma(d\gamma) = e^{\int (e^f - 1)\, d\sigma} \quad \text{for all } f \in C_0^\infty(X). \tag{2.15}$$

The respective "local" description of $\pi_\sigma$ is as follows: for $\Lambda \in \mathcal{O}_c(X)$ and $n \in \mathbb{N}$ the product measure $\sigma^{\otimes n}$ can be considered as a (finite) measure on $\tilde{\Lambda}_n$. Let

$$\sigma_{\Lambda,n} \;:= \sigma^{\otimes n} \circ (s_\Lambda^n)^{-1} \tag{2.16}$$

be the corresponding measure on $\Gamma_\Lambda^{(n)}$ where $s_\Lambda^n : \tilde{\Lambda}^n \to \Gamma_\Lambda^{(n)}$, $s_\Lambda^n((x_1, \ldots, x_n)) := \sum_{i=1}^n \varepsilon_{x_i}$. Define $p_\Lambda : \Gamma \to \Gamma_\Lambda$ by

$$p_\Lambda(\gamma) \;:=\; \gamma_\Lambda$$

Then

$$\pi_\sigma \circ p_\Lambda^{-1} \;=\; e^{-\sigma(\Lambda)} \sum_{n=0}^\infty \frac{1}{n!}\, \sigma_{\Lambda,n}, \tag{2.17}$$

where $\sigma_{\Lambda,0} := \varepsilon_\emptyset$ on $\Gamma_\Lambda^{(0)} = \{\emptyset\}$.

(2.17) also leads to an existence proof of $\pi_\sigma$ as follows: Defining probability measures $\mu_\Lambda$ on $(\Gamma_\Lambda, \mathcal{B}(\Gamma_\Lambda))$ by the right hand side of (2.17), a version of Kolmogorov's theorem (cf. e.g. [Pa67]) implies that there exists a unique probability measure $\pi_\sigma$ on $(\Gamma_X, \mathcal{B}(\Gamma_X))$ such that

$$\pi_\sigma \circ p_\Lambda^{-1} = \mu_\Lambda \text{ for all } \Lambda \in \mathcal{O}_c(X).$$

It is easy to check that $\pi_\sigma$ indeed satisfies (2.15).

E.g. by (2.17) it is easy to check that

$$\int_{\Gamma_X} \langle f, \gamma \rangle \; \pi_\sigma(d\gamma) = \int_X f \; d\sigma \text{ for all } f \in L^1(X;\sigma) \tag{2.18}$$

and

$$\int_{\Gamma_X} \langle f, \gamma \rangle^2 \; \pi_\sigma(d\gamma) = \int_X f^2 \; d\sigma + \left( \int_X f \; d\sigma \right)^2 \text{ for all } f \in L^1(X;\sigma) \cap L^2(X;\sigma) \tag{2.19}$$

We note that for $\sigma \equiv 0$, $\pi_\sigma$ is just the Dirac measure on $(\Gamma_X, \mathcal{B}(\Gamma_X))$ with total mass in the empty configuration $\gamma \equiv \emptyset \; (\subset X)$.

For a probability measure $\lambda$ on $(\mathbb{R}_+, \mathcal{B}(\mathbb{R}_+))$ the measure $\mu_{\lambda,\sigma}$ on $(\Gamma_X, \mathcal{B}(\Gamma_X))$ defined by

$$\mu_{\lambda,\sigma} := \int_{\mathbb{R}_+} \pi_{z\cdot\sigma} \; \lambda(dz),$$

is called *mixed Poisson measure.*

**Theorem 2.2** *Suppose $m(X) = \infty$ and let $\mu$ be a probability measure on $(\Gamma_X, \mathcal{B}(\Gamma_X))$ with Radon mean. Then the following assertions are equivalent:*

*(i) $\mu$ is a volume element on $\Gamma_X$.*

*(ii) $\mu$ is a mixed Poisson measure $\mu_{\lambda,m}$ for some probability measure $\lambda$ on $(\mathbb{R}_+, \mathcal{B}(\mathbb{R}_+))$.*

**Proof.** (ii) $\Rightarrow$ (i): (Here the assumption $m(X) = \infty$ is not needed.) We shall first prove that $\pi_m$ is a volume element on $\Gamma_X$.

Let $F = g_F(\langle f_1, \cdot \rangle, \ldots, \langle f_N, \cdot \rangle) \in \mathcal{F}C_b^\infty$, $V \in \mathcal{V}\mathcal{F}C_b^\infty$. By linearity we may assume that $V(\gamma) = G(\gamma)v$ for all $\gamma \in \Gamma_X$ for some $G = g_G(\langle g_1, \cdot \rangle, \ldots, \langle g_N, \cdot \rangle)$, $v \in V_0(X)$. Let $\Lambda \in \mathcal{O}_c(X)$ such that $\operatorname{supp} f_i$, $\operatorname{supp} g_i$, $\operatorname{supp} v \subset \Lambda$, $1 \le i \le N$. Then by (2.6), (2.7), (2.16), and (2.17)

$$\begin{aligned} \int_{\Gamma_X} \langle V, \nabla^\Gamma F \rangle_{T\Gamma_X} \; d\pi_m &= \int_{\Gamma_X} G \; \nabla_v^\Gamma F \; d\pi_m \\ &= \int_{\Gamma_X} G(\gamma_\Lambda) \; \nabla_v^\Gamma F(\gamma_\Lambda) \; \pi_m(d\gamma) \\ &= e^{-\sigma(\Lambda)} \sum_{n=0}^\infty \frac{1}{n!} \; A_{\Lambda,n} \end{aligned} \tag{2.20}$$

where

$$
\begin{aligned}
A_{\Lambda,n} &:= \int_{\Lambda^n} g_G\left(\sum_{j=1}^n g_1(x_j),\ldots,\sum_{j=1}^n g_N(x_j)\right)\cdot \\
&\qquad \sum_{i=1}^N\left(\partial_i\, g_F\left(\sum_{j=1}^n f_1(x_j),\ldots,\sum_{j=1}^n f_N(x_j)\right)\ \sum_{k=1}^n \nabla_v^X f_i(x_k)\right) \\
&\qquad\qquad m(dx_1)\ldots m(dx_n) \\
&= \sum_{k=1}^n \int_{\Lambda^n} g_G\left(\sum_{j=1}^n g_1(x_j),\ldots,\sum_{j=1}^n g_N(x_j)\right)\cdot \\
&\qquad \left\langle \nabla_{x_k}^X g_F\left(\sum_{j=1}^n f_1(x_j),\ldots,\sum_{j=1}^n f_N(x_j)\right), v(x_k)\right\rangle_{T_{x_k}X} \\
&\qquad\qquad m(dx_1)\ldots m(dx_n)
\end{aligned}
$$

where $\nabla_{x_k}^X$ denotes the gradient on $X$ acting w.r.t. the variable $x_k$. Integrating by parts w.r.t. $x_k$ we obtain that

$$
\begin{aligned}
A_{\Lambda,n} &= -\sum_{k=1}^n \int_{\Lambda^n} \Big[ \left\langle \nabla_{x_k}^X g_G\left(\sum_{j=1}^n g_1(x_j),\ldots,\sum_{j=1}^n g_N(x_j)\right), v(x_k)\right\rangle_{T_{x_k}X} \\
&\qquad + g_G\left(\sum_{j=1}^n g_1(x_j),\ldots,\sum_{j=1}^n g_N(x_j)\right)\ \mathrm{div}^X v(x_k)\Big] \\
&\qquad \cdot g_F\left(\sum_{j=1}^n f_1(x_j),\ldots,\sum_{j=1}^n f_N(x_j)\right)\ m(dx_1)\ldots m(dx_n) \\
&= -\int_{\Lambda^n} \Big[\sum_{i=1}^N \partial g_G\left(\sum_{j=1}^n g_1(x_j),\ldots,\sum_{j=1}^n g_N(x_j)\right)\ \sum_{k=1}^n \nabla_v^X g_i(x_k) \\
&\qquad + g_G\left(\sum_{j=1}^n g_1(x_j),\ldots,\sum_{j=1}^n g_N(x_j)\right)\sum_{k=1}^n \mathrm{div}^X v(x_k)\Big] \\
&\qquad g_F\left(\sum_{j=1}^n f_1(x_j),\ldots,\sum_{j=1}^n f_N(x_j)\right)\ m(dx_1)\ldots m(dx_N).
\end{aligned}
$$

Substituting back into (2.20) by (2.16), (2.17) yields

$$\int \langle V, \nabla^\Gamma F\rangle_{T\Gamma_X}\ d\pi_m = -\int \mathrm{div}^\Gamma V\ F\ d\pi_m\ .$$

Since we have not used a specific normalization of $m$, this proves that any $z \cdot m$, $z \in ]0, \infty[$, is a volume element on $\Gamma_X$. It is trivial to check that this is also true for $z = 0$, since $\pi_{0\cdot m}$ is the Dirac measure on $(\Gamma_X, \mathcal{B}(\Gamma_X))$ with total mass in the empty configuration $\emptyset \subset X$. Hence by Fubini's theorem this is also true for $\mu_{\lambda,m}$ as in the assertion, since integrability is ensured because $\mu_{\lambda,m}$ has a Radon mean (hence $\lambda$ has finite first moments).

For the proof of the converse (i) $\Rightarrow$ (ii) (which is a little harder) we refer to the detailed expositions in [AKR97a] and [AKR97b], more precisely Theorem 3.2, 4.1 respectively ibidem. (The first is analytic, the latter is purely probabilistic and works in more general situations, see the following remark.) $\square$

**Remark 2.3** (i) ("*$\sigma$–case*") Let $\rho \in L^1_{loc}(X; m)$ be such that $\rho^{\frac{1}{2}} \in H^{1,2}_{loc}(X; m)$, i.e., $\rho^{\frac{1}{2}}$ has locally $m$–square integrable weak derivatives. Define

$$\sigma := \rho \cdot m\ .$$

the *logarithmic derivative* of the measure $\sigma$ is given by the vector field

$$\beta^\sigma := \frac{\nabla^X \rho}{\rho} \tag{2.21}$$

where $\beta^\sigma := 0$ on $\{\rho = 0\}$. For every $v \in V_0(X)$ we set

$$\mathrm{div}^X_\sigma v := \langle \beta^\sigma, v\rangle_{TX} + \mathrm{div}^X v\ . \tag{2.22}$$

Note that $\langle \beta^\sigma, v\rangle \in L^2(X; \sigma) \cap L^1(X; \sigma)$ since it has compact support, hence the same is true for $\mathrm{div}^X_\sigma v$; consequently by (2.19)

$$\langle \mathrm{div}^X_\sigma v, \cdot\rangle \in L^2(\Gamma_X; \pi_\sigma)\ .$$

It is now easy to check that in the proof of Theorem 2.2, (ii) $\Rightarrow$ (i), $m$ can be replaced by $\sigma$ if in Definition 2.1 $\mathrm{div}^X$ is replaced by $\mathrm{div}^X_\sigma$. If $\sigma(X) = \infty$ the same is true for the converse (i) $\Rightarrow$ (ii), if, in addition, e.g. $|\beta^\sigma|_{TX} \in L^1_{loc}(X; m)$. We refer to [AKR97b, Theorem 4.1 and Remark 4.1] for details.

(ii) Since the left hand side of (2.14) is independent of the representation of $V \in \mathcal{V}\mathcal{F}C_b^\infty$ as $V = \sum_{i=1}^N F_i \cdot v_i$; $F_1, \ldots, F_N \in \mathcal{F}C_b^\infty$, $v_1, \ldots, v_N \in V_0(X)$, so is $\mathrm{div}^\Gamma V \in L^1(\Gamma_X; \mu)$ for any volume element $\mu$ on $\Gamma_X$, hence for any mixed Poisson measure $\mu_{\lambda,m}$ as in Theorem 2.2.

(iii) A proof of Theorem 2.2, (ii) $\Rightarrow$ (i), can also be derived from the well–known Mecke identity for $\pi_\sigma$ (cf. [Mec67, Satz 3.1]). This was observed and pointed out to us by V. Liebscher (private communication).

(iv) Theorem 2.2, (ii) $\Rightarrow$ (i), extends to so–called compound Poisson measures (cf. [dSKSt97]). The proof is similar.

## 2.4 Laplacian and classical pre–Dirichlet form

Now we proceed with our lifting operations. Since we have gradient and divergence on $\Gamma_X$, it is clear how to define the corresponding *Laplacian on* $\Gamma_X$, namely for $F \in \mathcal{F}C_b^\infty$

$$\Delta^\Gamma F := \mathrm{div}^\Gamma \ \nabla^\Gamma F \ . \tag{2.23}$$

(We here note that $\nabla^\Gamma F \in \mathcal{V}\mathcal{F}C_b^\infty$).

By (2.8), (2.13) it follows that for all $F = g(\langle f_1, \cdot\rangle, \ldots, \langle f_N, \cdot\rangle) \in \mathcal{F}C_b^\infty$

$$\begin{aligned} \Delta^\Gamma F(\gamma) \ &= \ \sum_{i,j=1}^N \partial_i \partial_j g(\langle f_1, \gamma\rangle, \ldots, \langle f_N, \gamma\rangle) \ \left\langle \langle \nabla^X f_i, \nabla^X f_j \rangle_{TX}, \gamma \right\rangle \\ &\quad + \sum_{i=1}^N \partial_i g(\langle f_1, \gamma\rangle, \ldots, \langle f_N, \gamma\rangle) \ \left\langle \Delta^X f_i, \gamma \right\rangle . \end{aligned} \tag{2.24}$$

Here $\Delta^X := \mathrm{div}^X \, \nabla^X$ is the Laplacian on $X$.

For any volume element $\mu_{\lambda,m}$ (cf. Theorem 2.2) we obtain the *classical pre–Dirichlet form on* $\Gamma_X$ as

$$\mathcal{E}^\Gamma_{\mu_{\lambda,m}}(F, G) := \int \langle \nabla^\Gamma F, \nabla^\Gamma G \rangle_{T\Gamma_X} \, d\mu_{\lambda,m}; \ F, G \in \mathcal{F}C_b^\infty. \tag{2.25}$$

"pre" is added because the form is not yet closed on $L^2(\Gamma_X; \mu_{\lambda,m})$. We come to the "problem of closability" in detail in Subsections 4.1 and 6.3 below (and also in the following subsection 3.1). We only emphasize here that only by passing to the closure of $(\mathcal{E}^\Gamma_{\mu_{\lambda,m}}, \mathcal{F}C_b^\infty)$ on $L^2(\Gamma_X; \mu_{\lambda,m})$ we start performing real infinite dimensional analysis on $\Gamma_X$. So far (apart from the construction of Poisson measures) we have been doing only essentially finite dimensional analysis, since our test function and vector fields on $\Gamma_X$ are finitely based, hence depend only on a finite part of each configuration $\gamma \in \Gamma_X$. We note that until now we also never used the full tangent space $T_\gamma\Gamma_X = L^2(X \to TX; \gamma)$ at $\gamma \in \Gamma_X$, but rather only $V_0(X)$ with the inherited inner product.

**Remark 2.4** The pre–Dirichlet form (2.25) seems to have appeared first in [BiGJ87] in the special case where $X$ is an open subset of $\mathbb{R}^d$ and with $\pi_\sigma$ replacing $\mu_{\lambda,m}$. But $T\Gamma_X$ and $\nabla^\Gamma$ were neither identified nor used there. Instead, $\langle \nabla^\Gamma F, \nabla^\Gamma G\rangle_{T\Gamma_X}$ was replaced by the corresponding formula obtained by means of (2.8). $\nabla^\Gamma$ and $T\Gamma_X$ at least in the special case $X = \mathbb{R}^1$ seem to have first been introduced in [Sm88]. We thank F. Hirsch resp. V.I. Bogachev for pointing this out to us.

# 3 Infinite dimensional analysis and Brownian motion on configuration spaces

In this section we shall develop a “real” infinite dimensional analysis on $\Gamma_X$ and discuss some applications. For the rest of this paper we shall use all terminology and notations from Section 2 without further notice. In particular, we fix $\mu_{\lambda,m}$ as in Theorem 2.2. It can be shown that $\mu_{\lambda,m}$ has full support on $\Gamma_X$, i.e., $\mu_{\lambda,m}(U) > 0$ for every open subset $U$ of $\Gamma_X$ (cf. [RSch97]). Therefore, we do not distinguish between $\mathcal{F}C_b^\infty$ and the corresponding $\mu_{\lambda,m}$ classes, since each of the latter has thus a *unique* representative in $\mathcal{F}C_b^\infty$. We additionally assume that

$$\int_{\mathbb{R}_+} z^2 \ \lambda(dz) < \infty \,. \tag{3.1}$$

It follows by (2.18), (2.19) that

$$\int_{\Gamma_X} \langle f, \gamma\rangle \ \mu_{\lambda,m}(d\gamma) = \int_{\mathbb{R}_+} z \ \lambda(dz) \ \int_X f \ dm \quad \text{for all } f \in L^1(X; m). \tag{3.2}$$

and

$$\int_{\Gamma_X} \langle f, \gamma\rangle^2 \, \mu_{\lambda,m}(d\gamma) \;=\; \int_{\mathbb{R}_+} z \ \lambda(dz) \ \int_X f^2 \, dm + \int_{\mathbb{R}_+} z^2 \, \lambda(dz) \left(\int_X f \, dm\right)^2 \quad \text{for all } f \in L^1(X; m) \cap L^2(X; m). \tag{3.3}$$

### 3.1 Dirichlet forms and operators

For the definitions of closability and closure we refer to Subsection 4.1 below.

**Proposition 3.1** *For all $F, G \in \mathcal{F}C_b^\infty$*

$$\mathcal{E}^\Gamma_{\mu_{\lambda,m}}(F,G) = -\int \Delta^\Gamma F \; G \; d\mu_{\lambda,m} \; . \tag{3.4}$$

*In particular, $(\mathcal{E}^\Gamma_{\mu_{\lambda,m}}, \mathcal{F}C_b^\infty)$ is closable on $L^2(\Gamma_X; \mu_{\lambda,m})$.*

**Proof.** (3.4) is immediate because $\mu_{\lambda,m}$ is a volume element on $\Gamma_X$ by Theorem 2.2. The closability follows then as a special case from [MR92, Chap. I, Proposition 3.3], since $\Delta^\Gamma F \in L^2(\Gamma_X; \mu_{\lambda,m})$ by (3.3) and (2.24). □

Let $(\mathcal{E}^\Gamma_{\mu_{\lambda,m}}, D(\mathcal{E}^\Gamma_{\mu_{\lambda,m}}))$ be the closure of $(\mathcal{E}^\Gamma_{\mu_{\lambda,m}}, \mathcal{F}C_b^\infty)$ on $L^2(\Gamma_X; \mu_{\lambda,m})$. $D(\mathcal{E}^\Gamma_{\mu_{\lambda,m}})$ is the analogue of a first order Sobolev space on $\Gamma_X$. We, therefore, set

$$H_0^{1,2}(\Gamma_X; \mu_{\lambda,m}) := D(\mathcal{E}^\Gamma_{\mu_{\lambda,m}}) \tag{3.5}$$

with norm

$$\| \cdot \|_{H_0^{1,2}(\Gamma_X;\mu_{\lambda,m})} := \Big( \mathcal{E}^\Gamma_{\mu_{\lambda,m}}(\cdot,\cdot) + (\cdot,\cdot)_{L^2(\Gamma_X;\mu_{\lambda,m})} \Big)^{1/2} . \tag{3.6}$$

$(\mathcal{E}_{\mu_{\lambda,m}}, D(\mathcal{E}^\Gamma_{\mu_{\lambda,m}}))$ is a closed, positive definite, symmetric bilinear form on $L^2(\Gamma_X; \mu_{\lambda,m})$, hence it is associated with a unique positive definite self–adjoint operator $(H^\Gamma_{\mu_{\lambda,m}}, D(H^\Gamma_{\mu_{\lambda,m}}))$ (called its *generator*) on $L^2(\Gamma_X; \mu_{\lambda,m})$, i.e.,

$$D\big(\sqrt{H^\Gamma_{\mu_{\lambda,m}}}\big) = D(\mathcal{E}^\Gamma_{\mu_{\lambda,m}}), \; \mathcal{E}^\Gamma_{\mu_{\lambda,m}}(F,G) = \Big( \sqrt{H^\Gamma_{\mu_{\lambda,m}}} F, \sqrt{H^\Gamma_{\mu_{\lambda,m}}} G \Big)_{L^2(\Gamma_X;\mu_{\lambda,m})}$$
$$\text{for all } F, G \in D(\mathcal{E}^\Gamma_{\mu_{\lambda,m}}) \tag{3.7}$$

(cf. e.g. [MR92, Chap. I] for details, except that the generator there is $(-1)$ times our generator here).

Clearly, by (3.4)

$$H^\Gamma_{\mu_{\lambda,m}} F = -\Delta^\Gamma F \quad \text{for all } F \in \mathcal{F}C_b^\infty,$$

i.e., $(H^\Gamma_{\mu_{\lambda,m}}, D(H^\Gamma_{\mu_{\lambda,m}}))$ extends $(-\Delta^\Gamma, \mathcal{F}C_b^\infty)$. It is nothing but the Friedrichs' extension of $(-\Delta^\Gamma, \mathcal{F}C_b^\infty)$ on $L^2(\Gamma_X; \mu_{\lambda,m})$ and is also called *Dirichlet operator.* In fact it is also a Dirichlet operator in the sense of [MR92, Chap. I, Definition 4.1]. This follows from [MR92, Chap. I, Theorem 4.4, Proposition 4.3] and the following result:

**Proposition 3.2** *$(\mathcal{E}^{\Gamma}_{\mu_{\lambda,m}}, D(\mathcal{E}^{\Gamma}_{\mu_{\lambda,m}}))$ is a symmetric* Dirichlet form *on $L^2(\Gamma_X; \mu_{\lambda,m})$, i.e.,*

$$F \in D(\mathcal{E}^{\Gamma}_{\mu_{\lambda,m}}) \;\Rightarrow\; F^{\#} := F^{+} \wedge 1 \in D(\mathcal{E}^{\Gamma}_{\mu_{\lambda,m}}) \text{ and } \mathcal{E}^{\Gamma}_{\mu_{\lambda,m}}(F^{\#}, F^{\#}) \leq \mathcal{E}^{\Gamma}_{\mu_{\lambda,m}}(F, F) \; . \tag{3.8}$$

**Proof.** The argument is entirely standard, since $\nabla^{\Gamma}$ satisfies the chain rule. (One proceeds exactly in the same way as e.g. in [MR92, Chap. II, Subsection 2c)].) □

E.g. by the spectral theorem $(H^{\Gamma}_{\mu_{\lambda,m}}, D(H^{\Gamma}_{\mu_{\lambda,m}}))$ generates a strongly continuous contraction semigroup $T^{\Gamma}_{\mu_{\lambda,m}}(t) := e^{-tH^{\Gamma}_{\mu_{\lambda,m}}}$, $t > 0$, on $L^2(\Gamma_X; \mu_{\lambda,m})$. Since (3.8) is equivalent to the *sub–Markov property* of $(T^{\Gamma}_{\mu_{\lambda,m}}(t))_{t>0}$ (i.e., $F \in L^2(\Gamma_X; \mu_{\lambda,m})$, $0 \leq F \leq 1 \Rightarrow 0 \leq T^{\Gamma}_{\mu_{\lambda,m}}(t)F \leq 1$ for all $t \geq 0$, cf. [MR92, Chap. II, Sect. 4]), it follows that $(T^{\Gamma}_{\mu_{\lambda,m}}(t))_{t>0}$ is a strongly continuous contraction semigroup on every $L^p(\Gamma_X; \mu_{\lambda,m})$, $p \in [1, \infty[$ (cf. [RS75, Theorem X.59]), which we shall denote by the same symbol. The corresponding generators we denote by $(H^{\Gamma}_{\mu_{\lambda,m},p}, D(H^{\Gamma}_{\mu_{\lambda,m},p}))$. We have in fact that $(T^{\Gamma}_{\mu_{\lambda,m}}(t))_{t\geq 0}$ is *Markovian*, i.e., in addition, it holds that $T^{\Gamma}_{\mu_{\lambda,m}}(t)\,1 = 1$ for all $t > 0$ (which follows form the fact that $1 \in \mathcal{F}C_b^{\infty} \subset D(H^{\Gamma}_{\mu_{\lambda,m}})$ and $H^{\Gamma}_{\mu_{\lambda,m}}\,1 = \Delta^{\Gamma}\,1 = 0$).

## 3.2 Markov uniqueness, essential self–adjointness and heat semigroup

We emphasize that in general the Friedrichs' extension is not the only self–adjoint extension of a symmetric operator on $L^2$. In Section 5 we shall consider Dirichlet forms $(\mathcal{E}, D(\mathcal{E}))$ on $L^2(\Gamma_X; \mu_{\lambda,m})$ which are a–priori extensions of $(\mathcal{E}^{\Gamma}_{\mu_{\lambda,m}}, D(\mathcal{E}^{\Gamma}_{\mu_{\lambda,m}}))$ and whose generators also extend $(\Delta^{\Gamma}, \mathcal{F}C_b^{\infty}))$. So, it is a non-trivial question whether these extensions are *strict* extensions or indeed coincide with $(\mathcal{E}^{\Gamma}_{\mu_{\lambda,m}}, D(\mathcal{E}^{\Gamma}_{\mu_{\lambda,m}}))$. Below we shall give a complete answer to this question and prove that, in fact, they all coincide (at least under conditions (A), (B) specified in Theorem 3.3). To this end we introduce the set $\underline{\underline{\mathcal{E}}}$ of all Dirichlet forms $(\mathcal{E}, D(\mathcal{E}))$ on $L^2(\Gamma_X; \mu_{\lambda,m})$ whose respective generators $(L(\mathcal{E}), D(L(\mathcal{E})))$ extend $(-\Delta^{\Gamma} F, \mathcal{F}C_b^{\infty})$. We note that $(\mathcal{E}^{\Gamma}_{\mu_{\lambda,m}}, D(\mathcal{E}^{\Gamma}_{\mu_{\lambda,m}})) \in \underline{\underline{\mathcal{E}}}$ and that any $(\mathcal{E}, D(\mathcal{E})) \in \underline{\underline{\mathcal{E}}}$ then automatically extends $(\mathcal{E}^{\Gamma}_{\mu_{\lambda,m}}, D(\mathcal{E}^{\Gamma}_{\mu_{\lambda,m}}))$,

since for all $F, G \in \mathcal{F}C_b^\infty$ ($\subset D(L(\mathcal{E}))$)

$$\begin{aligned}\mathcal{E}(F,G) &= \int L(\mathcal{E})F\; G\; d\mu_{\lambda,m} = -\int \Delta^\Gamma F\; G\; d\mu_{\lambda,m} \\ &\underset{(3.4)}{=} \mathcal{E}^\Gamma_{\mu_{\lambda,m}}(F,G),\end{aligned}$$

hence this equality extends to the closure. In this sense $(\mathcal{E}^\Gamma_{\mu_{\lambda,m}}, D(\mathcal{E}^\Gamma_{\mu_{\lambda,m}}))$ is minimal in $\underline{\underline{\mathcal{E}}}$ and the extensions mentioned above to be studied in Section 5 are in $\underline{\underline{\mathcal{E}}}$. If $\underline{\underline{\mathcal{E}}} = \{(\mathcal{E}^\Gamma_{\mu_{\lambda,m}}, D(\mathcal{E}^\Gamma_{\mu_{\lambda,m}}))\}$ one says that *Markov uniqueness* holds which we shall prove to be the case in many situations. Before we state the corresponding theorem and give a complete proof, we recall that a symmetric operator on a Hilbert space is called *essentially self–adjoint* if its closure is self–adjoint. Let $(H_m^X, D(H_m^X))$ be the Friedrichs' extension of $(\Delta^X, \mathcal{D})$ on $L^2(X;m)$.

**Theorem 3.3** *Assume that the following conditions hold:*

(A) $(\Delta^X, \mathcal{D})$ *is essentially self–adjoint on* $L^2(X;m)$ *(which, as is well-known, is e.g. the case, if X is complete w.r.t. the Riemannian metric).*

(B) $\lambda = \varepsilon_z$ *for some* $z \in [0,\infty[$ *or*

$$\int H_m^X f\; dm = 0 \qquad \text{for all } f \in D(H_m^X) \cap L^1(X;m) \quad \text{such that } H_m^X f \in L^1(X;m).$$

*(The latter just means that* $(H_m^X, D(H_m^X))$ *is* conservative*).*

*Then*

$$\underline{\underline{\mathcal{E}}} = \{(\mathcal{E}^\Gamma_{\mu_{\lambda,m}}, D(\mathcal{E}^\Gamma_{\mu_{\lambda,m}}))\},$$

*i.e., Markov uniqueness holds.*

**Remark 3.4** *("σ–case")* Let $\sigma = \rho \cdot m$ be as specified at the beginning of Remark 2.3 (i). Then Theorem 3.3 (as well as all subsequent results in Subsections 3.2, 3.3) extend to the case where $m$ is replaced by $\sigma$. The proof is word by word the same as the one given below. We only note here that condition (A) for $\sigma$ is e.g. fulfilled if $X$ is complete and if $|\beta^\sigma|_{TX} \in L^p_{loc}(X;m)$ for some $p > \dim(X)$, where $\beta^\sigma := \frac{\nabla \rho}{\rho}$. This immediately follows from the proof of Theorem 1 and Remark 4 (iii) in [BKR97] and Gaffney's Lemma (cf. e.g. [Ba87]).

Theorem 3.3 is a consequence of the following result:

**Theorem 3.5** *Suppose that conditions (A), (B) in Theorem 3.3 hold. Then* $\mathcal{F}C_b^\infty$ *is dense in* $D(H^\Gamma_{\mu_{\lambda,m},1})$ *w.r.t. the graph norm* $\|\cdot\|_1 := \|H_{\mu_{\lambda,m},1}\cdot\|_{L^1(\Gamma_X;\mu_{\lambda,m})} + \|\cdot\|_{L^1(X;\mu_{\lambda,m})}$.

Before we prove Theorem 3.5, we show that it implies Theorem 3.3:

**Proof of Theorem 3.3.** Let $(\mathcal{E}, D(\mathcal{E})) \in \underline{\underline{\mathcal{E}}}$ with corresponding generator $(L, D(L))$. Since $T_t := e^{-tL}$, $t > 0$, is Markovian, it extends to a strongly continuous contraction semigroup on $L^1(\Gamma_X;\mu_{\lambda,m})$. Let $(L_1, D(L_1))$ denote the corresponding generator. Since $\mu_{\lambda,m}$ has finite total mass, also $(L_1, D(L_1))$ extends $(-\Delta, \mathcal{F}C_b^\infty)$. Therefore, it extends its $L^1$–closure, which coincides by Theorem 3.5 with $(H^\Gamma_{\mu_{\lambda,m},1}, D(H^\Gamma_{\mu_{\lambda,m},1}))$. But one generator of a strongly continuous semigroup cannot strictly extend another, because this would contradict the fact that by the Hille–Yosida theorem both $1+L_1$ and $1+H^\Gamma_{\mu_{\lambda,m},1}$ have range all of $L^1(\Gamma_X;\mu_{\lambda,m})$. Thus, $(L_1, D(L_1)) = (H^\Gamma_{\mu_{\lambda,m},1}, D(H^\Gamma_{\mu_{\lambda,m},1}))$, hence $e^{-tL} = e^{-tH^\Gamma_{\mu_{\lambda,m}}}$ for all $t \geq 0$, consequently $(\mathcal{E}, D(\mathcal{E})) = (\mathcal{E}^\Gamma_{\mu_{\lambda,m}}, D(\mathcal{E}^\Gamma_{\mu_{\lambda,m}}))$. $\square$

We shall prove Theorem 3.5 along with two other related results. To formulate them we need some preparations.

Consider the classical Dirichlet form $(\mathcal{E}^X_m, D(\mathcal{E}^X_m))$ on $X$ defined as the closure of

$$\begin{aligned} \mathcal{E}^X_m(f,g) &:= \int_X (\nabla^X f, \nabla^X g)_{TX}\, dm \\ &= -\int_X \Delta^X f\; g\; dm \quad ;\; f,g \in \mathcal{D}. \end{aligned} \tag{3.9}$$

Let $(H^X_m, D(H^X_m))$ be the corresponding generator (i.e., the Friedrichs' extension of $(-\Delta^X, \mathcal{D})$ on $L^2(X;m)$ mentioned before). We note that if condition (A) holds $(H^X_m, D(H^X_m))$ coincides with the closure of $(-\Delta^X, \mathcal{D})$ since one self–adjoint operator cannot strictly extend another. Since also $(e^{-tH^X_m})_{t>0}$ is sub–Markovian, it is a strongly continuous contraction semigroup on every $L^p(X;m)$, $p \in [1,\infty[$ (cf. [Da89, Theorem 1.4.1]). Define

$$\begin{aligned} \mathcal{D}_1 := \{f \in D(H^X_m) \;\cap\; & L^1(X;m) \mid H^X_m f \in L^1(X;m) \text{ and} \\ & -\delta \leq f \leq 0 \text{ for some } \delta \in ]0,1[\} \end{aligned} \tag{3.10}$$

and for $\mathcal{D}_0 \subset \mathcal{D}_1$

$$E(\mathcal{D}_0) := \text{linear hull of} \{\exp(\langle \log(1+f), \cdot \rangle) \mid f \in \mathcal{D}_0\}. \tag{3.11}$$

**Remark 3.6** (i) It is an easy exercise to show that if $f \in D(H_m^X) \cap L^1(X;m)$ and $H_m^X f \in L^1(X;m)$, then $f \in D(H_{m,1}^X)$ and $H_{m,1}^X f = H_m^X f$ (cf. the remark following the proof of Theorem 1.4.1 in [Da89]). In particular, $\mathcal{D}_1 \subset D(H_{m,1}^X)$.

(ii) We note that obviously $e^{-tH_m^X} \mathcal{D}_1 \subset \mathcal{D}_1$.

**Proposition 3.7** *(i) Suppose condition (A) holds and let* $t \geq 0$, $f \in L^1(X;m)$, $-1 < f \leq 0$.
*If* $\lambda = \varepsilon_z$ *for some* $z \in [0, \infty[$, *then*

$$\begin{aligned} & T^\Gamma_{\pi_{z \cdot m}}(t) \left(\exp(\langle \log(1+f), \cdot \rangle)\right) \\ = \; & \exp\left( \langle \log(1 + e^{-tH_m^X} f), \cdot \rangle - z \int \left( e^{-tH_m^X} f - f \right) dm \right). \end{aligned}$$

*If* $(H_m^X, D(H_m^X))$ *is conservative, then*

$$T^\Gamma_{\mu_{\lambda,m}}(t) \left(\exp \langle \log(1+f), \cdot \rangle)\right) = \exp\left( \langle \log(1 + e^{-tH_m^X} f, \cdot \rangle \right).$$

*(ii) Suppose conditions (A) and (B) hold. Then* $E(\mathcal{D}_1)$ *is a dense subset of both* $D(H^\Gamma_{\mu_{\lambda,m}})$ *and* $D(H^\Gamma_{\mu_{\lambda,m},1})$ *w.r.t. the respective graph norms* $\| \cdot \|_2 := \left( \|H^\Gamma_{\mu_{\lambda,m}} \cdot \|^2_{L^2(\Gamma_X;\mu_{\lambda,m})} + \| \cdot \|^2_{L(\Gamma_X;\mu_{\lambda,m})} \right)^{1/2}$ *(i.e.,* $(H^\Gamma_{\mu_{\lambda,m}}, E(\mathcal{D}_1))$ *is essentially self–adjoint on* $L^2(\Gamma_X;\mu_{\lambda,m})$*) and* $\| \cdot \|_1$.
*(iii) Suppose conditions (A) and (B) hold. Then any function in* $E(\mathcal{D}_1)$ *can be approximated by a uniformly bounded sequence of functions in* $E(\mathcal{D}_1 \cap C_0^\infty(X))$ $(\subset \mathcal{F}C_b^\infty)$ *w.r.t.* $\| \cdot \|_1 + \| \cdot \|_{H_0^{1,2}(\Gamma_X;m)}$.

Clearly, Theorem 3.5 is an immediate consequence of Proposition 3.7 (ii) and (iii). so, it remains to prove 3.7 whose parts (ii), (iii) are essentially consequences of its first part (i).

**Remark 3.8** (i) In the preliminary preprint version of [AKR97a] it was claimed that assertions (i) and (ii) of Proposition 3.7 also hold if $(H_m^X, D(H_m^X))$ is not conservative and $\mu_{\lambda,m}$ is not (necessarily) a pure Poisson measure. However, there was a gap in the proof. This case is so far unproved.

(ii) In fact if (B) holds, the two formulae in Proposition 3.7 (i) coincide, because (B) implies that

$$\int \left( e^{-tH_m^X} f - f \right) dm = 0 \quad \text{for all } t > 0, f \in L^1(X;m).$$

(iii) If $\lambda = \varepsilon_z$ for some $z \in [0, \infty[$, then it can be shown using second quantization that $\mathcal{F}C_b^\infty$ is also dense in $D(H^\Gamma_{\pi_z \cdot m})$ w.r.t. $\| \cdot \|_2$ (cf. [AKR97a, Subsection 5.4]).

For the proof of Proposition 3.7 we need a number of lemmas. We first recall the following well–known result:

**Lemma 3.9** *Let $\mu$ be a $\sigma$–finite measure on a measurable space $(E, \mathcal{B})$ and $p \in [1, \infty)$. Let $\chi \in C_b^1(\mathbb{R})$ and let $\mathbb{R}_+ \ni t \mapsto f_t \in L^p(E;\mu)$ be a continuously differentiable map. Then so is $\mathbb{R}_+ \ni t \mapsto \chi \circ f_t \in L^p(E;\mu)$ and $\frac{d(\chi \circ f_t)}{dt} = \chi'(\psi_t) \frac{df_t}{dt}$.*

**Proof.** Since $L^p(E;\mu) \ni g \mapsto \int h \ \chi(g)\, d\mu$ is differentiable for all $h \in L^q(E;\mu)$, $q \in (1, \infty]$, $p^{-1} + q^{-1} = 1$, the assertion follows immediately by the fundamental theorem of calculus. □

For $f \in \mathcal{D}_1$ (cf. (3.10)) define

$$f_t := e^{-H_m^X} f, \quad t \geq 0.$$

**Lemma 3.10** *Let $f \in \mathcal{D}_1$. Then $\mathbb{R}_+ \ni t \mapsto \langle \log(1 + f_t), \cdot \rangle \in L^2(\Gamma_X; \mu_{\lambda,m})$ is differentiable with derivative equal to*

$$-\langle (1 + f_t)^{-1} H_m^X f_t, \cdot \rangle, t \in \mathbb{R}_+.$$

**Proof.** Let $t \in \mathbb{R}_+$, $\varepsilon \in \mathbb{R}$ with $\varepsilon > 0$ if $t = 0$. Then by (3.3) (since $H_m^X f_t \in L^1(X;m) \cap L^2(X;m)$)

$$\int \left( \frac{1}{\varepsilon} (\langle \log(1 + f_{t+\varepsilon}), \gamma \rangle - \langle \log(1 + f_t), \gamma \rangle) + \langle (1 + f_t)^{-1} H_m^X f_t, \gamma \rangle \right)^2 \mu_{\lambda,m}(d\gamma)$$

$$= \int z \lambda(dz) \int A(\varepsilon, t)^2 \, dm + \int z^2 \lambda(dz) \left( \int A(\varepsilon, t) \, dm \right)^2,$$

where

$$A(\varepsilon, t) := \frac{1}{\varepsilon} (\log(1 + f_{t+\varepsilon}) - \log(1 + f_t)) + (1 + f_t)^{-1} H_m^X f_t \ .$$

Now the assertion follows by Remark 3.6 and Lemma 3.9 □

Note that clearly the measure

$$\mathcal{B}(\Gamma_X) \ni B \mapsto \int_{\mathbb{R}_+} \int_{\Gamma_X} 1_B(\gamma)\ \pi_{z\cdot m}(d\gamma)\ z\ \lambda(dz)$$

is absolutely continuous w.r.t. $\mu_{\lambda,m}$. Let $\rho_{\lambda,m}$ denote the corresponding Radon–Nikodym derivative. Since by Hölder's inequality

$$\int \rho_{\lambda,m}^2\, d\mu_{\lambda,m} = \iint \rho_{\lambda,m}(\gamma)\, z\, \pi_{z\cdot m}(d\gamma)\, \lambda(dz) \leq \left(\int \rho_{\lambda,m}^2\, d\mu_{\lambda,m}\right)^{1/2} \left(\int z^2 \lambda(dz)\right)^{1/2},$$

it follows that

$$\int \rho_{\lambda,m}^2\, d\mu_{\lambda,m} \leq \int z^2\, \lambda(dz)\ . \tag{3.12}$$

**Lemma 3.11** *Suppose condition (A) holds and let $f \in \mathcal{D}_1$. Then:*

*(i)* $\exp(\langle \log(1+f), \cdot\rangle) \in D(H^\Gamma_{\mu_{\lambda,m}})$ *and*

$$H^\Gamma_{\mu_{\lambda,m}} \exp(\langle \log(1+f), \cdot\rangle) = \left(\langle \frac{H^X_m f}{1+f}, \cdot\rangle - \rho_{\lambda,m} \int H^X_m f\, dm\right) \exp(\langle \log(1+f), \cdot\rangle)\ . \tag{3.13}$$

*(ii) For every $F \in E(\mathcal{D}_1)$ there exist $F_n \in E(\mathcal{D}_1 \cap C_0^\infty(X))$ ($\subset \mathcal{F}C_b^\infty$), $n \in \mathbb{N}$, such that $\sup_n \|F_n\|_\infty < \infty$ and $F_n \to F$ as $n \to \infty$ w.r.t. $\|\cdot\|_1$.*

**Proof.** Since $\exp(\langle \log(1+f), \cdot\rangle)$ and the right hand side of (3.13) are in $L^2(\Gamma_X; \mu_{\lambda,m})$, to show (i) it suffices to prove that $\exp(\langle \log(1+f), \cdot\rangle) \in D(H^\Gamma_{\mu_{\lambda,m},1})$ and that (3.13) holds with $H^\Gamma_{\mu_{\lambda,m},1}$ replacing $H^\Gamma_{\mu_{\lambda,m}}$.

By condition (A) there exist $\psi_n \in \mathcal{D}$, $n \in \mathbb{N}$, such that

$$\psi_n \to f,\ \ \Delta^X \psi_n \to H^X_m f \text{ in } L^2(X; m) \text{ as } n \to \infty. \tag{3.14}$$

Let $\delta, \delta' \in (0,1)$ such that $-1 < -\delta' < -\delta \leq f$, and let $\chi \in C_0^\infty(\mathbb{R})$ such that $-\delta' \leq \chi$ and $\chi(s) = s$ for all $s \in [-\delta, 1]$. Then for $n \in \mathbb{N}$ we have that $f_n := \chi \circ \psi_n \in \mathcal{D}$, hence $\log(1+f_n) \in \mathcal{D}$. Clearly, by (3.2)

$$\begin{aligned} &\int \langle (\log(1+f_n) - \log(1+f))^2, \gamma\rangle\, \mu_{\lambda,m}(d\gamma) \\ =\ &\int (\log(1+f_n) - \log(1+f))^2 dm \int z\ \lambda(dz) \underset{n\to\infty}{\longrightarrow} 0. \end{aligned}$$

Hence for any $\xi \in C^\infty(\mathbb{R})$ such that $\xi(s) = s$ for all $s \in ]-\infty, 0]$ and $\operatorname{supp}\xi \cap [0,\infty[$ is compact, it easily follows that

$$\exp\circ\xi(\langle \log(1+f_n),\cdot\rangle) \to \exp(\langle \log(1+f),\cdot\rangle) \tag{3.15}$$

in $L^1(\Gamma_X;\mu_{\lambda,m})$.

Furthermore, by (2.24) (since $H^\Gamma_{\mu_{\lambda,m},1}$ extends $H^\Gamma_{\mu_{\lambda,m}}$, hence $\Delta^\Gamma$)

$$\begin{aligned}
&H^\Gamma_{\mu_{\lambda,m},1}\exp\circ\xi(\langle \log(1+f_n),\cdot\rangle)\\
&= \Big(\langle \frac{\Delta^X f_n}{1+f_n},\cdot\rangle\cdot \xi'(\langle \log(1+f_n),\cdot\rangle)\\
&\quad +[\xi'(\langle \log(1+f_n),\cdot\rangle) - \xi''(\langle \log(1+f_n),\cdot\rangle)\\
&\quad -(\xi')^2(\langle \log(1+f_n),\cdot\rangle)]\cdot\langle \frac{|\nabla^X f_n|^2_{TX}}{(1+f_n)^2},\cdot\rangle\Big)\\
&\qquad \exp\circ\xi(\langle \log(1+f_n),\cdot\rangle).
\end{aligned} \tag{3.16}$$

We have that for all $n\in\mathbb{N}$

$$\Delta^X f_n = \chi'(\psi_n)\Delta^X\psi_n + \chi''(\psi_n)|\nabla^X\psi_n|^2_{TX}. \tag{3.17}$$

Since $\chi''(f)=0$, (3.14) implies that

$$\frac{\chi''(\psi_n)}{1+f_n}|\nabla^X\psi_n|^2_{TX} \to 0 \text{ in } L^1(X;m) \text{ as } n\to\infty.$$

Hence by (3.2)

$$\langle \frac{\chi''(\psi_n)}{1+f_n}|\nabla^X\psi_n|^2_{TX},\cdot\rangle \to 0 \text{ in } L^1(\Gamma_X;\mu_{\lambda,m}) \text{ as } n\to\infty. \tag{3.18}$$

Furthermore, by (3.12) and (3.3) resp. (2.18)

$$\begin{aligned}
&\int \Big(\langle \frac{H^X_m f}{1+f},\gamma\rangle - \rho_{\lambda,m}\cdot\int H^X_m f\,dm - \langle \frac{\chi'(\psi_n)}{1+f_n}\,\Delta^X\psi_n,\gamma\rangle\Big)^2\,\mu_{\lambda,m}(d\gamma)\\
=\;& \int \rho^2_{\lambda,m}\,d\mu_{\lambda,m}\Big(\int H^X_m f\,dm\Big)^2 + \int \langle \frac{H^X_m f}{1+f} - \frac{\chi'(\psi_n)}{1+f_n}\Delta^X\psi_n,\gamma\rangle^2\mu_{\lambda,m}(d\gamma)\\
&-2\int H^X_m f\,dm\int\int \langle \frac{H^X_m f}{1+f} - \frac{\chi'(\psi_n)}{1+f_n}\Delta^X\psi_n,\gamma\rangle\;\pi_{z\cdot m}(d\gamma)\,z\,\lambda(dz)
\end{aligned}$$

$$\leq \int z^2\lambda(dz)\left(\int H_m^X f\,dm\right)^2 + \int z\,\lambda(dz)\int (\frac{H_m^X f}{1+f} - \frac{\chi'(\psi_n)}{1+\chi(\psi_n)}\,\Delta^X\psi_n)^2\,dm$$

$$+\int z^2\lambda(dz)\left(\int(\frac{H_m^X f}{1+f} - \frac{\chi'(\psi_n)}{1+\chi(\psi_n)}\Delta^X\psi_n)\,dm\right)^2$$

$$-2\int z^2\lambda(dz)\ \int H_m^X f\,dm\int(\frac{H_m^X f}{1+f} - \frac{\chi'(\psi_n)}{1+\chi(\psi_n)}\Delta^X\psi_n)\,dm\ .$$

But

$$\begin{aligned}
&\int \frac{\chi'(\psi_n)}{1+\chi(\psi_n)}\Delta^X\psi_n\,dm \\
= &\int \frac{\chi''(\psi_n)(1+\chi(\psi_n)) - \chi'(\psi_n)^2}{(1+\chi(\psi_n))^2}|\nabla^X\psi_n|_{TX}^2\,dm \\
\underset{n\to\infty}{\longrightarrow} & -\int \frac{|\nabla^X f|_{TX}^2}{(1+f)^2}\,dm = -\int \langle \nabla^X(\frac{f}{1+f}), \nabla^X f\rangle_{TX}\,dm \\
= &\int \frac{H_m^X f}{1+f}\,dm - \int H_m^X f\,dm
\end{aligned}$$

(because $\chi(f) = f$, $\chi'(f) = 1$). Therefore,

$$\langle \frac{\chi'(\psi_n)}{1+f_n}\Delta^X\psi_n, \cdot\rangle \to \langle \frac{H_m^X f}{1+f}, \cdot\rangle - \rho_{\lambda,m}\cdot\int H_m^X f\,dm \text{ in } L^2(\Gamma_X;\mu_{\lambda,m}) \text{ as } n\to\infty. \tag{3.19}$$

(3.17) – (3.19) and the fact that $\xi' \equiv 1$ on $]-\infty,0]$ imply that

$$\begin{aligned}
&\langle \frac{\Delta^X f_n}{1+f_n}, \cdot\rangle\,\xi'(\langle \log(1+f_n), \cdot\rangle)\exp\circ\xi(\langle\log(1+f_n),\cdot\rangle) \\
&\underset{n\to\infty}{\longrightarrow} \left(\langle \frac{H_m^X f}{1+f}, \cdot\rangle - \rho_{\lambda,m}\cdot\int H_m^X f\,dm\right)\ \exp(\langle\log(1+f),\cdot\rangle) \\
&\quad \text{in } L^1(\Gamma_X;\mu_{\lambda,m}).
\end{aligned}$$

Since by (3.14)

$$\frac{|\nabla^X f_n|_{TX}^2}{(1+f_n)^2} = \frac{\chi'(\psi_n)^2\,|\nabla^X\psi_n|_{TX}^2}{(1+f_n)^2} \to \frac{|\nabla^X f|_{TX}^2}{(1+f)^2} \quad \text{in } L^1(X;m) \text{ as } n\to\infty,$$

it follows by (3.3), because $\xi' = (\xi')^2 \equiv 1$ on $]-\infty,0]$ and hence $\xi'' \equiv 0$ on $]-\infty,0]$, that the remaining terms on the right hand side of (3.16) converge to zero in $L^1(\Gamma_X;\mu_{\lambda,m})$. Thus (i) is proved.

Assertion (ii) also follows then easily since
$\exp \circ \xi(\langle \log(1+f_n), \cdot\rangle) \in \mathcal{F}C_b^\infty$ and $|F_n| \le \max(1, \exp \|\xi\|_\infty)$ for all $n \in \mathbb{N}$.
□

**Proof of Proposition 3.7.** (i): By a simple approximation argument it suffices to prove (i), if $f \in \mathcal{D}_1$. Then by Lemmas 3.9, 3.10

$$\mathbb{R}_+ \ni t \mapsto F(t) := \exp(\langle \log(1+f_t), \cdot\rangle - z \int (f_t - f)\, dm) \in L^2(\Gamma_X; \mu_{\lambda,m})$$

is differentiable with derivative equal to

$$\begin{aligned} &\Big(-\langle \frac{H_m^X f_t}{1+f_t}, \cdot\rangle + z \int H_m^X f_t \, dm\Big) \\ &\qquad \cdot \exp(\langle \log(1+f_t), \cdot\rangle - z \int (f_t - f)\, dm), \\ &\qquad t \in \mathbb{R}_+ \ . \end{aligned} \tag{3.20}$$

Realizing that $\rho_{\lambda,m} \equiv z$ if $\lambda = \varepsilon_z$, we see that by Lemma 3.11 this coincides with

$$-H^\Gamma_{\mu_{\lambda,m}} \exp(\langle \log(1+f_t), \cdot\rangle - z \int (f_t - f)\, dm)$$

i.e.,

$$\frac{d}{dt} F(t) = -H^\Gamma_{\mu_{\lambda,m}} F(t), \quad F(0) = \exp(\langle \log(1+f), \cdot\rangle). \tag{3.21}$$

$\int H_m^X f_t \, dm = 0$ for all $t \in \mathbb{R}_+$, if $(H_m^X, D(H_m^X))$ is conservative. Hence (3.21) holds for general probability measures $\lambda$ on $(\mathbb{R}_+, \mathcal{B}(\mathbb{R}_+))$ satisfying (3.1) in this case. But also

$$\mathbb{R}_+ \ni t \mapsto T^\Gamma_{\mu_{\lambda,m}}(t) \exp(\langle \log(1+f), \cdot\rangle) \in L^2(\Gamma_X; \mu_{\lambda,m})$$

solves the Cauchy problem (3.21). But as is well-known (and can be easily proved using Duhamel's formula) the solution to (3.21) is unique. Hence assertion (i) is proved.

(ii): By part (i) and by Remark 3.6 (ii)

$$T^\Gamma_{\mu_{\lambda,m}}(t)\, (E(\mathcal{D}_1)) \subset E(\mathcal{D}_1) \ . \tag{3.22}$$

Since $E(\mathcal{D}_1)$ is obviously dense in $L^2(\Gamma_X; \mu_{\lambda,m})$, hence in $L^1(\Gamma_X; \mu_{\lambda,m})$ and since $E(\mathcal{D}_1) \subset D(H^\Gamma_{\mu_{\lambda,m}})$ by Lemma 3.11 (i), the assertion follows from a standard theorem on operator semigroups (cf. e.g. [RS75, Theorem X.49]).

(iii): Let $F \in E(\mathcal{D}_1)$ and $F_n$, $n \in \mathbb{N}$, as in Lemma 3.11(ii). Then $F_n \to F$ in $L^2(\Gamma_X; \mu_{\lambda,m})$ as $n \to \infty$ and

$$\int |\nabla^\Gamma F_n|^2_{T\Gamma_X} d\mu_{\lambda,m} = -\int F_n \ \Delta^\Gamma F_n \ d\mu_{\lambda,m} \underset{n\to\infty}{\longrightarrow} \int F \ H^\Gamma_{\mu_{\lambda,m}} F \, d\mu_{\lambda,m} \ .$$

In particular,

$$\sup_n \mathcal{E}^\Gamma_{\mu_{\lambda,m}}(F_n, F_n) < \infty.$$

Now by [MR92, Chap. I, Lemma 2.12] there exists a subsequence $(n_k)_{k\in\mathbb{N}}$ such that for $G_N := \frac{1}{N}\sum_{k=1}^N F_{n_k} \in E(\mathcal{D}_1 \cap C_0^\infty(X))$, $N \in \mathbb{N}$,

$$G_N \longrightarrow F \text{ as } N \to \infty \text{ w.r.t. } \| \cdot \|_{H_0^{1,2}(\Gamma_X;\mu_{\lambda,m})} \ .$$

□

## 3.3 Brownian motion

In this section we shall treat the diffusion process associated with our Dirichlet form $(\mathcal{E}^\Gamma_{\mu_{\lambda,m}}, D(\mathcal{E}^\Gamma_{\mu_{\lambda,m}}))$ from Subsection 3.1 and identify it (when started with $\mu_{\lambda,m}$) as the well–known independent infinite particle process on $X$ provided condition (A) in Theorem 3.3 is fulfilled and $(H^X_m, D(H^X_m))$ is conservative.

The desired process in general will live on the bigger state space $\ddot{\Gamma}$ consisting of all $\mathbb{Z}_+$-valued Radon measures on $X$ (which is Polish, see e.g. [MR97]). Since $\Gamma_X \subset \ddot{\Gamma}_X$ and $\mathcal{B}(\ddot{\Gamma}) \cap \Gamma_X = \mathcal{B}(\Gamma_X)$, we can consider $\mu_{\lambda,m}$ as a measure on $(\ddot{\Gamma}, \mathcal{B}(\ddot{\Gamma}))$ and correspondingly $(\mathcal{E}, D(\mathcal{E}))$ as a Dirichlet form on $L^2(\ddot{\Gamma}; \mu_{\lambda,m})$.

**Theorem 3.12** $(\mathcal{E}^\Gamma_{\mu_{\lambda,m}}, D(\mathcal{E}^\Gamma_{\mu_{\lambda,m}}))$ *is quasi–regular. In particular, there exists a conservative diffusion process*

$$\mathbf{M} = (\mathbf{\Omega}, \mathbf{F}, (\mathbf{F})_{t\geq 0}, (\mathbf{\Theta}_t)_{t\geq 0}, (\mathbf{X}_t)_{t\geq 0}, (\mathbf{P}_\gamma)_{\gamma\in\ddot{\Gamma}})$$

*on* $\ddot{\Gamma}$ *(cf. [Dy65]) which is properly associated with* $(\mathcal{E}^\Gamma_{\mu_{\lambda,m}}, D(\mathcal{E}^\Gamma_{\mu_{\lambda,m}}))$*, i.e., for all (* $\mu_{\lambda,m}$*-versions of)* $F \in L^2(\ddot{\Gamma}; \mu_{\lambda,m})$ *and all* $t > 0$ *the function*

$$\gamma \ \mapsto \ p_t F(\gamma) \ := \ \int_\Omega F(X_t) d\mathbf{P}_\gamma, \ \ \gamma \in \ddot{\Gamma}, \tag{3.23}$$

*is an* $\mathcal{E}^\Gamma_{\mu_{\lambda,m}}$*-quasi-continuous version of* $\exp(-tH^\Gamma_{\mu_{\lambda,m}})F$*.* $\mathbf{M}$ *is up to* $\mu_{\lambda,m}$*-equivalence unique (cf. [MR92, Ch. IV, Sect. 6]). In particular,* $\mathbf{M}$ *is* $\mu_{\lambda,m}$*–symmetric (i.e.,* $\int G \ p_t F d\mu_{\lambda,m} = \int F \ p_t G \ d\mu_{\lambda,m}$ *for all* $F, G : \ddot{\Gamma} \to \mathbb{R}_+, \mathcal{B}(\ddot{\Gamma})$*-measurable) and has* $\mu_{\lambda,m}$ *as an invariant measure.*

In the above theorem $\mathbf{M}$ is canonical, i.e., $\mathbf{\Omega} = C([0,\infty[\to \ddot{\Gamma})$, $\mathbf{X}_t(\omega) := \omega(t)$, $t \geq 0$, $\omega \in \mathbf{\Omega}$, $(\mathbf{F}_t)_{t\geq 0}$ together with $\mathbf{F}$ is the corresponding minimum completed admissible family (cf. [F80, Sect.4.1] and also [FOT94]) and $\mathbf{\Theta}_t$, $t \geq 0$, are the corresponding natural time shifts.

Theorem 3.12 is in fact a special case of Theorems 4.12, 4.14 in the next section where we shall also recall the definition of quasi–regularity, $\mathcal{E}^{\Gamma}_{\mu_{\lambda,m}}$–quasi–continuity etc. So, we refer to Subsection 4.2 for more details.

Let $\mathbf{M} = (\mathbf{\Omega}, \mathbf{F}, (\mathbf{F}_t)_{t\geq 0}, (\mathbf{\Theta}_t)_{t\geq 0}, (\mathbf{X}_t)_{t\geq 0}, (\mathbf{P}_\gamma)_{\gamma\in\ddot{\Gamma}})$ be as in Theorem 3.12. As usual we define

$$\mathbf{P}_{\mu_{\lambda,m}} := \int \mathbf{P}_\gamma \, \mu_{\lambda,m}(d\gamma). \tag{3.24}$$

Let $M = (\Omega, \mathcal{F}, (\mathcal{F}_t)_{t\geq 0}, (X_t)_{t\geq 0}, (P_x)_{x\in X})$ be the diffusion process on $X$ associated with the Dirichlet form $(\mathcal{E}^X_m, D(\mathcal{E}^X_m))$ (cf. e.g. [MR92, Ch.IV, Subsect.4 a) resp. Ch. II, Subsect. 2 a)]). Set

$$P_m := \int_X P_x \, m(dx). \tag{3.25}$$

We note that in the non–conservative case $\Omega = C([0,\infty[\to X_\Delta)$ rather than $C([0,\infty[\to X)$. Here $X_\Delta$ is the one–point compactification of $X$.

In order to identify $\mathbf{M}$ as the independent infinite particle process if $(H^X_m, D(H^X_m))$ is conservative, we have to show that $\mathbf{P}_{\mu_{\lambda,m}}$ has the same finite dimensional distributions as $\int_{\mathbb{R}_+} \pi_{zP_m} \, \lambda(dz)$, where $\pi_{zP_m}$ is the Poisson measure on $\mathbf{\Omega}$ with intensity $zP_m$, $z \geq 0$ (cf. e.g. [ST74], [M-L76]). For this it suffices to prove the following proposition whose proof is straight-forward by Proposition 3.7(i) (see also [ST74, Proposition 1.5]).

**Proposition 3.13** *Assume that condition* (A) *in Theorem 3.3 holds and that* $(H^X_m, D(H^X_m))$ *is conservative. Then for all* $f_1, \ldots, f_N \in L^1(X;m)$ *taking values in* $]-1,0]$ *and* $0 \leq t_1 < \cdots < t_N < \infty$

$$\begin{aligned}&\int_{\mathbf{\Omega}} \exp(\langle \log(1+f_1), \mathbf{X}_{t_1}\rangle) \cdots \exp(\langle \log(1+f_N), \mathbf{X}_{t_N}\rangle) d\mathbf{P}_{\mu_{\lambda,m}} \\ &= \int \exp\left( z \int_\Omega [(1+f_1(X_{t_1})) \cdots (1+f_N(X_{t_N})) - 1] dP_m \right) \lambda(dz). \end{aligned} \tag{3.26}$$

**Proof.** If $N = 1$, then by Proposition 3.7(i)

$$\int_{\mathbf{\Omega}} \exp(\langle \log(1+f_1), \mathbf{X}_{t_1}\rangle) d\mathbf{P}_{\mu_{\lambda,m}} = \int_{\Gamma_X} \exp(\langle \log(1+e^{-t_1 H_m^X} f_1), \gamma\rangle)\mu_{\lambda,m}(d\gamma)$$

$$= \int \exp\Big[z\int_X e^{-t_1 H_m^X} f_1\, dm\Big]\lambda(dz) = \int \exp\Big[z\int_\Omega f_1(X_{t_1}) dP_m\Big]\lambda(dz).$$

Suppose (3.26) holds for $N-1$. Then by Proposition 3.7 and the Markov property of $\mathbf{P}_{\mu_{\lambda,m}}$ the left hand side of (3.26) is equal to

$$\int_{\mathbf{\Omega}} \exp(\langle \log(1+f_1), \mathbf{X}_{t_1}\rangle)\cdots \exp(\langle \log(1+f_{N-1}), \mathbf{X}_{t_{N-1}}\rangle)$$
$$\times \exp(\langle \log(1+e^{-(t_N-t_{N-1})H_m^X} f_N), \mathbf{X}_{t_{N-1}}\rangle) d\mathbf{P}_{\mu_{\lambda,m}}.$$

This in turn by induction hypothesis is equal to

$$\int \exp\Big(z\int_\Omega [(1+f_1(X_{t_1}))\cdots(1+f_{N-1}(X_{t_{N-1}}))$$
$$\times (1+e^{-(t_N-t_{N-1})H_m^X} f_N(X_{t_{N-1}})) - 1] dP_m\Big)\lambda(dz)$$

which by the Markov property of $P_m$ is just the right hand side of (3.26). $\square$

In the non–conservative case the finite dimensional distributions of $P_{\mu_{\lambda,m}}$ and $\int_{\mathbb{R}_+} \pi_{zP_m}\, \lambda(dz)$ do not coincide in general as the following result shows:

**Proposition 3.14** *Assume condition (A) in Theorem 3.3 holds and let* $z \in [0,\infty[$*. Then for all* $f_1,\ldots,f_N \in L^1(X;m)$ *taking values in* $]-1,0]$ *and* $0 \le t_1 < \ldots < t_N < \infty$

$$\begin{aligned}
&\int_{\mathbf{\Omega}} \exp(\langle \log(1+f_1), \mathbf{X}_{t_1}\rangle)\ldots\exp(\langle \log(1+f_N), \mathbf{X}_{t_N}\rangle)\; d\mathbf{P}_{\pi_{z\cdot m}}\\
= \;& \exp(z\int_\Omega [(1+f_1(X_{t_1}))\ldots(1+f_N(X_{t_N})) - 1]\, dP_m\\
& \cdot \exp\Big(-z\int (e^{-t_1H_m^X}-1) f_1 dm - z\sum_{i=2}^N \int (e^{-t_1H_m^X}-1)\Big[(1+f_1)e^{-(t_2-t_1)H_m^X}\\
& (1+f_2)e^{-(t_3-t_2)H_m^X}\ldots(1+f_{i-1})e^{-(t_i-t_{i-1})H_m^X} f_i\Big]\, dm\Big)\;. \qquad (3.27)
\end{aligned}$$

By Proposition 3.7 (i) the proof of Proposition 3.14 is entirely the same as that of 3.13.

**Remark 3.15** (i) Since for $a_1, \dots, a_N \in \mathbb{R}$,

$$(1+a_1)\dots(1+a_N) - 1 = a_1 + \sum_{i=2}^{N}(1+a_1)\dots(1+a_{i-1})a_i$$

the right hand side of (3.27) obviously simplifies to

$$\exp\Big(z\int f_1\,dm + \sum_{i=1}^{N}\int(1+f_1)e^{-(t_2-t_1)H_m^X}\Big[(1+f_2)e^{-(t_3-t_2)H_m^X}[\dots(1+f_{i-1})e^{-(t_i-t_{i-1})H_m^X}f_i]\Big]\,dm\Big).$$

(ii) Since $\mathbf{M}$ is always conservative, $M$, however, if and only if so is $(H_m^X, D(H_m^X))$, Proposition 3.14 is, of course, not surprising.

(iii) ("$\sigma$–case") Let $\sigma$ be as in the beginning of Remark 2.3(i). As mentioned before all the above extends to the case where $\sigma$ replaces $m$. The diffusion process $M$ associated to $(\mathcal{E}_\sigma^X, D(\mathcal{E}_\sigma^X))$ is usually called *distorted Brownian motion on* $X$. Therefore, we call the diffusion associated to $(\mathcal{E}_{\mu_{\lambda,\sigma}}^\Gamma, D(\mathcal{E}_{\mu_{\lambda,\sigma}}^\Gamma))$ *distorted Brownian motion on* $\Gamma_X$ because of the complete analogy. Likewise, since $M$ is just ordinary Brownian motion on $X$ if $\sigma = m$, it is justified to call $\mathbf{M}$ *Brownian motion* on $\Gamma_X$ in this case. Therefore, Proposition 3.13 shows that, provided $M$ or equivalently $(H_\sigma^X, D(H_\sigma^X))$ is conservative and condition (A) holds, the usual independent particle process known for more than 40 years (cf. [D53], [Do56]) is nothing but (distorted) Brownian motion on $\Gamma_X$.

We conclude this section with a result concerning the question whether the bigger state space $\ddot{\Gamma}_X$ rather than just $\Gamma_X$ is really needed.

**Theorem 3.16** *Suppose $X = \mathbb{R}^d$, $d \geq 2$. Let $\sigma = \rho \cdot m$ be as specified at the beginning of Remark 2.3 (i) and assume, in addition, that $\sigma \in L^2_{loc}(\mathbb{R}^d; m)$. Then $\ddot{\Gamma}_X \setminus \Gamma_X$ is $\mathcal{E}_{\mu_{\lambda,\sigma}}^\Gamma$–exceptional.*

This result was obtained in [RS97] to which we refer for its proof. Theorem 3.16 just says that for $\mathcal{E}_{\mu_{\lambda,\sigma}}$–q.e. $\gamma \in \Gamma_X$ under $\mathbf{P}_\gamma$, $(\mathbf{X}_t)_{t\geq 0}$ never leaves $\Gamma_X$ (cf. [MR92, Proposition 5.30 (i)]).

## 3.4 Ergodicity in the free case

Let $\mathbf{M} = (\mathbf{\Omega}, \mathbf{F}, (\mathbf{F}_t)_{t\geq 0}, (\mathbf{\Theta}_t)_{t\geq 0}, (\mathbf{X}_t)_{t\geq 0}, (\mathbf{P}_\gamma)_{\gamma\in\ddot{\Gamma}_X})$ be as in Theorem 3.16. Concerning ergodicity / irreducibility we recall the following well–known result:

**Proposition 3.17** *The following assertions:*

(i) $\mathbf{P}_{\mu_{\lambda,m}}$ *is (time–)ergodic (i.e., every bounded* $\mathbf{F}$*–measurable function* $G : \mathbf{\Omega} \to \mathbb{R}$ *which is* $\mathbf{\Theta}_t$*–invariant for all* $t \geq 0$*, is constant* $\mathbf{P}_{\mu_{\lambda,m}}$*–a.e.).*

(ii) $(\mathcal{E}^\Gamma_{\mu_{\lambda,m}}, D(\mathcal{E}^\Gamma_{\mu_{\lambda,m}}))$ *is irreducible (i.e., for* $F \in D(\mathcal{E}^\Gamma_{\mu_{\lambda,m}})$*,* $\mathcal{E}^\Gamma_{\mu_{\lambda,m}}(F,F) = 0$ *implies that* $F = \mathrm{const}$*).*

(iii) $(T^\Gamma_{\mu_{\lambda,m}}(t))_{t>0}$ *is irreducible (i.e., if* $G \in L^2(\Gamma_X;\mu_{\lambda,m})$ *such that* $T^\Gamma_{\mu_{\lambda,m}}(t)(GF) = GT^\Gamma_{\mu_{\lambda,m}}(t)F$ *for all* $F \in L^\infty(\Gamma_X;\mu_{\lambda,m}), t > 0$*, then* $G = \mathrm{const}$*).*

(iv) *If* $F \in L^2(\Gamma_X;\mu_{\lambda,m})$ *such that* $T^\Gamma_{\mu_{\lambda,m}}(t)F = F$ *for all* $t > 0$*, then* $F = \mathrm{const}$*.*

(v) $(T^\Gamma_{\mu_{\lambda,m}}(t))_{t>0}$ *is ergodic (i.e.,*

$$\int (T^\Gamma_{\mu_{\lambda,m}}(t)F - \int F d\mu_{\lambda,m})^2 d\mu_{\lambda,m} \to 0 \ \ as \ \ t\to\infty \ \ for \ \ all \ \ F \in L^2(\Gamma_X;\mu_{\lambda,m})).$$

(vi) *If* $F \in D(H^\Gamma_{\mu_{\lambda,m}})$ *with* $H^\Gamma_{\mu_{\lambda,m}}F = 0$*, then* $F = \mathrm{const}$*.*

**Proof.** The equivalence of (ii) – (vi) can be shown entirely analogously to [AKR97d, Proposition 2.3]. The equivalence of (i) and (ii) is a direct consequence of the regularization method in [MR92, Chap. VI, Sect. 1] and the main results in [F80], [F83]. □

The main result here is, however, the following:

**Theorem 3.18** *Assume that* $m(X) = \infty$*. If condition (A) in Theorem 3.3 holds and* $(H^X_m, D(H^X_m))$ *is conservative, then any of the assertions (i)–(vi) of Proposition 3.17 is equivalent to:*

(vii) $\mu_{\lambda,m} = \pi_{z\cdot m}$ *for some* $z \in [0,\infty[$*.*

This theorem is a special case of Theorem 9.9 below. So, we refer to the detailed discussion in Section 9.

**Remark 3.19** *("$\sigma$–case")* There are generalizations of Theorem 3.18 and Proposition 3.17 in the case where $m$ is replaced by $\sigma = \rho\cdot m$ with $\sigma(X) = \infty$ for a large class of even not weakly differentiable functions $\rho : X \to \mathbb{R}_+$ (cf. [AKRR98]). However, if $\rho$ satisfies all assumptions specified in the beginning of Remark 2.3 (i), then the generalization to the "$\sigma$–case is as before straightforward (cf. [AKR97a].)

# 4 Classical Dirichlet forms on configuration spaces w.r.t. general measures

Many results on Dirichlet forms, defined in terms of the gradient $\nabla^\Gamma$ and the "tangent bundle" $(T_\gamma\Gamma_X)_{\gamma\in\Gamma_X} = (L^2(X \to TX;\gamma))_{\gamma\in\Gamma_X}$ introduced before, can be proved for quite arbitrary probability measures $\mu$ on $(\Gamma_X, \mathcal{B}(\Gamma_X))$ replacing the mixed Poisson measure $\mu_{\lambda,m}$. The purpose of this section is to summarize some of these results which are of particular importance. In this section we fix a probability measure $\mu$ on $(\Gamma_X, \mathcal{B}(\Gamma_X))$ with Radon mean, i.e.,

$$\int \gamma(K)\ \mu(d\gamma) < \infty \text{ for all compact } K \subset X. \tag{4.1}$$

We define

$$\mathcal{E}^\Gamma_\mu(F,G) := \int_{\Gamma_X} \langle \nabla^\Gamma F, \nabla^\Gamma G\rangle_{T\Gamma_X}\ d\mu\ ; \quad F,G \in \mathcal{F}C_b^\infty. \tag{4.2}$$

Note that by (2.8), the Cauchy–Schwarz inequality and (4.1) the integral in (4.2) exists for all $F,G \in \mathcal{F}C_b^\infty$, hence $\mathcal{E}^\Gamma_\mu$ is well-defined. Furthermore, assume that

$$\nabla^\Gamma F = \nabla^\Gamma G\ \mu\text{–a.e. if } F,G \in \mathcal{F}C_b^\infty \text{ such that } F = G\ \mu\text{-a.e.} \tag{4.3}$$

Clearly, (4.3) holds if supp $\mu = \Gamma_X$. Let $\mathcal{F}C_b^{\infty,\mu}$ denote the set of $\mu$–equivalence classes determined by $\mathcal{F}C_b^\infty$. (4.3) implies that $(\mathcal{E}^\Gamma_\mu, \mathcal{F}C_b^{\infty,\mu})$ is a well–defined symmetric positive defined bilinear form on $L^2(\Gamma_X;\mu)$. It is also densely defined, since by an easy monotone class argument its domain $\mathcal{F}C_b^{\infty,\mu}$ is dense in $L^2(\Gamma_X;\mu)$.

## 4.1 Closability, domains and Sobolev spaces

We first recall the definition of closability: Equip $\mathcal{F}C_b^{\infty,\mu}$ with the norm (cf. (3.6))

$$\|\cdot\|_{H_0^{1,2}(\Gamma_X;\mu)} := \left(\mathcal{E}_\mu^\Gamma(\cdot,\cdot) + (\ ,\ )_{L^2(\Gamma_X;\mu)}\right)^{1/2}. \tag{4.4}$$

Then the embedding map

$$i : (\mathcal{F}C_b^{\infty,\mu}, \|\cdot\|_{H_0^{1,2}(\Gamma_X;\mu)}) \to (L^2(\Gamma_X;\mu), \|\cdot\|_{L^2(\Gamma_X;\mu)})$$

is continuous, hence $i$ uniquely extends to a continuous (linear) map $\bar{i}$ from the abstract completion $\overline{\mathcal{F}C_b^{\infty,\mu}}$ of $\mathcal{F}C_b^{\infty,\mu}$ w.r.t. $\|\cdot\|_{H_0^{1,2}(\Gamma_X;\mu)}$ to $L^2(\Gamma_X;\mu)$.

**Definition 4.1** $(\mathcal{E}_\mu^\Gamma, \mathcal{F}C_b^{\infty,\mu})$ *is called closable on* $L^2(\Gamma_X;\mu)$ *if the above map*

$$\bar{i} := \overline{\mathcal{F}C_b^{\infty,\mu}} \longrightarrow L^2(\Gamma_X;\mu)$$

*is one–to–one.*

**Remark 4.2** (i) Clearly, $(\mathcal{E}_\mu^\Gamma, \mathcal{F}C_b^{\infty,\mu})$ is by definition closable on $L^2(\Gamma_X;\mu)$ if for all $F_n \in \mathcal{F}C_b^{\infty,\mu}$, $n \in \mathbb{N}$, such that

$$\mathcal{E}_\mu^\Gamma(F_n - F_m, F_n - F_m) \underset{n,m\to\infty}{\longrightarrow} 0, \ \|F_n\|_{L^2(\Gamma_X;\mu)} \underset{n\to\infty}{\longrightarrow} 0,$$

it follows that $\mathcal{E}_\mu^\Gamma(F_n, F_n) \underset{n\to\infty}{\longrightarrow} 0$.

(ii) If $(\mathcal{E}_\mu^\Gamma, \mathcal{F}C_b^{\infty,\mu})$ is closable on $L^2(\Gamma_X;\mu)$, $\mathcal{E}_\mu^\Gamma$ uniquely extends to a continuous bilinear map on all of $\overline{\mathcal{F}C_b^{\infty,\mu}} =: D(\mathcal{E}_\mu^\Gamma)$. $(\mathcal{E}_\mu^\Gamma, D(\mathcal{E}_\mu^\Gamma))$ is the smallest closed extension of $(\mathcal{E}_\mu^\Gamma, \mathcal{F}C_b^{\infty,\mu})$.

In general it is not true that $(\mathcal{E}_\mu^\Gamma, \mathcal{F}C_b^{\infty,\mu})$ is closable on $L^2(\Gamma_X;\mu)$, but, in fact, quite weak assumptions ensure this. We refer to Subsection 6.3 below and to [MR97], [AKRR98], [dSKR98] for other general methods and other concrete examples. In the rest of this subsection we shall mainly concentrate on measures satisfying an integration by parts formula of type (2.14), such as the mixed Poisson measures. Therefore, we introduce the following condition:

(IbP) $\int \gamma(K)^2 \ \mu(d\gamma) < \infty$ for all compact $K \subset X$ and for every $v \in V_0(X)$ there exists $B_v^\mu \in L^2(\Gamma_X;\mu)$ such that

$$\int \nabla_v^\Gamma F \, d\mu = -\int F \ B_v^\mu \, d\mu.$$

Since $\nabla_v^\Gamma$ satisfies the product rule, (IbP) implies that

$$\int \nabla_v^\Gamma F\ G\ d\mu \ = \ -\int F\ \nabla_v^\Gamma G\ d\mu - \int F\ G B_v^\mu\ d\mu$$
$$\text{for all } F, G \in \mathcal{F}C_b^\infty \text{ and all } v \in V_0(X). \qquad (4.5)$$

Furthermore, (IbP) , obviously, implies (4.3).

¿From now on we assume that (IbP) holds. The mixed Poisson measures of Section 3 and Ruelle measure to be introduced in Subsection 7.1 below are typical examples.

In analogy with (2.13), (2.24) for $V = \sum_{i=1}^N F_i\, v_i \in \mathcal{V}\mathcal{F}C_b^\infty$ we define

$$\operatorname{div}_\mu^\Gamma V := \sum_{i=1}^N (\nabla_{v_i}^\Gamma F_i + B_{v_i}^\mu F_i) \qquad (4.6)$$

and for $F \in \mathcal{F}C_b^\infty$

$$\Delta_\mu F := \operatorname{div}_\mu^\Gamma \nabla^\Gamma F\ . \qquad (4.7)$$

The following result is the analogue of Remark 2.3 (ii) and Propositions 3.1, 3.2. By (4.5)–(4.7) its proof is the same , and therefore, omitted.

**Proposition 4.3** *Suppose $\mu$ satisfies* (IbP) *. Then:*
*(i) The right hand side of* (4.6) *is $\mu$–a.e. independent of the representation of $V$.*
*(ii) For all $F, G \in \mathcal{F}C_b^\infty$*

$$\mathcal{E}_\mu^\Gamma(F, G) = -\int \Delta_\mu^\Gamma F\ G\, d\mu\ . \qquad (4.8)$$

*In particular, $(\mathcal{E}_\mu^\Gamma, \mathcal{F}C_b^{\infty,\mu})$ is closable on $L^2(\Gamma_X;\mu)$. Its closure is a symmetric Dirichlet form on $L^2(\Gamma_X;\mu)$ which is* conservative *(i.e., $1 \in D(\mathcal{E}_\mu^\Gamma)$ and $\mathcal{E}_\mu^\Gamma(1,1) = 0$). Its generator, denoted by $(H_\mu^\Gamma, D(H_\mu^\Gamma))$ is the Friedrichs' extension of $(-\Delta_\mu, \mathcal{F}C_b^{\infty,\mu})$ on $L^2(\Gamma_X;\mu)$.*

Clearly, $\nabla^\Gamma$ extends to all of $D(\mathcal{E}_\mu^\Gamma)$. We shall denote its extension by $\nabla_\mu^\Gamma$. As for $\mu = \mu_{\lambda,m}$ we set

$$H_0^{1,2}(\Gamma_X;\mu) := D(\mathcal{E}_\mu^\Gamma) \qquad (4.9)$$

i.e., the $(1,2)$–*Sobolev space on* $\Gamma_X$ (*w.r.t.* $\mu$).

We shall now construct an extension of $(\mathcal{E}_\mu^\Gamma, H_0^{1,2}(\Gamma_X;\mu))$ which is the exact analogue on $\Gamma_X$ of a weak Sobolev space on a finite dimensional manifold. This is essentially a special case of a construction by A. Eberle in [Eb97].

We first note that, as is easy to check, the $\mu$–equivalence classes $\mathcal{VFC}_b^{\infty,\mu}$ determined by $\mathcal{VFC}_b^\infty$ are dense in $L^2(\Gamma_X \to T\Gamma_X;\mu)$, i.e., the Hilbert space of $\mu$–square integrable sections in $T\Gamma_X$. Let $((\mathrm{div}_\mu^\Gamma)^{*,\mu}, D((\mathrm{div}_\mu^\Gamma)^{*,\mu}))$ be the adjoint of $(\mathrm{div}_\mu^\Gamma, \mathcal{VFC}_b^{\infty,\mu})$ as an operator from $L^2(\Gamma_X \to T\Gamma_X;\mu)$ to $L^2(\Gamma_X;\mu)$.

**Remark 4.4** As an adjoint the operator $((\mathrm{div}_\mu^\Gamma)^{*,\mu}, D((\mathrm{div}_\mu^\Gamma)^{*,\mu}))$ is automatically closed, and by definition it is an operator from $L^2(\Gamma_X;\mu)$ to $L^2(\Gamma \to T\Gamma_X;\mu)$. Furthermore, again by definition, $G \in L^2(\Gamma_X;\mu)$ belongs to $D((\mathrm{div}_\mu^\Gamma)^{*,\mu})$ if and only if there exists $V_G \in L^2(\Gamma_X \to T\Gamma_X;\mu)$ such that

$$\int G \ \mathrm{div}_\mu^\Gamma \ V \ d\mu = -\int \langle V_G, V\rangle_{T\Gamma_X} \ d\mu \text{ for all } V \in \mathcal{VFC}_b^\infty. \tag{4.10}$$

In this case $(\mathrm{div}_\mu^\Gamma)^{*,\mu} G = V_G$.

Because of Remark 4.4 we set

$$W^{1,2}(\Gamma_X;\mu) := D((\mathrm{div}_\mu^\Gamma)^{*,\mu}), \quad d^\mu := (\mathrm{div}_\mu^\Gamma)^{*,\mu} \tag{4.11}$$

and think of $W^{1,2}(\Gamma_X;\mu)$ as a *weak* $(1,2)$*–Sobolev space on* $\Gamma_X$ with norm

$$W^{1,2}(\Gamma_X;\mu) \ni G \mapsto \|G\|_{W^{1,2}(\Gamma_X;\mu)} := \left(\int \langle d^\mu G, d^\mu G\rangle_{T\Gamma_X} d\mu + \int G^2 \, d\mu\right)^{1/2}. \tag{4.12}$$

By (4.5) it follows from Remark 4.4 that

$$\mathcal{F}C_b^{\infty,\mu} \subset W^{1,2}(\Gamma_X;\mu) \text{ and } d^\mu = \nabla^\Gamma \text{ on } \mathcal{F}C_b^{\infty,\mu}, \tag{4.13}$$

hence the densely defined positive definite symmetric bilinear form

$$(F,G) \mapsto \int \langle d^\mu F, d^\mu G\rangle_{T\Gamma_X} d\mu,$$

with domain $W^{1,2}(\Gamma_X;\mu)$ extends $(\mathcal{E}_\mu^\Gamma, H_0^{1,2}(\Gamma_X;\mu))$ and is, therefore, denoted by $(\mathcal{E}_\mu^\Gamma, W^{1,2}(\Gamma_X;\mu))$.

**Proposition 4.5** *Let* $F \in H_0^{1,2}(\Gamma_X;\mu) \cap L^\infty(\Gamma_X;\mu)$ *and* $G \in W^{1,2}(\Gamma_X;\mu) \cap L^\infty(\Gamma_X;\mu)$. *Then* $FG \in W^{1,2}(\Gamma_X;\mu)$ *and*

$$d^\mu(FG) = F \ d^\mu G + G \ d^\mu F. \tag{4.14}$$

**Proof.** Let $V \in \mathcal{VFC}_b^\infty$ and assume first that $F \in \mathcal{FC}_b^\infty$. Then by the product rule for $\nabla^\Gamma$ and (4.6) for all $G \in W^{1,2}(\Gamma_X;\mu) \cap L^\infty(\Gamma_X;\mu)$

$$\begin{aligned}\int FG \operatorname{div}_\mu V \, d\mu &= \int G \operatorname{div}_\mu^\Gamma(FV)\, d\mu - \int G\langle \nabla^\Gamma F, V\rangle_{T\Gamma_X}\, d\mu \\ &= -\int \langle d^\mu G, FV\rangle_{T\Gamma_X}\, d\mu - \int \langle G d^\mu F, V\rangle_{T\Gamma_X}\, d\mu\end{aligned}$$

where we used (4.10), (4.11) and (4.13) in the last step. Since $\mathcal{VFC}_b^{\infty,\mu}$ is dense in $L^2(\Gamma_X \to T\Gamma_X;\mu)$, the assertion follows provided $F \in \mathcal{FC}_b^{\infty,\mu}$. The case where $F \in H_0^{1,2}(\Gamma_X;\mu) \cap L^\infty(\Gamma_X;\mu)$ is then an immediate consequence, since $\mathcal{FC}_b^{\infty,\mu}$ is dense in $H_0^{1,2}(\Gamma_X;\mu)$ w.r.t. $\| \; \|_{H_0^{1,2}(\Gamma_X;\mu)}$ by definition. □

**Proposition 4.6** *Suppose conditions (A) and (B) in Theorem 3.3 hold and let $\mu_{\lambda,m}$ (or even $\mu_{\lambda,\sigma}$) be a mixed Poisson measure satisfying* (3.1). *Then*

$$\mathcal{E}^\Gamma_{\mu_{\lambda,m}}(F,G) = \int H^\Gamma_{\mu_{\lambda,m}} F \; G \, d\mu_{\lambda,m} \text{ for all } F \in E(\mathcal{D}_1), G \in W^{1,2}(\Gamma_X;\mu) \tag{4.15}$$

*(cf.* (3.10), (3.11)*). In particular,*

$$H_0^{1,2}(\Gamma_X;\mu_{\lambda,m}) = W^{1,2}(\Gamma_X;\mu_{\lambda,m}). \tag{4.16}$$

**Proof.** Let $G \in W^{1,2}(\Gamma_X;\mu_{\lambda,m})$ and $F \in E(\mathcal{D}_1)$. By Proposition 3.7 (ii) there exist $F_n \in \mathcal{FC}_b^\infty$, $n \in \mathbb{N}$, such that $F_n \to F$ as $n \to \infty$ w.r.t. $\| \cdot \|_1 + \| \cdot \|_{H_0^{1,2}(\Gamma_X;\mu_{\lambda,m})}$. Hence, because $(d^\mu, W^{1,2}(\Gamma_X;\mu_{\lambda,m}))$ extends $(\nabla^\Gamma_\mu, H_0^{1,2}(\Gamma_X;\mu_{\lambda,m}))$ by definition of $d^\mu$

$$\begin{aligned}\mathcal{E}^\Gamma_\mu(F,G) &= \lim_{n\to\infty} \int \langle \nabla^\Gamma F_n, d^\mu G\rangle_{T\Gamma_X}\, d\mu_{\lambda,m} \\ &= -\lim_{n\to\infty} \int \Delta^\Gamma F_n \; G \, d\mu_{\lambda,m} = \int H^\Gamma_{\mu_{\lambda,m}} F \; G\, d\mu_{\lambda,m}.\end{aligned}$$

This implies that the generator of $(\mathcal{E}^\Gamma_\mu; W^{1,2}(\Gamma_X;\mu_{\lambda,m}))$, which is self–adjoint, extends $(H^\Gamma_{\mu_{\lambda,m}}, E(\mathcal{D}_1))$ (cf. e.g. [MR92, Chap. I, Proposition 2.16]), hence it also extends its closure on $L^2(\Gamma_X;\mu_{\lambda,m})$, which by Proposition 3.7 (ii) is self–adjoint, too. Since self–adjoint operators cannot strictly extend each other, they must coincide. But by Proposition 3.7 (ii) the closure of $(H^\Gamma_{\mu_{\lambda,m}}, E(\mathcal{D}_1))$ is equal to its Friedrichs' extension, i.e., equal to the generator of $(\mathcal{E}^\Gamma_\mu, H_0^{1,2}(\Gamma_X;\mu_{\lambda,m}))$, hence the assertion follows. □

In general, it might be that $H_0^{1,2}(\Gamma_X;\mu) \subsetneqq W^{1,2}(\Gamma_X;\mu)$ (though we expect that equality holds for Ruelle measures to be discussed in Subsection 7.1 below under quite weak assumptions) and that $(\mathcal{E}_\mu^\Gamma, W^{1,2}(\Gamma_X;\mu)$ might not be a Dirichlet form. Therefore, in [RSch97] two useful "intermediate" spaces $H_0^{1,2}(\Gamma_X;\mu) \subset \mathcal{F}^{(c)} \subset \mathcal{F} \subset W^{1,2}(\Gamma_X;\mu)$ were introduced so that $(\mathcal{E}_\mu^\Gamma, \mathcal{F}^{(c)})$, $(\mathcal{E}_\mu^\Gamma, \mathcal{F})$ are Dirichlet forms on $L^2(\Gamma_X;\mu)$. We shall briefly recall the corresponding definitions and results. We need the following additional assumption on $\mu$:

(QI)(i) For any $n \in \mathbb{N}$ either $\mu(\{\gamma \in \Gamma_X \mid \gamma(X) = n\}) > 0$ ore $\mu(\{\gamma \in \Gamma_X \mid \gamma(X) \geq n\}) > 0$ corresponding to whether $X$ is compact or non–compact respectively and for all $v \in V_0(X)$ and $t \in \mathbb{R}$, $\mu$ is quasi–invariant w.r.t. $\psi_t^v$, i.e., $\mu \circ (\psi_t^v)^{-1} \approx \mu$.

(ii) $\operatorname*{essinf}_{r \leq s \leq t} \Phi_s > 0$ $\mu$–a.e. for all $r, t \in \mathbb{R}$, $r < t$, where

$$\Phi_s := \frac{d(\mu \circ (\psi_s^v)^{-1}) \otimes ds}{d\mu \otimes ds}.$$

By [AKR97a, Proposition 2.2] it follows that every mixed Poisson measure satisfies (QI). The same is true for a class of Ruelle measures to be studied in Section 7 below (cf. Remark 7.6 (iii)).

**Remark 4.7** It can be shown that (QI) (i) implies that $\operatorname{supp}\mu = \Gamma_X$ (cf. [RSch97, Proposition 5.6]). Hence, in particular, (4.3) holds in this case without assumption (IbP).

We now assume that $\mu$ satisfies both (IbP) and (QI)(i). Let $\mathbb{F}$ denote the set of all $F \in L^\infty(\Gamma_X;\mu)$ such that there exists $\nabla^\Gamma F \in L^2(\Gamma_X \to T\Gamma_X;\mu)$ such that for all $v \in V_0(X)$, $s \in \mathbb{R}$

$$\frac{F \circ \psi_{s+t}^v - F \circ \psi_s^v}{t} \underset{t \to 0}{\longrightarrow} \langle \nabla^\Gamma F, v \rangle_{T\Gamma_X} \circ \psi_s^v \quad \text{in } L^2(\Gamma_X;\mu).$$

Let $\mathbb{F}^{(c)}$ denote the subset of $\mathbb{F}$ having a $\mu$–version in $C(\Gamma_X)$, where $C(\Gamma_X)$ denotes the set of all continuous functions on $\Gamma_X$. Clearly, $\mathcal{F}C_b^{\infty,\mu} \subset \mathbb{F}^{(c)}$ and $\nabla_\mu^\Gamma = \nabla^\Gamma$ on $\mathcal{F}C_b^{\infty,\mu}$. But we also have:

**Lemma 4.8** *Let $\mu$ satisfy (IbP) and (QI). Then $\mathbb{F} \subset W^{1,2}(\Gamma_X;\mu)$ and $\nabla^\Gamma = d^\mu$ on $\mathbb{F}$.*

**Proof.** [RSch97, Lemma 6.3]. $\square$

**Proposition 4.9** *(i) Let $\mu$ satisfy (IbP) and (QI)(i). Then $(\mathcal{E}_\mu^\Gamma, \mathbb{F})$, $(\mathcal{E}_\mu^\Gamma, \mathbb{F}^{(c)})$ are closable on $L^2(\Gamma_X;\mu)$. Their respective closures $(\mathcal{E}_\mu^\Gamma, \mathcal{F})$, $(\mathcal{E}_\mu^\Gamma, \mathcal{F}^{(c)})$ are symmetric Dirichlet forms on $L^2(\Gamma_X;\mu)$. Furthermore, $D(\mathcal{E}_\mu^\Gamma) \subset \mathcal{F}^{(c)} \subset \mathcal{F} \subset W^{1,2}(\Gamma_X;\mu)$.*

*(ii) If $\mu = \mu_{\lambda,m}$ (or $\mu_{\lambda,\sigma}$) is a mixed Poisson measure satisfying* (3.1) *and if conditions (A) and (B) in Theorem 3.3 hold, then*

$$D(\mathcal{E}_\mu^\Gamma) = \mathcal{F}^{(c)} = \mathcal{F} = W^{1,2}(\Gamma_X;\mu_{\lambda,m}).$$

**Proof.** (i): The closability of both forms and the last assertion are immediate consequences of Lemma 4.8. The proof that both $(\mathcal{E}_\mu^\Gamma, \mathcal{F})$ and $(\mathcal{E}_\mu^\Gamma, \mathcal{F}^{(c)})$ are Dirichlet forms is quite straightforward. We refer to [RSch97, Proof of Proposition 1.4(i)].

(ii): This follows by Proposition 4.6. $\square$

We would like to emphasize that even only trying to prove directly that $D(\mathcal{E}_{\mu_{\lambda,m}}^\Gamma) = \mathcal{F}^{(c)}$ (i.e., that $\mathcal{F}C_b^\infty$ is dense in $\mathcal{F}^{(c)}$ w.r.t. $\|\cdot\|_{W^{1,2}(\Gamma_X;\mu_{\lambda,m})}$) appears extremely hard. The elaborate machinery about uniqueness of semigroup generators when restricted to subsets of their full domain, on which the proof of Proposition 4.9 is based, seems unavoidable even in this simple case where $\mu$ is a mixed Poisson measure.

There are many reasons to study $(\mathcal{E}_\mu^\Gamma, \mathcal{F})$ and $(\mathcal{E}_\mu^\Gamma, \mathcal{F}^{(c)})$. One is that one can prove a Rademacher theorem on $\Gamma_X$ (cf. Section 5 below). Others are presented in [AKRR98].

## 4.2 Quasi–regularity and corresponding diffusions

Suppose that $\mu$ is a probability measure on $(\Gamma_X, \mathcal{B}(\Gamma_X))$ satisfying (4.1) and (4.3) such that

$$(\mathcal{E}_\mu^\Gamma, \mathcal{F}C_b^{\infty,\mu}) \text{ is closable on } L^2(\Gamma_X;\mu). \tag{4.17}$$

As before the closure is denoted by $(\mathcal{E}_\mu^\Gamma, D(\mathcal{E}_\mu^\Gamma))$. The generators of $(\mathcal{E}_\mu^\Gamma, D(\mathcal{E}_\mu^\Gamma))$ we denote again by $(H_\mu^\Gamma, D(H_\mu^\Gamma))$ but emphasize that because we do not assume that (IbP) holds, it is in general not equal to the Friedrichs' extension of an operator such as $\Delta_\mu^\Gamma$ explicitly given on $\mathcal{F}C_b^{\infty,\mu}$. It might be in general not even true that $\mathcal{F}C_b^{\infty,\mu} \subset D(H_\mu^\Gamma)$.

As in Subsection 3.3 the desired processes in general will live on the bigger state space $\ddot{\Gamma}_X$. To state the corresponding existence result precisely we need Definitions 4.10, 4.11 below. As in Subsection 3.3 we can consider $\mu$ as a measure on $(\ddot{\Gamma}_X, \mathcal{B}(\ddot{\Gamma}_X))$ and correspondingly $(\mathcal{E}_\mu^\Gamma, D(\mathcal{E}_\mu^\Gamma))$ as a Dirichlet form on $L^2(\ddot{\Gamma}_X; \mu)$.

**Definition 4.10** (*cf.* [MR92, Ch.III, Definitions 2.1 and 3.2]). *Let $(\mathcal{E}, D(\mathcal{E}))$ be a Dirichlet form on $L^2(\ddot{\Gamma}_X; \mu)$.*

(i) *A sequence $(F_n)_{n\in\mathbb{N}}$ of closed subsets of $\ddot{\Gamma}_X$ is called an $\mathcal{E}$-nest, if*

$$\bigcup_{n\in\mathbb{N}} \{F \in D(\mathcal{E}) | F = 0 \ \mu\text{-a.e. on } \ddot{\Gamma}_X \setminus F_n\}$$

*is dense in $D(\mathcal{E})$ w.r.t. $\mathcal{E}_1^{1/2} (:= [\mathcal{E} + (\ ,\ )_{L^2(\mu)}]^{1/2})$.*

(ii) *A set $N \subset \ddot{\Gamma}_X$ is called $\mathcal{E}$-exceptional, if $N \subset \ddot{\Gamma}_X \setminus \bigcup_{n\in N} F_n$ for some $\mathcal{E}$-nest $(F_n)_{n\in\mathbb{N}}$. We say that a property of points in $\ddot{\Gamma}_X$ holds $\mathcal{E}$-quasi-everywhere (abbeviated $\mathcal{E}$-q.e.) if it holds outside some $\mathcal{E}$-exceptional set.*

(iii) *A function $f : \ddot{\Gamma}_X \to \mathbb{R}$ is called $\mathcal{E}$-quasi-continuous if there exists an $\mathcal{E}$-nest $(F_n)_{n\in\mathbb{N}}$ such that the restriction $f_{|F_n}$ of $f$ to $F_n$ is continuous for all $n \in \mathbb{N}$.*

We recall that the notion "$\mathcal{E}$-q.e." is much finer than "$\mu$-a.e." (cf. [MR92, Exercise 2.6 (ii)]).

**Definition 4.11** (cf. [MR92, Chap.IV, Definition 3.1]) A Dirichlet form $(\mathcal{E}, D(\mathcal{E}))$ on $L^2(\ddot{\Gamma}_X; \mu)$ is called *quasi–regular* if:

(i) There exists an $\mathcal{E}$–nest $(F_k)_{k\in\mathbb{N}}$ consisting of compact subsets of $\ddot{\Gamma}_X$.

(ii There exists a $\mathcal{E}_1^{1/2}$–dense subset of $D(\mathcal{E})$ whose elements have $\mathcal{E}$–quasi–continuous $\mu$–versions.

(iii) There exist $u_n \in D(\mathcal{E})$, $n \in \mathbb{N}$, having $\mathcal{E}$–quasi–continuous $\mu$–versions $\tilde{u}_n$, $n \in \mathbb{N}$, and an $\mathcal{E}$–exceptional set $N \subset \ddot{\Gamma}_X$ such that $\{\tilde{u}_n \mid n \in \mathbb{N}\}$ separates the points of $\ddot{\Gamma}_X \setminus N$.

The following is a special case of the main result in [MR97].

**Theorem 4.12** *$(\mathcal{E}_\mu^\Gamma, D(\mathcal{E}_\mu^\Gamma))$ is a quasi–regular Dirichlet form on $L^2(\ddot{\Gamma}_X;\mu)$.*

We do not include a proof of Theorem 4.12 here, but refer to [MR97].

**Corollary 4.13** *$(\mathcal{E}_\mu^\Gamma, D(\mathcal{E}_\mu^\Gamma))$ is* local *(i.e., $\mathcal{E}_\mu^\Gamma(F,G) = 0$ provided $F,G \in D(\mathcal{E}_\mu^\Gamma)$ with $\operatorname{supp}(|F|\mu) \cap \operatorname{supp}(|G|\mu) = \emptyset$).*

**Proof.** Since $(\mathcal{E}_\mu^\Gamma, D(\mathcal{E}_\mu^\Gamma))$ is quasi–regular by Theorem 4.12 and since $\nabla_\mu^\Gamma$ satisfies the product rule on bounded functions in $D(\mathcal{E}_\mu^\Gamma)$, the proof for locality is entirely analogous to that in [MR92, Chap. V, Examples 1.12 (ii)]. (We also refer to the beautiful work [S95] for a complete discussion of locality of Dirichlet forms.) □

As a consequence of Theorem 4.12, Corollary 4.13 and [MR92, Chap. IV, Theorem 3.5 and Chap. V, Theorem 1.11] we immediately obtain

**Theorem 4.14** *There exists a conservative diffusion process*

$$\mathbf{M} = (\mathbf{\Omega}, \mathbf{F}, (\mathbf{F}_t)_{t\geq 0}, (\mathbf{\Theta}_t)_{t\geq 0}, (\mathbf{X}_t)_{t\geq 0}, (\mathbf{P}_\gamma)_{\gamma\in\ddot{\Gamma}_X})$$

*on $\ddot{\Gamma}_X$ (cf. [Dy65]) which is properly associated with $(\mathcal{E}_\mu^\Gamma, D(\mathcal{E}_\mu^\Gamma))$, i.e., for all ($\mu$-versions of) $F \in L^2(\ddot{\Gamma}_X;\mu)$ and all $t > 0$ the function*

$$\gamma \ \mapsto \ p_t F(\gamma) \ := \ \int_\Omega F(\mathbf{X}_t) \; d\mathbf{P}_\gamma, \ \ \gamma \in \ddot{\Gamma}_X, \qquad (4.18)$$

*is an $\mathcal{E}_\mu^\Gamma$-quasi-continuous version of $\exp(-tH_\mu^\Gamma)F$. $\mathbf{M}$ is up to $\mu$-equivalence unique* (*cf.* [MR92, Chap. IV, Sect. 6]). *In particular, $\mathbf{M}$ is $\mu$–symmetric (i.e., $\int G \; p_t F d\mu = \int F \; p_t G \; d\mu$ for all $F,G : \ddot{\Gamma}_X \to \mathbb{R}_+, \mathcal{B}(\ddot{\Gamma}_X)$-measurable) and has $\mu$ as an invariant measure.*

**Proof.** By Theorem 4.12 and Corollary 4.13 the proof follows directly from [MR92, Chap. V, Theorem 1.11]. □

In the above theorem $\mathbf{M}$ is canonical, i.e., $\mathbf{\Omega} = C([0,\infty) \to \ddot{\Gamma}_X)$, $\mathbf{X}_t(\omega) := \omega(t)$, $t \geq 0$, $\omega \in \mathbf{\Omega}$, $(\mathbf{F}_t)_{t\geq 0}$ together with $\mathbf{F}$ is the corresponding minimum completed admissible family (cf. [F80, Sect.4.1]) and $\mathbf{\Theta}_t$, $t \geq 0$, are the corresponding natural time shifts.

**Remark 4.15** (i) Theorems 4.12, 4.14 and Corollary 4.13 generalize corresponding results in [Os96], [Y96].

(ii) $\mathbf{M}$ is also called the corresponding *stochastic dynamics*, or *stochastic quantization corresponding to* $\mu$.

**Theorem 4.16** $\mathbf{M}$ *from Theorem 4.14 is the (up to $\mu$–equivalence, cf. [MR92, Chap. IV, Definition 6.3]) unique diffusion process having $\mu$ as an invariant measure and solving the martingale problem for $(-H_\mu^\Gamma, D(H_\mu^\Gamma))$, i.e., for all $G \in D(H_\mu^\Gamma)$ $(\supset \mathcal{F}C_b^\infty(\mathcal{D}, \Gamma))$*

$$\tilde{G}(\mathbf{X}_t) - \tilde{G}(\mathbf{X}_0) + \int_0^t H_\mu^\Gamma G(\mathbf{X}_s)\, ds, \ t \geq 0,$$

*is an $(\mathbf{F}_t)$–martingale under $\mathbf{P}_\gamma$ (hence starting at $\gamma$) for $\mathcal{E}_\mu^\Gamma$–q.e. $\gamma \in \ddot{\Gamma}_X$. (Here $\tilde{G}$ denotes a quasi–continuous version of $G$, cf. [MR92, Chap. IV, Proposition 3.3]).*

**Proof.** This follows immediately from [AR95, Theorem 3.5]. □

**Remark 4.17** In fact, the uniqueness statement in Theorem 4.16 can be strengthened, as follows:
$\mathbf{M}$ is (up to $\mu$–equivalence) unique among all right processes $\mathbf{M}' = (\mathbf{\Omega}', \mathbf{F}', (\mathbf{F}'_t)_{t\geq 0}, (\mathbf{\Theta}'_t)_{t\geq 0}, (\mathbf{X}')_{t\geq 0}, (\mathbf{P}'_\gamma)_{\gamma\in\ddot{\Gamma}_X})$ on $\ddot{\Gamma}_X$ having $\mu$ as a sub–invariant measure and being such that for all $G \in D(H_\mu^\Gamma)$

$$\tilde{G}(\mathbf{X}'_t) - \tilde{G}(\mathbf{X}'_0) + \int_0^t H_\mu^\Gamma G(\mathbf{X}'_s)\, ds, \ t \geq 0,$$

is an $(\mathbf{F}'_t)$–martingale under $\mathbf{P}'_\mu := \int \mathbf{P}'_\gamma \ \mu(d\gamma)$.

The proof is analogous to that of [AR95, Theorem 3.5] except for its ending where one has to use the recent result of P. Fitzsimmons in [Fi96] in order to be able to apply the arguments in [MR92, Chap. IV, Section 6] to deduce the result.

# 5 Intrinsic metric on configuration spaces

In this section we give a summary of the main results in the recent work [RSch97]. Below we fix a probability measure $\mu$ on $(\Gamma_X, \mathcal{B}(\Gamma_X))$ satisfying the conditions (IbP) and (QI) specified in Subsection 4.1, to which we also refer for examples. Let also $(\mathcal{E}_\mu^\Gamma, D(\mathcal{E}_\mu^\Gamma))$, $(\mathcal{E}_\mu^\Gamma, \mathcal{F}^{(c)})$, $(\mathcal{E}_\mu^\Gamma, \mathcal{F})$ be as in Subsection 4.1, considered as Dirichlet forms on $L^2(\ddot{\Gamma}_X; \mu)$.

## 5.1 Identification of the intrinsic metric for a class of measures

The following $L^2$–*Wasserstein type distance* $\rho$ on $\ddot{\Gamma}_X$ was introduced in [RSch97]:

$$\rho(\gamma,\omega) := \inf\Big\{ \Big( \int d(x,y)^2 \; \eta(dxdy) \Big)^{1/2} \mid \eta \in \ddot{\Gamma}_{X\times X} \text{ with marginals } \omega,\gamma \text{ respectively} \Big\};\ \gamma,\omega \in \ddot{\Gamma}_X \tag{5.1}$$

where $d$ is the Riemannian distance on $X$. Note that $\rho$ is a pseudo–metric, in particular, $\rho(\omega,\gamma) = 0$ implies $\omega = \gamma$. But $\rho(x,y)$ will be infinite if $\omega(X) \neq \gamma(X)$, since there is no $\eta \in \ddot{\Gamma}_{X\times X}$ with marginals $\omega$, $\gamma$ respectively. But also if both $\omega$ and $\gamma$ are infinite configurations one will find that possibly $\rho(\omega,\gamma) = \infty$ (e.g. take $X = \mathbb{R}$, $\omega = \sum_{z\in\mathbb{Z}} \varepsilon_z$, and $\gamma = \omega - \delta_0$). Obviously, the topology induced on $\Gamma_X$ by $\rho$ is finer than the vague topology.

**Theorem 5.1** *The intrinsic metrics of the Dirichlet forms* $(\mathcal{E}^\Gamma_\mu, \mathcal{F})$ *and* $(\mathcal{E}^\Gamma_\mu, \mathcal{F}^{(c)})$ *are both equal to the* $L^2$*–Wasserstein type distance* $\rho$ *defined in* (5.1)*, more precisely:*

$$\begin{aligned}
& \sup\{F(\gamma) - F(\omega) \mid F \in \mathcal{F} \cap C(\Gamma_X) \text{ and } |\nabla^\Gamma F|_{T\Gamma_X} \le 1 \ \mu\text{–a.e. on } \ddot{\Gamma}_X\} \\
= \ & \sup\{F(\gamma) - F(\omega) \mid F \in \mathcal{F}^{(c)} \cap C(\Gamma_X) \text{ and } |\nabla^\Gamma F|_{T\Gamma_X} \le 1 \ \mu\text{–a.e. on } \ddot{\Gamma}_X\} \\
= \ & \sup\{F(\gamma) - F(\omega) \mid F \in \mathbb{F}^{(c)} \text{ and } |\nabla^\Gamma F|_{T\Gamma_X} \le 1 \ \mu\text{–a.e. on } \ddot{\Gamma}_X\} \\
= \ & \rho(\gamma,\omega).
\end{aligned}$$

**Proof.** See [RSch97, Sect. 8]. □

## 5.2 A Rademacher theorem

We recall that a function $F : \ddot{\Gamma}_X \to \mathbb{R}$ is called $\rho$–Lipschitz continuous if

$$\mathrm{Lip}\,(F) := \sup\Big\{ \frac{|F(\gamma) - F(\omega)|}{\rho(\gamma,\omega)} \mid \gamma,\omega \in \ddot{\Gamma}_X,\ \gamma \neq \omega \Big\} < \infty. \tag{5.2}$$

Concerning generalizations, different from the following theorem, of the classical Rademacher theorem (cf. [Ra19]) to infinite dimensional spaces (though in contrast to ours only “flat” ones), we refer to the introduction of [RSch97] and the references therein.

**Theorem 5.2** *(i) Suppose $F \in L^2(\ddot{\Gamma}_X;\mu)$ is $\rho$–Lipschitz continuous. Then $F \in \mathcal{F}$ and*

$$|\nabla^\Gamma F|_{T\Gamma_X} \leq \mathrm{Lip}\,(F) \quad \mu - a.e.$$

*(ii) Suppose $F \in \mathcal{F}$ such that $|\nabla^\Gamma F|_{T\Gamma_X} \leq c$ $\mu$–a.e. for some $c \in ]0,\infty[$ and such that $F$ has a $\rho$–continuous $\mu$–version. Then there exists a $\rho$–Lipschitz continuous $\mu$–version $\widetilde{F}$ of $F$ satisfying*

$$\mathrm{Lip}\,(\widetilde{F}) \leq c\ .$$

**Proof.** [RSch97, Sections 7 and 8]. □

## 5.3 Potential-theoretic consequences

We now state some of the results on the potential theory on $\ddot{\Gamma}_X$, induced by the Dirichlet forms $(\mathcal{E}^\Gamma_\mu, D(\mathcal{E}^\Gamma_\mu))$ and $(\mathcal{E}^\Gamma_\mu, \mathcal{F}^{(c)})$ on $L^2(\ddot{\Gamma}_X;\mu)$, obtained as consequences of Theorems 5.1, 5.2 respectively their proofs in [RSch97]. For proofs we refer to [RSch97, Sect. 3].

Define for $A \subset \ddot{\Gamma}_X$

$$\rho_A(\gamma) := \inf\{\rho(\omega,\gamma) \mid \omega \in A\}. \tag{5.3}$$

Let $\bar{A}^\rho$ denote the closure of $A$ w.r.t. the topology induced by $\rho$.

**Proposition 5.3** *If $K \subset \ddot{\Gamma}_X$ is compact and $c \geq 0$, then $c \wedge \rho_K$ is an $(\mathcal{E}^\Gamma_\mu, \mathcal{F}^{(c)})$–quasi–continuous function in $\mathcal{F}^{(c)}$.*

**Corollary 5.4** *If $A \in \mathcal{B}(\ddot{\Gamma}_X)$ such that $\mu(A) = 1$, then $\ddot{\Gamma}_X \setminus \bar{A}^\rho = \{\rho_A > 0\}$ is $(\mathcal{E}^\Gamma_\mu, \mathcal{F}^{(c)})$–exceptional.*

**Example 5.5** If $A = \{\gamma \in \ddot{\Gamma}_X \mid \gamma(X) = \infty\}$ one sees immediately that $A = \{\rho_A < \infty\}$. Hence, as follows from Corollary 5.4, $\mu(A) = 1$ implies that $\ddot{\Gamma}_X \setminus A$ is $(\mathcal{E}^\Gamma_\mu, \mathcal{F}^{(c)})$–exceptional, whereas $\mu(A) = 0$ implies that $A$ is $(\mathcal{E}^\Gamma_\mu, \mathcal{F}^{(c)})$–exceptional. This result has first been proven by Byron Schmuland (private communication).

We note that since $\mathcal{F}C_b^{\infty,\mu} \subset D(\mathcal{E}^\Gamma_\mu) \subset \mathcal{F}^{(c)} \subset \mathcal{F}$, both $(\mathcal{E}^\Gamma_\mu, \mathcal{F}^{(c)})$ and $(\mathcal{E}^\Gamma_\mu, \mathcal{F})$ obviously fulfill property (iii) in the definition of quasi–regularity (cf. Definition 4.11) with $N = \emptyset$. (Just pick $u_n \in \mathcal{F}C_b^\infty$, $n \in \mathbb{N}$, separating the points of $\ddot{\Gamma}_X$). Furthermore, by [MR92, Chap. III, Theorem 2.11(i) and Excercise 2.10] it follows that any sequence $(F_k)_{k\in\mathbb{N}}$ which is a nest w.r.t. $(\mathcal{E}_\mu^\Gamma, D(\mathcal{E}_\mu^\Gamma))$ is also a nest w.r.t. $(\mathcal{E}_\mu^\Gamma, \mathcal{F}^{(c)})$ and $(\mathcal{E}_\mu^\Gamma, \mathcal{F})$. In particular, then by Theorem 4.12 both $(\mathcal{E}_\mu^\Gamma, \mathcal{F}^{(c)})$ and $(\mathcal{E}_\mu^\Gamma, \mathcal{F})$ also have property (i) in Definition 4.11. Since $(\mathcal{E}_\mu^\Gamma, \mathcal{F}^{(c)})$ by definition also satisfies 4.11 (ii). We have:

**Proposition 5.6** *$(\mathcal{E}_\mu^\Gamma, \mathcal{F}^{(c)})$ is quasi–regular on $L^2(\ddot{\Gamma}_X; \mu)$. Moreover, $(\mathcal{E}_\mu^\Gamma, \mathcal{F}^{(c)})$ is local.*

**Proof.** We have just seen that $(\mathcal{E}_\mu^\Gamma, \mathcal{F}^{(c)})$ is quasi–regular. Clearly, $\nabla^\Gamma$ satisfies the product rule on $\mathbb{F}^{(c)}$, hence on $\mathcal{F}^{(c)}$. Consequently, the proof of locality is entirely analogous to that in [MR92, Chap. V, Examples 1.12 (ii)]. □

But in the present situation one can even give an "explicit" nest consisting of compacts:

**Proposition 5.7** *Let $K_n \subset \ddot{\Gamma}_X$, $K_n \subset K_{n+1}$, $n \in \mathbb{N}$, be compact such that $\mu(K_n) \uparrow 1$. Then*

$$F_n := \{\rho_{K_n} \leq \frac{1}{2}\},\ n \in \mathbb{N},$$

*are compact and form a $(\mathcal{E}_\mu^\Gamma, \mathcal{F}^{(c)})$–nest.*

**Proof.** [RSch97, Corollary 3.4]. □

**Theorem 5.8** *The analogous of Theorems 4.14 and 4.16 hold with $(\mathcal{E}_\mu^\Gamma, \mathcal{F}^{(c)})$ replacing $(\mathcal{E}_\mu^\Gamma, D(\mathcal{E}_\mu^\Gamma))$ resp. with the generator of $(\mathcal{E}_\mu^\Gamma, \mathcal{F}^{(c)})$ replacing $(H_\mu^\Gamma, D(H_\mu^\Gamma))$.*

**Proof.** Analogous to the proofs of Theorems 4.14 and 4.16. □

# 6 Gibbs measures on configuration spaces

So far, the only concrete examples of probability measures $\mu$ on $(\Gamma_X, \mathcal{B}(\Gamma_X))$ have been the mixed Poisson measures. In this section we will describe another class of measures on $(\Gamma_X, \mathcal{B}(\Gamma_X))$, so–called *Gibbs measures.* Our presentation is (apart from slight notational modifications) very much oriented along the beautiful work by C. Preston (see [P79], [P80], but also [P76] and [Ge79]).

In this section we formulate all results right away for the more general intensity measures $\sigma = \rho \cdot m$ specified at the beginning of Remark 2.3 (i). So, the underlying (reference) Poisson measure is $\pi_\sigma$ rather than $\pi_m$.

## 6.1 Grand canonical Gibbs measures

Let $\mathcal{B}_b(\Gamma_X)$ denote the set of all bounded $\mathcal{B}(\Gamma_X)$ - measurable functions $G : \Gamma_X \to \mathbb{R}$. If for $B \in \mathcal{B}(X)$ we define $N_B : \Gamma_X \to \mathbb{Z}_+ \cup \{+\infty\}$ by

$$N_B(\gamma) \;:=\; \gamma(B), \tag{6.1}$$

then

$$\mathcal{B}(\Gamma_X) \;=\; \sigma(\{N_\Lambda | \Lambda \in \mathcal{O}_c(X)\}). \tag{6.2}$$

For $A \in \mathcal{B}(X)$ we define

$$\mathcal{B}_A(\Gamma_X) \;:=\; \sigma(\{N_B | B \in \mathcal{B}_c(X), B \subset A\}). \tag{6.3}$$

Let $\Phi$ be a *potential*, i.e., a function $\Phi : \Gamma_X \to \mathbb{R} \cup \{+\infty\}$ such that

$$\begin{gathered} \Phi \;=\; 1_{\{N_X < \infty\}}\Phi, \;\; \Phi(\emptyset) \;=\; 0, \;\; \text{and} \\ \gamma \;\mapsto\; \Phi(\gamma_\Lambda) \;\; \text{is} \;\; \mathcal{B}_\Lambda(\Gamma_X) - \text{measurable} \\ \text{for all} \;\; \Lambda \in \mathcal{B}_c(X). \end{gathered} \tag{6.4}$$

For $\Lambda \in \mathcal{B}_c(X)$ the *conditional energy* $E_\Lambda^\Phi : \Gamma_X \to \mathbb{R} \cup \{+\infty\}$ is defined by

$$E_\Lambda^\Phi(\gamma) \;:=\; \begin{cases} \sum_{\gamma' \subset \gamma, \gamma'(\Lambda) > 0} \Phi(\gamma') & \text{if } \sum_{\gamma' \subset \gamma, \gamma'(\Lambda) > 0} |\Phi(\gamma')| < \infty, \\ +\infty & \text{otherwise,} \end{cases} \tag{6.5}$$

(where the sum over the empty set is defined to be zero).

Now we can define grand canonical Gibbs measures:

**Definition 6.1** *For $\Lambda \in \mathcal{O}_c(X)$ define for $\gamma \in \Gamma_X, \Delta \in \mathcal{B}(\Gamma_X)$*

$$\Pi_\Lambda^{\sigma,\Phi}(\gamma, \Delta) := 1_{\{Z_\Lambda^{\sigma,\Phi}<\infty\}}(\gamma)[Z_\Lambda^{\sigma,\Phi}(\gamma)]^{-1} \int_{\Gamma_X} 1_\Delta(\gamma_{X\setminus\Lambda} + \gamma'_\Lambda) \times \exp[-E_\Lambda^\Phi(\gamma_{X\setminus\Lambda} + \gamma'_\Lambda)]\ \pi_\sigma(d\gamma'),$$

*where*

$$Z_\Lambda^{\sigma,\Phi}(\gamma) := \int_{\Gamma_X} \exp[-E_\Lambda^\Phi(\gamma_{X\setminus\Lambda} + \gamma'_\Lambda)]\ \pi_\sigma(d\gamma').$$

*A probability measure $\mu$ on $(\Gamma_X, \mathcal{B}(\Gamma_X))$ is called a* grand canonical Gibbs measure with interaction potential $\Phi$ *if*

$$\mu\Pi_\Lambda^{\sigma,\Phi} = \mu \text{ for all } \Lambda \in \mathcal{O}_c(X). \tag{6.6}$$

*Let $\mathcal{G}_{gc}(\sigma, \Phi)$ denote the set of all such probability measures $\mu$.*

**Remark 6.2** (i) It is well-known that $(\Pi_\Lambda^{\sigma,\Phi})_{\Lambda\in\mathcal{O}_c(X)}$ is a $(\mathcal{B}_{X\setminus\Lambda}(\Gamma_X))_{\Lambda\in\mathcal{O}_c(X)}$ - specification in the sense of [Fö75], [P76] (cf. [P76, Section 6] or [P79]), i.e., for all $\Lambda, \Lambda' \in \mathcal{O}_c(X)$

(S 1) $\Pi_\Lambda^{\sigma,\Phi}(\gamma, \Gamma_X) \in \{0, 1\}$ for all $\gamma \in \Gamma_X$.

(S 2) $\Pi_\Lambda^{\sigma,\Phi}(\cdot, \Delta)$ is $\mathcal{B}_{X\setminus\Lambda}(\Gamma_X)$-measurable for all $\Delta \in \mathcal{B}(\Gamma_X)$.

(S 3) $\Pi_\Lambda^{\sigma,\Phi}(\cdot, \Delta' \cap \Delta) = 1_{\Delta'}\Pi_\Lambda(\cdot, \Delta)$ for all $\Delta \in \mathcal{B}(\Gamma_X), \Delta' \in \mathcal{B}_{X\setminus\Lambda}(\Gamma_X)$.

(S 4) $\Pi_{\Lambda'}^{\sigma,\Phi} = \Pi_{\Lambda'}^{\sigma,\Phi}\Pi_\Lambda^{\sigma,\Phi}$ if $\Lambda \subset \Lambda'$.

Here for $\gamma \in \Gamma_X, \Delta \in \mathcal{B}(\Gamma_X)$

$$(\Pi_{\Lambda'}^{\sigma,\Phi}\Pi_\Lambda^{\sigma,\Phi})(\gamma, \Delta) := \int \Pi_\Lambda^{\sigma,\Phi}(\gamma', \Delta)\ \Pi_{\Lambda'}^{\sigma,\Phi}(\gamma, d\gamma')$$

(and $\mu\Pi_\Lambda^{\sigma,\Phi}$ in (6.6) is defined correspondingly). (6.6) above are called *Dobrushin-Lanford-Ruelle equations.* We also note that as a direct consequence of (S 3) we have for all $\Lambda \in \mathcal{O}_c(X)$ that

$$\mu(\{Z_\Lambda^{\sigma,\Phi} < \infty\}) = 1$$

for all $\mu \in \mathcal{G}_{gc}(\sigma, \Phi)$.

(ii) It is an easy exercise (see e.g. [P80, Proposition 10.3]) to show that because of (2.17) for all $\Lambda \in \mathcal{O}_c(X)$, $\gamma \in \Gamma_X$, $\Delta \in \mathcal{B}(\Gamma_X)$

$$\Pi_\Lambda^{\sigma,\Phi}(\gamma, \Delta) = 1_{\{\tilde{Z}_\Lambda^{\sigma,\Phi}<\infty\}}(\gamma)[\tilde{Z}_\Lambda^{\sigma,\Phi}(\gamma)]^{-1}[1_\Delta(\gamma_{X\setminus\Lambda}) +$$

$$\sum_{n=1}^{\infty} \frac{1}{n!} \int_\Lambda \dots \int_\Lambda 1_\Delta(\gamma_{X\setminus\Lambda} + \sum_{k=1}^{n} \varepsilon_{x_k}) \exp[-E_\Lambda^\Phi(\gamma_{X\setminus\Lambda} + \sum_{k=1}^{n} \varepsilon_{x_k})]\, \sigma(dx_1)\dots\sigma(dx_n)],$$

where

$$\tilde{Z}_\Lambda^{\sigma,\Phi}(\gamma) := 1 + \sum_{n=1}^{\infty} \frac{1}{n!} \int_\Lambda \dots \int_\Lambda \exp[-E_\Lambda^\Phi(\gamma_{X\setminus\Lambda} + \sum_{k=1}^{n} \varepsilon_{x_k})]\ \sigma(dx_1)\dots\sigma(dx_n)$$

$$= \exp(\sigma(\Lambda)) Z_\Lambda^{\sigma,\Phi}(\gamma).$$

(iii) Clearly, by properties (S 2), (S 3) a probability measure $\mu$ on $(\Gamma_X, \mathcal{B}(\Gamma_X))$ is a grand canonical Gibbs measure if and only if for all $\Lambda \in \mathcal{O}_c(X)$ and all $G \in \mathcal{B}_b(\Gamma_X)$

$$E_\mu[G \mid \mathcal{B}_{X\setminus\Lambda}(\Gamma_X)] \ = \ \Pi_\Lambda^{\sigma,\Phi} G \ \ \mu-\text{a.e.},$$

where for a sub-$\sigma$-algebra $\Sigma \subset \mathcal{B}(\Gamma_X)$, $E_\mu[\cdot \mid \Sigma]$ denotes the conditional expectation w.r.t. $\mu$ given $\Sigma$ and for $G \in \mathcal{B}_b(\Gamma_X)$ we set

$$\Pi_\Lambda^{\sigma,\Phi}(\gamma, G) \ := \ (\Pi_\Lambda^{\sigma,\Phi} G)(\gamma) \ = \ \int G(\gamma')\ \Pi_\Lambda^{\sigma,\Phi}(\gamma, d\gamma'). \tag{6.7}$$

(iv) Above as well as in the next section we may always replace $\mathcal{O}_c(X)$ by $\mathcal{B}_c(X)$ without any changes (in particular, $\mathcal{G}_{gc}(\sigma, \Phi)$ remains unchanged). We chose $\Lambda \in \mathcal{O}_c(X)$ only for simplicity.

## 6.2 Canonical Gibbs measures

Consider the situation of Subsection 6.1.

**Definition 6.3** *For $\Lambda \in \mathcal{O}_c(X)$ define for $\gamma \in \Gamma_X, \Delta \in \mathcal{B}(\Gamma_X)$*

$$\hat{\Pi}_\Lambda^{\sigma,\Phi}(\gamma, \Delta) \ := \ \begin{cases} \frac{\Pi_\Lambda^{\sigma,\Phi}(\gamma, \Delta\cap\{N_\Lambda=\gamma(\Lambda)\})}{\Pi_\Lambda^{\sigma,\Phi}(\gamma,\{N_\Lambda=\gamma(\Lambda)\})} & \text{if } \Pi_\Lambda^{\sigma,\Phi}(\gamma, \{N_\Lambda = \gamma(\Lambda)\}) > 0 \\ 0 & \text{otherwise} \end{cases} \tag{6.8}$$

*A probability measure $\mu$ on $(\Gamma_X, \mathcal{B}(\Gamma_X))$ is called a* canonical Gibbs measure with interaction potential $\Phi$ *if*

$$\mu\hat{\Pi}_\Lambda^{\sigma,\Phi} \ = \ \mu \ \textit{for all} \ \Lambda \in \mathcal{O}_c(X). \tag{6.9}$$

*Let $\mathcal{G}_c(\sigma, \Phi)$ denote the set of all such probability measures $\mu$.*

**Remark 6.4** For $\Lambda \in \mathcal{O}_c(X)$ let $\hat{\mathcal{B}}_{X\setminus\Lambda}(\Gamma_X)$ denote the $\sigma$ - algebra generated by $\mathcal{B}_{X\setminus\Lambda}(\Gamma)$ and $\sigma(N_\Lambda)$. It has been shown in [P79] that $(\hat{\Pi}_\Lambda^{\sigma,\Phi})_{\Lambda\in\mathcal{O}_c(X)}$ is indeed a $(\hat{\mathcal{B}}_{X\setminus\Lambda}(\Gamma_X))_{\Lambda\in\mathcal{O}_c(X)}$ - specification (i.e., satisfies (S 1) - (S 4) in Remark 6.2 (i) with $\hat{\Pi}_\Lambda^{\sigma,\Phi}$ replacing $\Pi_\Lambda^{\sigma,\Phi}$ and $\hat{\mathcal{B}}_{X\setminus\Lambda}(\Gamma_X)$ replacing $\mathcal{B}_{X\setminus\Lambda}(\Gamma_X)$ for $\Lambda \in \mathcal{O}_c(X)$). So, the above definition makes sense. In particular, for all $\Lambda \in \mathcal{O}_c(X)$ we have that

$$\mu(\{\gamma \in \Gamma_X \mid \Pi_\Lambda^{\sigma,\Phi}(\gamma, \{N_\Lambda = \gamma(\Lambda)\}) > 0\}) \ = \ 1 \tag{6.10}$$

for all $\mu \in \mathcal{G}_c(\sigma, \Phi)$. Clearly, also the analogue of Remark 6.2 (iii) holds.

It can be checked that

$$\mathcal{G}_{gc}(\sigma, \Phi) \ \subset \mathcal{G}_c(\sigma, \Phi) \tag{6.11}$$

(cf. [P79, Proposition 2.1]). The precise relation between the convex sets $\mathcal{G}_c(\sigma, \Phi)$ and $\mathcal{G}_{gc}(\sigma, \Phi)$ resp. between their extreme points $ex\mathcal{G}_c(\sigma, \Phi)$ and $ex\mathcal{G}_{gc}(\sigma, \Phi)$ is one of the major problems in Statistical Mechanics, known under the key word *equivalence of ensembles* (cf. [P79], [Ge79]).

**Remark 6.5** By a general result due to E.B. Dynkin and H. Föllmer (cf. [Dy78], [Fö75], and for a beautiful exposition [P76, Theorem 2.2]) we have $\mathrm{ex}\mathcal{G}_c(\sigma, \Phi) \neq \emptyset$, $\mathrm{ex}\mathcal{G}_{gc}(\sigma, \Phi) \neq \emptyset$ provided $\mathcal{G}_c(\sigma, \Phi) \neq \emptyset$, $\mathcal{G}_{gc}(\sigma, \Phi) \neq \emptyset$ respectively, and that any $\mu$ in $\mathcal{G}_c(\sigma, \Phi)$ resp. in $\mathcal{G}_{gc}(\sigma, \Phi)$ has an integral representation in terms of the respective extreme points (see Theorem 6.6 below for the special case $\Phi \equiv 0$).

In the sequel we shall mainly analyze $\mathcal{G}_c(\sigma, \Phi)$, but using known results on the above problem we shall deduce new results also for $\mathcal{G}_{gc}(\sigma, \Phi)$ (cf. Subsection 9.2). Concrete examples for the potential $\Phi$ will be treated in Subsection 7.1 below.

In fact the mixed Poisson measures from the previous sections are special canonical Gibbs measures with potential $\Phi \equiv 0$, i.e., in the *free case*. Let us set $\mathcal{G}_c(\sigma) := \mathcal{G}_c(\sigma, 0)$, $\mathcal{G}_{gc}(\sigma) := \mathcal{G}_{gc}(\sigma, 0)$. Then we have the following well known result (cf. [NZ77, Sect. 2] or [Ge79, Sect. 4.2], resp. [MKM78]).

**Theorem 6.6** *(i) Suppose $\sigma(X) = \infty$ and let $\mu$ be a probability measure on $(\Gamma_X, \mathcal{B}(\Gamma_X))$. Then $\mu \in \mathcal{G}_c(\sigma)$ if and only if there exists a probability measure $\lambda$ on $(\mathbb{R}_+, \mathcal{B}(\mathbb{R}_+))$ such that*

$$\mu \;=\; \int_{\mathbb{R}_+} \pi_{z\sigma}\; \lambda(dz)$$

*(i.e., if and only if $\mu$ is a mixed Poisson measure with mean proportional to $\sigma$).*

*(ii)* $\qquad ex\mathcal{G}_c(\sigma) = \{\pi_{z\sigma} \mid z \in \mathbb{R}_+\}$.

*(iii)* $\qquad ex\mathcal{G}_c(\sigma) = \cup_{z\in\mathbb{R}_+} ex\mathcal{G}_{gc}(z\sigma)$.

We note that since $\sigma$ is diffuse, the set $\mathfrak{M}^{\cdot\cdot}$ in [NZ77, Sect. 2] (resp. $M$ in [Ge79, Sect. 4.2]) can really be replaced by $\Gamma_X$ as done in Theorem 6.6 above.

## 6.3 Closability of the corresponding Dirichlet forms

In this subsection we want to prove the closability of $(\mathcal{E}_\mu^\Gamma, \mathcal{F}C_b^\infty)$ defined in (4.2) on $L^2(\Gamma_X;\mu)$ for all $\mu \in \mathcal{G}^1_{gc}(\sigma,\Phi)$ where

$$\mathcal{G}^1_{gc}(\sigma,\Phi) := \{\mu \in \mathcal{G}_{gc}(\sigma,\Phi) \mid \mu \text{ satisfies (4.1)}\} \tag{6.12}$$

(with $\mathcal{G}_{gc}(\sigma,\Phi)$ as defined in Subsection 6.1) under quite general assumptions on the potential $\Phi$. The presentation below is a self–contained and (w.r.t. technical details) completed version of [dSKR98, Sections 5, 6].

We start with giving a simple direct proof for the following identity for $\mu \in \mathcal{G}^1_{gc}(\sigma,\Phi)$ which was derived in [NZ79], [MMW79] from the well–known Mecke identity for $\pi_\sigma$ (cf. [Mec67, Satz 3.1]).

**Lemma 6.7** *Let $\mu \in \mathcal{G}^1_{gc}(\sigma,\Phi)$ and let $h : \Gamma_X \times X \to \mathbb{R}_+$ be $\mathcal{B}(\Gamma_X)\otimes\mathcal{B}(X)$–measurable. Then*

$$\int_{\Gamma_X}\int_X h(\gamma,x)\;\gamma(dx)\mu(d\gamma) = \int_{\Gamma_X}\int_X h(\gamma+\varepsilon_x,x)e^{-E^\Phi_{\{x\}}(\gamma+\varepsilon_x)}\sigma(dx)\mu(d\gamma), \tag{6.13}$$

*where*

$$E^\Phi_{\{x\}}(\gamma+\varepsilon_x) := \begin{cases} \displaystyle\sum_{\{x\}\subset\gamma'\subset\gamma\cup\{x\}} \Phi(\gamma') & \textit{if} \displaystyle\sum_{\{x\}\subset\gamma'\subset\gamma\cup\{x\}} |\Phi(\gamma')| < \infty \\ +\infty & \textit{else} \end{cases} \tag{6.14}$$

**Proof.** Because of standard monotone class arguments it suffices to prove (6.13) for

$$h(\gamma, x) := e^{\langle g,\gamma\rangle}\ f(x), \quad \gamma \in \Gamma_X, x \in X,$$

where $f, g \in C_0^\infty(X)$; $(-g), f \geq 0$. Let $\Lambda \in \mathcal{O}_c(X)$ such that $\operatorname{supp} f, \operatorname{supp} g \subset \Lambda$. By Remark 6.2 (ii) the left hand side of (6.13) is equal to

$$\begin{aligned}
& \int e^{\langle g,\gamma\rangle}\ \langle f,\gamma\rangle\ \mu(d\gamma) \\
= & \int \Pi_\Lambda^{\sigma,\Phi}(e^{\langle g,\cdot\rangle}\langle f,\cdot\rangle)\, d\mu \\
= & \int_{\Gamma_X} (\widetilde{Z}_\Lambda^{\sigma,\Phi}(\gamma))^{-1} \sum_{n=1}^\infty \frac{1}{n!} \int_{\Lambda^n} e^{\sum_{k=1}^n g(x_k)} \sum_{i=1}^n f(x_i) \exp[-E_\Lambda^\Phi(\gamma_{X\setminus\Lambda} + \sum_{k=1}^n \varepsilon_{x_k})] \\
& \qquad \sigma(dx_1)\ldots\sigma(dx_n)\,\mu(d\gamma) \\
= & \int_\Lambda e^{g(x)}\ f(x) \int_{\Gamma_X} (\widetilde{Z}_\Lambda^{\sigma,\Phi}(\gamma))^{-1} \sum_{n=1}^\infty \frac{1}{(n-1)!} \int_{\Lambda^{n-1}} e^{\sum_{k=1}^{n-1} g(x_k)}\ \cdot \\
& \qquad \cdot \exp[-E_\Lambda^\Phi(\gamma_{X\setminus\Lambda} + \sum_{k=1}^{n-1} \varepsilon_{x_k} + \varepsilon_x)]\ \sigma(dx_1)\ldots\sigma(dx_{n-1})\,\mu(d\gamma)\,\sigma(dx).
\end{aligned}$$

It is easy to see that

$$E_\Lambda^\Phi(\gamma + \varepsilon_x) = E_\Lambda^\Phi(\gamma) + E_{\{x\}}^\Phi(\gamma + \varepsilon_x) \text{ for all } \gamma \in \Gamma_X, x \in \Lambda \text{ such that } x \notin \gamma, \tag{6.15}$$

and that

$$(\gamma, x) \mapsto E_{\{x\}}^\Phi(\gamma + \varepsilon_x) \text{ is } \mathcal{B}(\Gamma_X) \times \mathcal{B}(X)\text{–measurable.} \tag{6.16}$$

Consequently, the left hand side of (6.13) is equal to

$$\begin{aligned}
& \int_\Lambda e^{g(x)}\ f(x) \int_{\Gamma_X} (\widetilde{Z}_\Lambda^{\sigma,\Phi}(\gamma))^{-1} \sum_{n=0}^\infty \frac{1}{n!} \int_{\Lambda^n} e^{\sum_{k=1}^n g(x_k)} \\
& \cdot \exp[-E_{\{x\}}^\Phi(\gamma_{X\setminus\Lambda} + \varepsilon_x + \sum_{k=1}^n \varepsilon_x)] \cdot \exp[-E_\Lambda^\Phi(\gamma_{X\setminus\Lambda} + \sum_{k=1}^n \varepsilon_{x_k})] \\
& \quad \sigma(dx_1)\ldots\sigma(dx_n)\ \mu(d\gamma)\ \sigma(dx) \\
= & \int_\Lambda e^{g(x)} f(x) \int_{\Gamma_X} \Pi_\Lambda^{\sigma,\Phi}\left(e^{\langle g,\cdot\rangle} e^{-E_{\{x\}}^\Phi(\cdot\,+\varepsilon_x)}\right)(\gamma)\,\mu(d\gamma)\sigma(dx) \\
= & \int_\Lambda e^{g(x)} f(x) \int_{\Gamma_X} e^{\langle g,\gamma\rangle} e^{-E_{\{x\}}^\Phi(\gamma+\varepsilon_x)} \mu(d\gamma)\ \sigma(dx)
\end{aligned}$$

which by (6.16) and Fubini's theorem is equal to the right hand side of (6.13). $\square$

**Remark 6.8** (i) We note that if $x \notin \gamma$ then the definition of $E^{\Phi}_{\{x\}}(\gamma + \varepsilon_x)$ in (6.14) coincides with (6.5) for $\Lambda = \{x\}$. On the other hand if

$$M_x := \{\gamma \in \Gamma_X \mid x \in \gamma\},$$

then for $\Lambda \in \mathcal{O}_c(X)$ such that $x \in \Lambda$ and every $\mu \in \mathcal{G}^1_{gc}(\sigma, \Phi)$.

$$\begin{aligned}
& \mu(M_x) \\
= & \int \Pi^{\sigma,\Phi}_{\Lambda}(\gamma, M_x)\ \mu(d\gamma) \\
= & \int_{\Gamma_X} (\widetilde{Z}^{\sigma,\Phi}_{\Lambda}(\gamma))^{-1} \sum_{n=1}^{\infty} \frac{1}{n!} \int_{\Lambda^n} 1_{M_x}(\gamma_{X\setminus\Lambda} + \sum_{k=1}^{n} \varepsilon_{x_k}) \exp[-E^{\Phi}_{\Lambda}(\gamma_{X\setminus\Lambda} + \sum_{k=1}^{n} \varepsilon_{x_k})] \\
& \qquad \sigma(dx_1)\ldots\sigma(dx_n)\ \mu(d\gamma) \\
= & \int_{\Gamma_X} (\widetilde{Z}^{\sigma,\Phi}_{\Lambda}(\gamma))^{-1} \sum_{n=1}^{\infty} \frac{1}{n!} \int_{\tilde{\Lambda}^n} \sum_{i=1}^{n} 1_{\{x\}}(x_i) \exp[-E^{\Phi}_{\Lambda}(\gamma_{X\setminus\Lambda} + \sum_{k=1}^{n} \varepsilon_{x_k})] \\
& \qquad \sigma(dx_1)\ldots\sigma(dx_n)\ \mu(d\gamma) \\
= & \ 0\ ,
\end{aligned} \tag{6.17}$$

where $\tilde{\Lambda}^n$ is as defined in (2.2). Hence (6.14) and (6.5) coincide for $\mu$–a.e. $\gamma \in \Gamma_X$.

(ii) We emphasize that Lemma 6.7, of course, holds for arbitrary positive Radon measures $\sigma$ on $(X, \mathcal{B}(X))$. We never used our specific form $\sigma = \rho \cdot m$ in the above proof.

Now we can prove the following representation of $(\mathcal{E}^{\Gamma}_{\mu}, \mathcal{F}C^{\infty,\mu}_b)$ (cf. [dSKR98, Theorem 5.1]).

**Theorem 6.9** *Let $\mu \in \mathcal{G}^1_{gc}(\sigma, \Phi)$ and let $F, G \in \mathcal{F}C^{\infty}_b$. Then*

$$\begin{aligned}
& \int_{\Gamma_X} \langle \nabla^{\Gamma} F, \nabla^{\Gamma} G \rangle_{T\Gamma_X}\ d\mu \\
= & \int_{\Gamma_X} \int_X \big\langle \nabla^X (F(\gamma + \varepsilon_\cdot) - F(\gamma)), \nabla^X (G(\gamma + \varepsilon_\cdot) - G(\gamma)) \big\rangle_{T\Gamma_X}(x) \\
& \qquad \sigma(dx)\ \mu(d\gamma).
\end{aligned} \tag{6.18}$$

*In particular,* (4.3) *holds, so $(\mathcal{E}^{\Gamma}_{\mu}, \mathcal{F}C^{\infty,\mu}_b)$ is a well-defined positive definite symmetric bilinear form on $L^2(\Gamma_X; \mu)$.*

**Proof.** By (2.8) and Lemma 6.7, if $F = g_F(\langle f_1, \cdot\rangle, \dots, \langle f_N, \cdot\rangle)$, $G = g_G(\langle g_1, \cdot\rangle, \dots, \langle g_M, \cdot\rangle)$, where without loss of generality we may assume $N = M$, we have

$$\begin{aligned}
&\int_{\Gamma_X} \langle \nabla^\Gamma F, \nabla^\Gamma G\rangle_{T\Gamma_X} d\mu \\
= &\int_{\Gamma_X} \int_X \sum_{i,j=1}^N \partial_i g_F(\langle f_1, \gamma\rangle, \dots, \langle f_N, \gamma\rangle)\, \partial_j g_G(\langle(g_1, \gamma\rangle, \dots, \langle g_N, \gamma\rangle) \\
&\qquad \langle \nabla^X f_i, \nabla^X g_j\rangle_{TX}(x)\ \gamma(dx)\ \mu(d\gamma) \\
= &\int_{\Gamma_X} \int_X \sum_{i,j=1}^N \partial_i g_F(\langle f_1, \gamma\rangle + f_1(x), \dots, \langle f_N, \gamma\rangle + f_N(x)) \\
&\qquad \partial_j g_G(\langle g_1, \gamma\rangle + g_1(x), \dots, \langle g_N, \gamma\rangle + g_N(x)) \\
&\qquad \langle \nabla^X f_i, \nabla^X g_j\rangle_{TX}(x)\ \sigma(dx)\ \mu(d\gamma).
\end{aligned}$$

The latter clearly is equal to the right hand side of (6.18).
If now $F = G$ $\mu$–a.e., then by Lemma 6.7 for all $\Lambda \in \mathcal{O}_c(X)$ and $F_0 := F - G$

$$\begin{aligned}
0 &= \int_{\Gamma_X} \gamma(\Lambda)\ F_0(\gamma)\mu(d\gamma) = \int_{\Gamma_X} \int_X 1_\Lambda(x) F_0(\gamma)\ \gamma(dx)\mu(d\gamma) \\
&= \int_{\Gamma_X} \int_X 1_\Lambda(x) F_0(\gamma + \varepsilon_x)\sigma(dx)\ \mu(d\gamma).
\end{aligned}$$

Since $\Lambda \in \mathcal{O}_c(X)$ was arbitrary, it follows that for $\mu$–a.e. $\gamma \in \Gamma_X$

$$F_0(\gamma + \varepsilon_x) = 0 \text{ for } \sigma\text{–a.e. } x \in X.$$

This implies that the right hand side of (6.18) applied to $F := G := F_0$ is zero, hence by (6.18)

$$\nabla^\Gamma G - \nabla^\Gamma F = \nabla^\Gamma F_0 = 0 \quad \mu\text{–a.e.}$$

and (4.3) is proven. □

**Remark 6.10** The reader might wonder why we took $F(\gamma+\varepsilon.) - F(\gamma)$ rather than only $F(\gamma + \varepsilon.)$ in (6.18). The reason is that for all $F \in \mathcal{F}C_b^\infty$

$$x \mapsto F(\gamma + \varepsilon_x) - F(\gamma) \in C_0^\infty(X)$$

for all $\gamma \in \Gamma_X$, whereas $x \mapsto F(\gamma + \varepsilon_x)$ is merely in $C_b^\infty(X)$.

Now we are going to prove the main result of this subsection, i.e., the said closability. This was first proved in [dSKR98, Sect. 6]. Because this result is of fundamental importance for the analysis on configurations spaces we repeat the complete proof here:

We define new intensity measures on $X$ by $\sigma_\gamma := \rho_\gamma \cdot m$, where

$$\rho_\gamma(x) := e^{-E^{\Phi}_{\{x\}}(\gamma+\varepsilon_x)}\rho(x),\ x \in X, \gamma \in \Gamma_X\ . \tag{6.19}$$

It was shown in [AR90, Theorem 5.3] (in the case $X = \mathbb{R}^d$) that the components of the Dirichlet form $(\mathcal{E}^X_{\sigma_\gamma}, \mathcal{D}^{\sigma_\gamma})$ corresponding to the measure $\sigma_\gamma$ are closable on $L^2(\mathbb{R}^d; \sigma_\gamma)$ if and only if $\sigma_\gamma$ is absolutely continuous with respect to Lebesgue measure on $\mathbb{R}^d$ and the Radon-Nikodym derivative satisfies some condition (see (6.20) below for details). This result allows us to prove the closability of $(\mathcal{E}^\Gamma_\mu, \mathcal{F}C^{\infty,\mu}_b)$ on $L^2(\Gamma_X; \mu)$. Let us first recall the above mentioned result.

**Theorem 6.11** (cf. [AR90, Theorem 5.3]) *Let $\nu$ by a probability measure on $(\mathbb{R}^d, \mathcal{B}(\mathbb{R}^d))$, $d \in \mathbb{N}$, and let $\mathcal{D}^\nu$ denote the $\nu$-classes determined by $\mathcal{D}$. Then the forms $(\mathcal{E}_{\nu,i}, \mathcal{D}^\nu)$ defined by*

$$\mathcal{E}_{\nu,i}(u,v) := \int_{\mathbb{R}^d} \frac{\partial u}{\partial x_i}\frac{\partial v}{\partial x_i} d\nu,\ u, v \in \mathcal{D},$$

*are well-defined and closable on $L^2(\mathbb{R}^d; \nu)$ for $1 \le i \le d$ if and only if $\nu$ is absolutely continuous with respect to Lebesgue measure $\lambda^d$ on $\mathbb{R}^d$, and the Radon-Nikodym derivative $\rho = d\nu/d\lambda^d$ satisfies the condition:*

$$\begin{aligned}&\text{for any } 1 \le i \le d \text{ and } \lambda^{d-1}\text{–a.e. } x \in \left\{ y \in \mathbb{R}^{d-1} |\ \int_{\mathbb{R}} \rho^{(i)}_y(s)\lambda^1(ds) > 0 \right\},\\ &\rho^{(i)}_x = 0\ \lambda^1\text{–a.e. on } \mathbb{R}\backslash R(\rho^{(i)}_x),\ \text{where } \rho^{(i)}_x(s) := \rho(x_1, \ldots, x_{i-1}, s, x_i, \ldots, x_d),\\ &s \in \mathbb{R},\ \text{if } x = (x_1, \ldots, x_{d-1}) \in \mathbb{R}^{d-1}\end{aligned} \tag{6.20}$$

*and where*

$$R(\rho^{(i)}_x) := \left\{ t \in \mathbb{R}\ |\ \int_{t-\varepsilon}^{t+\varepsilon} \frac{1}{\rho^{(i)}_x(s)}\lambda^1(ds) < \infty \textit{ for some } \varepsilon > 0 \right\}. \tag{6.21}$$

There is an obvious generalization of Theorem 6.11 to the case where a Riemannian manifold $X$ is replacing $\mathbb{R}^d$, to be formulated in terms of local charts. Since here we are only interested in the "if part" of Theorem 6.11, we now recall a slightly weaker sufficient condition for closability in the general case where $X$ is a manifold as before.

**Theorem 6.12** *Suppose $\sigma_1 = \rho_1 \cdot m$, where $\rho_1 : X \to \mathbb{R}_+$ is $\mathcal{B}(X)$-measurable such that*

$$\rho_1 = 0 \ m\text{−a.e. on } X \setminus \Big\{ x \in X \mid \int_{\Lambda_x} \frac{1}{\rho_1} dm < \infty \textit{ for some open neighbourhood } \Lambda_x \textit{ of } x \Big\}. \tag{6.22}$$

*Then $(\mathcal{E}^X_{\sigma_1}, \mathcal{D}^{\sigma_1})$ defined by*

$$\mathcal{E}^X_{\sigma_1}(u,v) := \int_X \langle \nabla^X u, \nabla^X v \rangle_{TX} \ d\sigma_1; \ u, v \in \mathcal{D},$$

*is closable on $L^2(X;\sigma_1)$.*

The proof is a straightforward adaptation of the line of arguments in [MR92, Chap. II, Subsection 2a] (see also [ABR89, Theorem 4.2] for details). We emphasize that (6.22) (hence (6.20)) e.g. always holds, if $\rho_1$ is lower semi-continuous, and that (6.20) holds if e.g. $\rho$ is weakly differentiable. Neither $\nu$ in Theorem 6.11 nor $\sigma_1$ in Theorem 6.12 is required to have full support, so e.g. $\rho_1$ is not necessarily strictly positive $m$-a.e. on $X$.

**Theorem 6.13** *Let $\mu \in \mathcal{G}^1_{gc}(\sigma, \Phi)$. Suppose that for $\mu$-a.e. $\gamma \in \Gamma_X$ the function $\rho_\gamma$ defined in (6.19) satisfies (6.22) (resp. (6.20) in case $X = \mathbb{R}^d$). Then the form $(\mathcal{E}^\Gamma_\mu, \mathcal{F}C^{\infty,\mu}_b)$ is closable on $L^2(\Gamma_X;\mu)$.*

**Proof.** Let $(F_n)_{n\in\mathbb{N}}$ be a sequence in $\mathcal{F}C^{\infty,\mu}_b$ such that $F_n \to 0$ in $L^2(\Gamma_X;\mu)$ as $n \to \infty$ and

$$\mathcal{E}^\Gamma_\mu(F_n - F_m, F_n - F_m) \underset{n,m\to\infty}{\longrightarrow} 0. \tag{6.23}$$

We have to show that

$$\mathcal{E}^\Gamma_\mu(F_{n_k}, F_{n_k}) \underset{k\to\infty}{\longrightarrow} 0 \tag{6.24}$$

for some subsequence $(n_k)_{k\in\mathbb{N}}$. Let $(n_k)_{k\in\mathbb{N}}$ be a subsequence such that

$$\left( \int_{\Gamma_X} F^2_{n_k} d\mu \right)^{1/2} + \mathcal{E}^\Gamma_\mu(F_{n_{k+1}} - F_{n_k}, F_{n_{k+1}} - F_{n_k})^{1/2} < \frac{1}{2^k} \text{ for all } k \in \mathbb{N}.$$

Then

$$\begin{aligned}
\infty &> \sum_{k=1}^{\infty} \mathcal{E}_{\mu}^{\Gamma}(F_{n_{k+1}} - F_{n_k}, F_{n_{k+1}} - F_{n_k})^{1/2} \\
&\geq \sum_{k=1}^{\infty} \int_{\Gamma_X} \Big( \int_X |\nabla^X((F_{n_{k+1}} - F_{n_k})(\gamma + \varepsilon.) - (F_{n_{k+1}} - F_{n_k})(\gamma))|_{TX}^2(x) \\
&\qquad\qquad e^{-E_{\{x\}}^{\Phi}(\gamma+\varepsilon_x)}\,\sigma(dx)\Big)^{1/2} \mu(d\gamma) \\
&= \int_{\Gamma_X} \sum_{k=1}^{\infty} \Big( \int_X |\nabla^X((F_{n_{k+1}} - F_{n_k})(\gamma + \varepsilon.) - (F_{n_{k+1}} - F_{n_k})(\gamma))|_{TX}^2(x) \\
&\qquad\qquad \sigma_\gamma(dx)\Big)^{1/2} \mu(d\gamma),
\end{aligned}$$

where we used Theorem 6.9 and (6.19). From the last expression we obtain that

$$\begin{aligned}
&\sum_{k=1}^{\infty} \mathcal{E}_{\sigma_\gamma}^{X}(u_{n_{k+1}}^{(\gamma)} - u_{n_k}^{(\gamma)}, u_{n_{k+1}}^{(\gamma)} - u_{n_k}^{(\gamma)})^{1/2} \\
&= \sum_{k=1}^{\infty} \Big( \int_X |\nabla^X((F_{n_{k+1}} - F_{n_k})(\gamma + \varepsilon.) - (F_{n_{k+1}} - F_{n_k})(\gamma))|_{TX}^2(x) \\
&\qquad\qquad \cdot e^{-E_{\{x\}}^{\Phi}(\gamma+\varepsilon_x)}\,\sigma(dx)\Big)^{1/2} \\
&< \infty \qquad \text{for } \mu\text{-a.e. } \gamma \in \Gamma_X,
\end{aligned} \tag{6.25}$$

where for $k \in \mathbb{N}$, $\gamma \in \Gamma_X$,

$$u_{n_k}^{(\gamma)}(x) := F_{n_k}(\gamma + \varepsilon_x) - F_{n_k}(\gamma),\ x \in X.$$

We note that $u_{n_k}^{(\gamma)} \in \mathcal{D}$. (6.25) implies that for $\mu$-a.e. $\gamma \in \Gamma_X$

$$\mathcal{E}_{\sigma_\gamma}^{X}(u_{n_k}^{(\gamma)} - u_{n_l}^{(\gamma)}, u_{n_k}^{(\gamma)} - u_{n_l}^{(\gamma)}) \underset{k,l\to\infty}{\longrightarrow} 0. \tag{6.26}$$

Let $\Lambda \subset \mathcal{O}_c(X)$.

**Claim 1**: For $\mu$-a.e. $\gamma \in \Gamma_X$

$$\int_X (u_{n_k}^{(\gamma)}(x))^2\, 1_\Lambda(x)\sigma_\gamma(dx) \underset{k\to\infty}{\longrightarrow} 0.$$

To prove Claim 1 we first note that for $\mu$-a.e. $\gamma \in \Gamma_X$

$$\sigma_\gamma(\Lambda) < \infty,$$

as follows immediately from (6.13) (taking $h(\gamma, x) := 1_\Lambda(x)$ for x$\in X$, $\gamma \in \Gamma_X$), since $\mu \in \mathcal{G}^1_{gc}(\sigma, \Phi)$. Therefore, for $\mu$-a.e. $\gamma \in \Gamma_X$

$$\int_X F^2_{n_k}(\gamma)\, 1_\Lambda(x)\sigma_\gamma(dx) = F^2_{n_k}(\gamma)\sigma_\gamma(\Lambda) \underset{k\to\infty}{\longrightarrow} 0. \tag{6.27}$$

Furthermore, by (6.13)

$$\begin{aligned} & \int_{\Gamma_X}\int_X F^2_{n_k}(\gamma+\varepsilon_x)\, 1_\Lambda(x)\sigma_\gamma(dx)(1+\gamma(\Lambda))^{-1}\mu(d\gamma) \\ = & \int_{\Gamma_X} F^2_{n_k}(\gamma)\int_X \frac{1_\Lambda(x)}{1+\gamma(\Lambda)-1_\Lambda(x)}\gamma(dx)\mu(d\gamma) \\ \le & \int_{\Gamma_X} F^2_{n_k}(\gamma)\mu(d\gamma) < \frac{1}{2^k}, \end{aligned}$$

because the integral w.r.t. $\gamma$ is dominated by 1 for all $\gamma \in \Gamma_X$. Hence

$$\begin{aligned} \infty &> \sum_{k=1}^\infty \left(\int_{\Gamma_X}\int_X F^2_{n_k}(\gamma+\varepsilon_x)\, 1_\Lambda(x)\sigma_\gamma(dx)(1+\gamma(\Lambda))^{-1}\mu(d\gamma)\right)^{1/2} \\ &\ge \int_{\Gamma_X}\sum_{k=1}^\infty\left(\int_X F^2_{n_k}(\gamma+\varepsilon_x)\, 1_\Lambda(x)\sigma_\gamma(dx)\right)^{1/2}(1+\gamma(\Lambda))^{-1}\mu(d\gamma). \end{aligned}$$

Therefore, for $\mu$-a.e. $\gamma \in \Gamma_X$

$$\int_X F^2_{n_k}(\gamma+\varepsilon_x)\, 1_\Lambda(x)\sigma_\gamma(dx) \underset{k\to\infty}{\longrightarrow} 0. \tag{6.28}$$

Now Claim 1 follows by (6.27) and (6.28).

**Claim 2**: For $\mu$-a.e. $\gamma \in \Gamma_X$

$$|\nabla^X u^{(\gamma)}_{n_k}|_{TX} \underset{k\to\infty}{\longrightarrow} 0 \quad \sigma_\gamma\text{−a.e.}$$

To prove Claim 2 we first note that clearly (6.25) implies that for $\mu$-a.e. $\gamma \in \Gamma_X$

$$\mathcal{E}^X_{1_\Lambda\sigma_\gamma}(u^{(\gamma)}_{n_k} - u^{(\gamma)}_{n_l}, u^{(\gamma)}_{n_k} - u^{(\gamma)}_{n_l}) \underset{k,l\to\infty}{\longrightarrow} 0. \tag{6.29}$$

Hence we can apply Theorem 6.12 (resp. 6.11) to $\rho_1 := 1_\Lambda \rho_\gamma$ and conclude by Claim 1 and (6.29) that for $\mu$-a.e. $\gamma \in \Gamma_X$

$$\mathcal{E}^X_{1_\Lambda \sigma_\gamma}(u^{(\gamma)}_{n_k}, u^{(\gamma)}_{n_k}) \underset{k\to\infty}{\longrightarrow} 0,$$

hence by (6.25)

$$1_\Lambda |\nabla^X u^{(\gamma)}_{n_k}|_{TX} \underset{k\to\infty}{\longrightarrow} 0 \quad \sigma_\gamma\text{−a.e.}$$

Since $\Lambda$ was arbitrary, Claim 2 is proven.

¿From Claim 2 we now easily deduce (6.24) by (6.18) and Fatou's Lemma as follows:

$$\begin{aligned} \mathcal{E}^\Gamma_\mu(F_{n_k}, F_{n_k}) &\le \int_{\Gamma_X} \liminf_{l\to\infty} \int_X |\nabla^X (u^{(\gamma)}_{n_k} - u^{(\gamma)}_{n_l})|^2_{TX}\ \sigma_\gamma(dx)\, \mu(d\gamma) \\ &\le \liminf_{l\to\infty} \mathcal{E}^\Gamma_\mu(F_{n_k} - F_{n_l}, F_{n_k} - F_{n_l}), \end{aligned}$$

which by (6.23) can be made arbitrarily small for $k$ large enough. □

**Remark 6.14** (i) The above method to prove closability of pre-Dirichlet forms on configuration spaces $\Gamma_X$ extends immediately to the case where $X$ is replaced by an infinite dimensional "manifold" such as the loop space (cf. [MR97]).

(ii) We also emphasize that all results in Subsection 4.2 now apply to every $\mu \in \mathcal{G}^1_{gc}(\sigma, \Phi)$ that satisfies the condition of Theorem 6.13.

# 7 Integration by parts characterization of canonical Gibbs measures

In this section we shall first recall the definition of a class of grand canonical Gibbs measures, so–called Ruelle measures. They serve as examples for the measures in $\mathcal{G}^1_{gc}(\sigma, \Phi)$ studied in the previous section, with explicit (pair) potentials $\Phi$. Subsequently, we shall present the analogue of Theorem 2.2 for canonical Gibbs measures which was proven in [AKR97b, Subsection 4.3]. In this section we consider for simplicity the case $X = \mathbb{R}^d$ with the Euclidean metric and $\sigma = z \cdot m$, $z \in ]0, \infty[$ (with $m$ = Lebesgue measure). We also simplify notations by setting $\Gamma := \Gamma_{\mathbb{R}^d}$, $\pi^z := \pi_{z\cdot m}$, $\nabla := \nabla^{\mathbb{R}^d}$. We follow the presentation of [AKR97b], but omit the technically rather involved proofs and, for simplicity, only consider the "finite range case" in Subsection 7.2.

## 7.1 Ruelle measures

We will now describe (following essentially [Ru70]) a class of grand canonical Gibbs measures which appears in Classical Statistical Mechanics of continuous systems, see also [Do70], [KY93], [Ru69].

A *pair potential* is a $\mathcal{B}(\mathbb{R}^d)$–measurable function $\phi : \mathbb{R}^d \to \mathbb{R} \cup \{+\infty\}$ such that $\phi(-x) = \phi(x)$. Any pair potential $\phi$ defines a potential $\Phi = \Phi_\phi$ as follows: we set $\Phi(\gamma) := 0,\ |\gamma| \neq 2$ and $\Phi(\gamma) := \phi(x-y)$ for $\gamma = \{x,y\} \subset \mathbb{R}^d$. It is useful to rewrite the conditional energy $E_\Lambda^\Phi$ (see (6.5)) in the following form for $\Lambda \in \mathcal{B}_c(\mathbb{R}^d)$ :

$$E_\Lambda^\Phi(\gamma) = E_\Lambda^\Phi(\gamma_\Lambda) \ + \ W(\gamma_\Lambda|\gamma_{\Lambda^c}), \tag{7.1}$$

where the term

$$W(\gamma_\Lambda|\gamma_{\Lambda^c}) \ := \sum_{x\in\gamma_\Lambda, y\in\gamma_{\Lambda^c}} \phi(x-y) \tag{7.2}$$

describes the *interaction energy* between $\gamma_\Lambda$ and $\gamma_{\Lambda^c}$. (Here as usual $\Lambda^c := \mathbb{R}^d \setminus \Lambda$). Analogously, we define $W(\gamma'|\gamma'')$ when $\gamma', \gamma''$ are located in disjoint regions.

For every $r = (r^1, \ldots, r^d) \in \mathbb{Z}^d$ we define a cube

$$Q_r \ := \ \{x \in \mathbb{R}^d \mid r^i - \frac{1}{2} \leq x^i < r^i + \frac{1}{2}\}.$$

These cubes form a partition of $\mathbb{R}^d$. For any $\gamma \in \Gamma$ we set $\gamma_r := \gamma_{Q_r} = \gamma \cap Q_r,\ r \in \mathbb{Z}^d$. For $N \in \mathbb{N}$ let $\Lambda_N$ be the cube with side length $2N-1$ centered at the origin in $\mathbb{R}^d$; $\Lambda_N$ is then a union of $(2N-1)^d$ unit cubes of the form $Q_r$.

Now we are able to formulate conditions on the interaction.

(S) (*Stability*) There exists $B \geq 0$ such that for any $\Lambda \in \mathcal{B}_c(\mathbb{R}^d)$ and for all $\gamma \in \Gamma_\Lambda$

$$E_\Lambda^\Phi(\gamma) \ \geq \ -B|\gamma|.$$

A condition stronger than stability is the following.

(SS) (*Superstability*) There exist $A > 0,\ B \geq 0$ such that if $\gamma \in \Gamma_{\Lambda_N}$ for some $N$, then

$$E_{\Lambda_N}^\Phi(\gamma) \ \geq \ \sum_{r\in\mathbb{Z}^d} [A|\gamma_r|^2 \ - \ B|\gamma_r|].$$

(LR) (*Lower regularity*) There exists a decreasing positive function $a : \mathbb{N} \to \mathbb{R}_+$ such that

$$\sum_{r \in \mathbb{Z}^d} a(\|r\|) \ < \ \infty$$

and for any $\Lambda', \Lambda''$ which are each finite unions of unit cubes of the form $Q_r$ and disjoint, with $\gamma' \in \Gamma_{\Lambda'}$, $\gamma'' \in \Gamma_{\Lambda''}$,

$$W(\gamma'|\gamma'') \geq \ - \sum_{r',r'' \in \mathbb{Z}^d} a(\|r' - r''\|) \, |\gamma'_{r'}| \cdot |\gamma''_{r''}|.$$

Here $\| \cdot \|$ denotes the maximum norm on $\mathbb{R}^d$

We say that the pair potential $\phi$ is *stable* (resp. *superstable*, *lower regular*) if (S) (resp (SS), (LR)) is satisfied. In [Ru70] there are given many criteria for stability, superstability and lower regularity of a pair potential. For example, if $\phi_1 \geq 0$ is continuous and $\phi_1(0) > 0$, then $\phi_1$ is superstable. If $\phi_2$ is stable than $\phi = \phi_1 + \phi_2$ is superstable etc. Let us mention also the following Dobrushin-Fisher-Ruelle criterium.

**Proposition 7.1** *Let* $0 < d_1 < d_2 < +\infty$ *and let*

$$s_1 : [0, d_1] \to \mathbb{R} \cup \{+\infty\}, \ \ s_2 : [d_2, +\infty[ \to \mathbb{R}$$

*be positive, decreasing and such that*

$$\int_0^{d_1} t^{d-1} s_1(t) dt = +\infty, \ \int_{d_2}^{\infty} t^{d-1} s_2(t) dt < +\infty.$$

*If the pair potential* $\phi$ *is bounded below and satisfies*

$$\phi(x) \ \geq \ s_1(\|x\|) \ \textit{for} \ \|x\| \leq d_1,$$

$$|\phi(x)| \ \leq \ s_2(\|x\|) \ \textit{for} \ \|x\| \geq d_2 \ ,$$

*then* $\phi$ *is superstable and lower regular.*

**Definition 7.2** *A probability measure* $\mu$ *on* $(\Gamma, \mathcal{B}(\Gamma))$ *is called* tempered *if* $\mu$ *is supported by*

$$S_\infty := \cup_{n=1}^{\infty} S_n,$$

*where*

$$S_n := \{\gamma \in \Gamma \mid \forall N \in \mathbb{N} \sum_{r \in \Lambda_N \cap \mathbb{Z}^d} |\gamma_r|^2 \leq n^2 |\Lambda_N \cap \mathbb{Z}^d| \ \}.$$

By $\mathcal{G}^t_{gc}(z,\phi) \subset \mathcal{G}_{gc}(z,\phi) := \mathcal{G}_{gc}(zm,\Phi_\phi)$ (cf. Subsection 6.1) we denote the set of all tempered grand canonical Gibbs measures (*Ruelle measures* for short).

**Remark 7.3** Let $z \in ]0,\infty[$. Due to [Ru70, Sect. 5]

$$\mathcal{G}^t_{gc}(z,\phi) \neq \emptyset$$

and for all $\mu \in \mathcal{G}^t_{gc}(z,\phi)$ and all $q \in [1,\infty[$

$$\int \gamma(K)^q \, \mu(d\gamma) < \infty \quad \text{for all compact } K \subset \mathbb{R}^d$$

if $\phi$ is a super-stable, lower regular potential and satisfies the following condition:

$$\text{(I)} \;\; (\textit{Integrability}) \qquad \int_{\mathbb{R}^d} |1 - e^{-\phi(x)}| \; m(dx) \; < \; +\infty \; .$$

For a given probability measure $\mu$ on $(\Gamma, \mathcal{B}(\Gamma))$ and any $\Lambda \in \mathcal{B}_c(\mathbb{R}^d)$ we can introduce a probability measure $\mu_\Lambda$ on $(\Gamma_\Lambda, \mathcal{B}(\Gamma_\Lambda))$ as $\mu_\Lambda = \mu \circ p_\Lambda^{-1}$, cf. (2.17). Since there is a decomposition of $\Gamma_\Lambda$ analogous to (2.13) for all $\Lambda \in \mathcal{B}_c(\mathbb{R}^d)$ (and not only for $\Lambda \in \mathcal{O}_c(\mathbb{R}^d)$), we have an induced decomposition of $\mu_\Lambda$: $\mu_\Lambda = \sum_{n=0}^\infty \mu_\Lambda^{(n)}$. In the case of a Ruelle measure $\mu$ which corresponds to a potential $\phi$ with properties (SS), (LR) and (I) we have the representations

$$\mu_\Lambda^{(n)}(d\gamma_n) = \frac{1}{n!} \tau^n_{\Lambda,\mu}(\gamma_n) m_{\Lambda,n}(d\gamma_n), \; \gamma_n \in \Gamma_\Lambda^{(n)}, \; n \geq 0,$$

where $\tau^n_{\Lambda,\mu}(\cdot)$ are positive measurable functions on $\Gamma_\Lambda^{(n)}$. They can be considered as positive Lebesgue integrable functions $\tau^n_{\Lambda,\mu}(x_1,\ldots,x_n)$ on $\Lambda^n$ invariant w.r.t. all permutations from the group $S_n$. These functions are called the *system of density distributions.* We refer to [Ru70] for details.

The following theorem gives necessary properties of the so-called correlation functions associated to Ruelle measures, for the proof see [Ru70].

**Theorem 7.4** *Let $\phi$ be a pair potential satisfying conditions (SS), (LR) and (I). For a given $\mu \in \mathcal{G}^t_{gc}(\phi,z)$ let $\{\tau^n_{\Lambda,\mu} | n \geq 0\}$ be the associated system of density distributions. We define correlation functions by*

$$\varrho_\mu(x_1,\ldots,x_n) = \sum_{k=0}^\infty \frac{1}{k!} \int_{\Lambda^n} \tau^{n+k}_{\Lambda,\mu}(x_1,\ldots,x_{n+k}) m(dx_{n+1}) \ldots m(dx_{n+k})$$

*for $x_1, \dots, x_n \in \Lambda$. The correlation functions satisfy the Mayer equations*

$$\varrho_\mu(x_1, \dots, x_n) = \exp[-\sum_{i<j} \phi(x_j - x_i)] \times$$

$$\sum_{p=0}^{\infty} \frac{1}{p!} K((x_1, \dots, x_n), (y_1, \dots, y_p)) \varrho_\mu(y_1, \dots, y_p) m(dy_1) \dots m(dy_p),$$

*where*

$$K((x_1, \dots, x_n), (y_1, \dots, y_p)) = \prod_{j=1}^{p} K((x_1, \dots, x_n), y_j),$$

$$K((x_1, \dots, x_n), y) = \exp[-\sum_{i=1}^{n} \phi(y - x_j)] \; - \; 1.$$

*Moreover, there exists a constant $\xi > 0$ such that for all $n$ and $x_1, \dots, x_n \in \mathbb{R}^d$ we have:*

$$\varrho_\mu(x_1, \dots, x_n) \; \leq \; \xi^n.$$

**Remark 7.5** For pair potentials $\phi$ the condition in Theorem 6.13 ensuring closability of $(\mathcal{E}_\mu^\Gamma, \mathcal{F}C_b^{\infty,\mu})$ on $L^2(\Gamma; \mu)$ for $\mu \in \mathcal{G}_{gc}^1(z, \phi)$ can be now easily formulated as follows: for $\mu$–a.e. $\gamma \in \Gamma$ and $m$–a.e. $x \in \{y \in \mathbb{R}^d \mid \sum_{y' \in \gamma} |\phi(y - y')| < \infty\}$ it holds that

$$\int_{V_x} e^{\sum_{y' \in \gamma} \phi(y - y')} m(dy) < \infty$$

for some open neighbourhood $V_x$ of $x$. This condition trivially holds e.g. if $\operatorname{supp} \phi$ is compact, $\{\phi < \infty\}$ is open, and $\phi^+ \in L_{loc}^\infty(\{\phi < \infty\}; m)$. If even $\mu \in \mathcal{G}_{gc}^t(z, \phi)$ and $\phi$ satisfies the assumptions in Proposition 7.1, then it suffices to merely assume that $\{\phi < \infty\}$ is open and $\phi^+ \in L_{loc}^\infty(\{\phi < \infty\}; m)$. This follows by an elementary consideration. This generalizes the closability result in [Os96]. The a priori bigger domain for $\mathcal{E}_\mu^\Gamma$ considered there, can also be treated by our method to prove closability in Subsection 6.3.

## 7.2 Integration by parts characterization

In this subsection we shall prove an integration by parts formula for a certain convex subset of canonical Gibbs measures, to be introduced below, which contains the Ruelle measures. In fact, we shall characterize this subset by integration by parts in a way analogous to the free case in Theorem 2.2. For simplicity we assume that

(C) $\quad$ $\operatorname{supp}\phi$ is compact.

The general case is more complicated to formulate. We refer to [AKR97b, Subsection 4.3] for this case including all proofs of the results stated below.

Let us introduce another condition on the potential $\phi$:

(D) (*Differentiability*) $e^{-\phi}$ is weakly differentiable on $\mathbb{R}^d$, $\phi$ is weakly differentiable on $\mathbb{R}^d \setminus \{0\}$ and the weak gradient $\nabla\phi$ (which is a locally $m$-integrable function on $\mathbb{R}^d \setminus \{0\}$) considered as an $m$-a.e. defined function on $\mathbb{R}^d$ satisfies

$$\nabla\phi \in L^1(\mathbb{R}^d; e^{-\phi}m) \cap L^2(\mathbb{R}^d; e^{-\phi}m).$$

(We refer to [Y96] for a similar condition). Note that for many typical potentials in Statistical Physics we have $\phi \in C^\infty(\mathbb{R}^d \setminus \{0\})$. For such "regular outside the origin" potentials condition (D) nevertheless does not exclude a singularity at the point $0 \in \mathbb{R}^d$.

**Remark 7.6** (i) Suppose that (D) holds. Then $e^{-\phi}$ is weakly differentiable, hence, in particular, $e^{-\phi} \in L^1_{loc}(\mathbb{R}^d; m)$. So, (C) and (D) imply (I).

(ii) It is an easy exercise to check that (D) implies that

$$\nabla e^{-\phi} \;=\; -\nabla\phi\; e^{-\phi} \;\; m-\text{a.e. on } \mathbb{R}^d.$$

(iii) By [RSch97, Proposition 2.1] any $\mu \in \mathcal{G}^t_{gc}(z,\phi)$ satisfies condition (QI) in Subsection 4.1, provided $\phi$ satisfies (SS), (LR), (C) and (D).

**Example 7.7** A concrete example of a pair potential which is especially important in atomic and molecular physics is given for $d = 3$ by the so-called *Lennard-Jones potential*

$$\phi_{a,b}(x) := \frac{a}{|x|^{12}} - \frac{b}{|x|^6}, \quad x \in \mathbb{R}^d \setminus \{0\},$$

where $a, b > 0$. This potential obviously satisfies all conditions of Proposition 7.1 and is, therefore, super-stable and lower regular. Furthermore, we have

$$\nabla\phi_{a,b}(x) \;=\; \frac{6x}{r^8(x)} - \frac{12x}{r^{14}(x)}$$

and it is easy to see that (D) is also true. To also satisfy (C) we e.g. just have to multiply $\phi_{a,b}$ by some $\chi \in C_0^\infty(\mathbb{R}^d)$, $\chi \geq 0$.

**Lemma 7.8** *Let $\phi$ be a pair potential satisfying conditions (SS), (LR), (C) and (D). For any vector field $v \in V_0(\mathbb{R}^d)$ we consider the function*

$$\Gamma \ni \gamma \mapsto L_v^\phi(\gamma) \;:=\; -\sum_{\{x,y\}\subset\gamma} \langle \nabla\phi(x-y), v(x) - v(y)\rangle_{\mathbb{R}^d}. \qquad (7.3)$$

*Then for any Ruelle measure $\mu$ and all $v \in V_0(\mathbb{R}^d)$ we have that*

$$L_v^\phi \in L^2(\Gamma;\mu).$$

**Remark 7.9** An analysis of the proof of Lemma 7.8 (cf. [AKR97b, Lemma 4.1]) shows that condition (D) can be weakened by replacing $L^1(\mathbb{R}^d; e^{-\phi}m)$, $L^2(\mathbb{R}^d; e^{-\phi}m)$ by the weak $L^p$–spaces $L^1_w(\mathbb{R}^d; e^{-\phi}m)$, $L^2_w(\mathbb{R}^d; e^{-\phi}m)$ respectively.

In order to prove the main result of this section (i.e., Theorem 7.11 below) we need to introduce the convex set $\mathcal{P}_q, q \in [1,\infty[$, of all probability measures $\mu$ on $(\Gamma, \mathcal{B}(\Gamma))$ satisfying the following properties:

($\mathcal{P}$ 1) $\mu(S_\infty) = 1$, i.e., $\mu$ is tempered.

($\mathcal{P}$ 2) $\int \gamma(K)^q \, \mu(d\gamma) < \infty$ for all compact $K \subset \mathbb{R}^d$.

($\mathcal{P}$ 3) For all $v \in V_0(\mathbb{R}^d)$, $L_v^\phi \in L^q(\Gamma;\mu)$, where $L_v^\phi$ is as defined in (7.3).

For $\mu \in \mathcal{P}_q$ and $v \in V_0(\mathbb{R}^d)$ we define

$$B_v^\phi \;:=\; L_v^\phi + \langle \operatorname{div} v, \cdot\rangle (\in L^q(\Gamma;\mu)). \qquad (7.4)$$

and for $V = \sum_{i=1}^N F_i\, v_i \in \mathcal{VF}C_b^\infty$

$$\operatorname{div}_\phi^\Gamma V := \sum_{i=1}^N (\nabla_{v_i}^\Gamma F_i + B_{v_i}^\phi F_i). \qquad (7.5)$$

**Lemma 7.10** *$\mathcal{G}^t_{gc}(z,\phi) \subset \mathcal{G}_c(m,\Phi) \cap \mathcal{P}_2$ for all $z > 0$ (where $\Phi = \Phi_\phi$ as above).*

**Proof.** Let $z > 0$, $\mu \in \mathcal{G}_{gc}(z,\phi)$. We already know that $\mu \in \mathcal{G}_c(m,\Phi)$ by (6.11) and by Remark 7.3 that $\mu$ satisfies ($\mathcal{P}$ 2). By Lemma 7.8 also ($\mathcal{P}$ 3) with $q = 2$ holds.

$\square$.

Now we are prepared to prove the main result of this section.

**Theorem 7.11** *Let $\phi$ be a pair potential satisfying properties (SS), (LR), (C), (D). Then for $\mu \in \mathcal{P}_1$ the following assertions are equivalent:*

*(i) $\mu \in \mathcal{G}_c(m,\Phi)$ (where $\Phi = \Phi_\phi$ as above).*

*(ii) For all $v \in V_0(\mathbb{R}^d)$ and all $F, G \in \mathcal{F}C_b^\infty$ the following integration by parts formula holds:*

$$\int \nabla^\Gamma_v F \; G \; d\mu \;=\; -\int F \; \nabla^\Gamma_v G \; d\mu \;-\; \int F \; G \; B^\phi_v \; d\mu.$$

*In particular, any Ruelle measure (with interaction given by $\phi$) satisfies (ii).*

**Remark 7.12** (i) We note that if even $\mu \in \mathcal{P}_2$, then (ii) in Theorem 7.11 is equivalent to

$$(\nabla^\Gamma_v)^* \;=\; -\,\nabla^\Gamma_v \;-\; B^\phi_v \quad \text{for all } v \in V_0(\mathbb{R}^d) \tag{7.6}$$

as an operator equality on the domain $\mathcal{F}C_b^\infty$ in $L^2(\Gamma;\mu)$, where $(\nabla^\Gamma_v)^*$denotes the adjoint of $\nabla^\Gamma_v$ on $L^2(\Gamma;\mu)$.

(ii) If $\mu \in \mathcal{P}_1$, then Theorem 7.11 (ii) is obviously equivalent to:

$$\int \langle \nabla^\Gamma F, V\rangle_{T\Gamma} \, d\mu = -\int F \; \mathrm{div}^\Gamma_\phi V \, d\mu \text{ for all } F \in \mathcal{F}C_b^\infty,\ V \in \mathcal{V}\mathcal{F}C_b^\infty. \tag{7.7}$$

This shows the analogy of Theorem 7.11 with Theorem 2.2.

(iii) As mentioned before, condition (C) is not really necessary for Theorem 7.11 to hold. The more general analogue of this result, i.e., Theorem 4.3 in [AKR97b] applies to the Lennard–Jones potential $\phi_{a,b}$ in Example 7.7.

# 8 Infinite interacting particle systems

The purpose of this section is to briefly summarize all results on the stochastic dynamics associated with Gibbs measures $\mu$, i.e., on the corresponding *infinite interacting particle systems* (cf. below). By all the previous preparations in Sections 6 and 7, these results are direct applications of the general results obtained in Subsection 4.2.

## 8.1 Stochastic dynamics corresponding to Gibbs states

Let $X$ be a Riemannian manifold as before. Let $\mathcal{G}^1_{gc}(\sigma, \Phi)$ (cf. (6.12)) for a general potential $\Phi$ and $\sigma = \rho \cdot m$ as specified at the beginning of Remark 2.3 (i).

**Theorem 8.1** *Suppose that $\mu$, $\Phi$ satisfy the conditions of Theorem 6.13. Then $(\mathcal{E}^\Gamma_\mu, \mathcal{F}C_b^{\infty,\mu})$ is (well–defined and) closable on $L^2(\Gamma_X;\mu)$. The closure $(\mathcal{E}^\Gamma_\mu, D(\mathcal{E}^\Gamma_\mu))$ is a symmetric quasi–regular Dirichlet form on $L^2(\ddot{\Gamma}_X;\mu)$ and, thus has an associated conservative diffusion process* **M** *satisfying all properties in Theorems 4.14, 4.16. In particular, all this holds if $X = \mathbb{R}^d$ with the Euclidean metric and if $\mu \in \mathcal{G}^t_{gc}(z, \phi)$ for some pair potential $\phi$ satisfying (SS), (LR), (I) and the condition specified in Remark 7.5.*

**Proof.** Theorems 6.13, 4.12, 4.14, 4.16. □

## 8.2 Solutions of the martingale problem resp. of the corresponding "heuristic" SDE

Let $X = \mathbb{R}^d$ with the Euclidean metric, $\Gamma := \Gamma_{\mathbb{R}^d}$, $\ddot{\Gamma} := \ddot{\Gamma}_{\mathbb{R}^d}$. Let $\mu \in \mathcal{G}_c(m, \Phi) \cap \mathcal{P}_2$ (cf. Subsection 7.2), (e.g. $\mu \in \mathcal{G}^t_{gc}(z \cdot m, \phi)$, $z \in ]0, \infty[$ where $\Phi = \Phi_\phi$ is associated with a pair potential $\phi$).

**Theorem 8.2** *Suppose $\phi$ satisfies (SS), (LR), (C) and (D). Then:*

*(i) $\mu$ satisfies (IbP) in Subsection 4.1 and for all $F, G \in \mathcal{F}C_b^\infty$*

$$\mathcal{E}^\Gamma_\mu(F, G) = -\int \Delta^\Gamma_\phi F \; G \; d\mu$$

*where*

$$\Delta^\Gamma_\phi = div^\Gamma_\phi \nabla^\Gamma$$

*and* $div_\phi^\Gamma$ *is given by* (7.5) *(see also* (7.4), (7.3)*).*

*(ii)* $(\mathcal{E}_\mu^\Gamma, \mathcal{F}C_b^{\infty,\mu})$ *is (well–defined and) closable on* $L^2(\Gamma;\mu)$. *The closure* $(\mathcal{E}_\mu^\Gamma, D(\mathcal{E}_\mu^\Gamma))$ *is a symmetric quasi–regular Dirichlet form on* $L^2(\ddot{\Gamma};\mu)$ *and is thus associated with a conservative diffusion process* $\mathbf{M}$ *satisfying all properties in Theorems 4.14, 4.16.*

*(iii) If* $\mu \in \mathcal{G}_{gc}^t(z\cdot m, \phi)$, *then* $\mu$ *satisfies (QI), hence* $(\mathcal{E}_\mu^\Gamma, \mathcal{F}^{(c)})$ *(as defined in Proposition 4.9) also is a symmetric quasi–regular Dirichlet form on* $L^2(\ddot{\Gamma};\mu)$ *and has thus an associated conservative diffusion process* $\mathbf{M}$ *satisfying all properties in the analogues of Theorems 4.14, 4.16 with* $(\mathcal{E}_\mu^\Gamma, \mathcal{F}^{(c)})$ *replacing* $(\mathcal{E}_\mu^\Gamma, D(\mathcal{E}_\mu^\Gamma))$ *resp. with the generator of* $(\mathcal{E}_\mu^\Gamma, \mathcal{F}^{(c)})$ *replacing* $(H_\mu^\Gamma, D(H_\mu^\Gamma))$.

*(iv) If* $\mu \in \mathcal{G}_{gc}^t(z\cdot m, \phi)$, *then* $\ddot{\Gamma}\setminus\Gamma$ *is* $(\mathcal{E}_\mu^\Gamma, D(\mathcal{E}_\mu^\Gamma))$*–exceptional, hence also* $(\mathcal{E}_\mu^\Gamma, \mathcal{F}^{(c)})$*–exceptional.*

**Proof.** (i), (ii): Theorem 7.11, Proposition 4.3 (ii), Theorems 4.12, 4.14, 4.16.
(iii): By [RSch97, Proposition 2.1] $\mu$ satisfies (QI). Hence, the assertions follow from Propositions 4.9, 5.6, and Theorem 5.8.
(iv): The first part is just [RS97, Corollary 1] the second part then follows from the first (cf. the discussion in Subsection 5.3 which led to the proof of Proposition 5.6. □

Theorem 4.16 states that both in the situation of Theorem 8.2 and 8.1 the process $\mathbf{M} = (\mathbf{\Omega}, \mathbf{F}, (\mathbf{F}_t)_{t\geq 0}, (\mathbf{\Theta}_t)_{t\geq 0}, (\mathbf{X}_t)_{t\geq 0}, (\mathbf{P}_\gamma)_{\gamma\in\ddot{\Gamma}_X})$ solves the martingale problem for $(-H_\mu^\Gamma, D(H_\mu^\Gamma))$, in particular, implies that for all $F \in \mathcal{F}C_b^\infty$

$$F(\mathbf{X}_t) - F(\mathbf{X}_0) + \int_0^t H_\mu^\Gamma F(\mathbf{X}_s)\, ds\ ,\ t \geq 0, \tag{8.1}$$

is an $(\mathbf{F}_t)$–martingale under $\mathbf{P}_\gamma$ for $\mathcal{E}_\mu^\Gamma$–q.e. $\gamma \in \ddot{\Gamma}_X$. (The reader should note that if this only holds for all $F \in \mathcal{F}C_b^\infty$, rather than all $F \in D(H_\mu^\Gamma)$, it is not clear whether this uniquely determines the diffusion $\mathbf{M}$. This is in fact an open question which so far has only been answered positively if $\mu$ is a pure Poisson measure resp. a mixed Poisson measure and $(H_\sigma^X, D(H_\sigma^X))$ is conservative. The latter follows immediately by Theorem 3.3). In the case of Theorem 8.2 we can, however, write an explicit formula for $H_\mu^\Gamma F$ for all $F \in \mathcal{F}C_b^\infty$ (which is not possible in the situation of Theorem 8.1). It namely follows from Theorem 8.2(i), (2.13) and (2.23) that

$$H_\mu^\Gamma F(\gamma) = \Delta^\Gamma F(\gamma) - 2\sum_{y\in\gamma}\langle\nabla\phi(\cdot - y), \nabla^\Gamma F(\gamma)\rangle_{T\Gamma}, \quad \gamma \in \Gamma.$$

Since the closure of $(\Delta^\Gamma, \mathcal{F}C_b^\infty)$ on $L^2(\ddot{\Gamma}; \pi^z)$ generates the Brownian motion on $\mathbb{R}^d$ (cf. Theorems 3.3, 3.12 and Remark 3.15(iii)), (8.1) can be interpreted that $\mathbf{M}$ solves the following "heuristic" stochastic differential equation

$$dX_t^i = dW_t^i + \sum_{j:j\neq i} \nabla\phi(X_t^i - X_t^j)dt, \quad i \in \mathbb{N}, \tag{8.2}$$

where $(W_t^i)_{t\geq 0}$, $i \in \mathbb{N}$, are independent Brownian motions on $\mathbb{R}^d$ and $(X_t^i)_{t\geq 0}$ are such that $\mathbf{X}_t = \sum_{i=1}^\infty \varepsilon_{X_t^i}$, $t \geq 0$ (cf. [La77], [Fr87] for the "smooth case", i.e., with $\phi \in C_0^3(\mathbb{R}^d)$). The solution to (8.2) is interpreted as an infinite particle system undergoing interactions given by the drift determined by $\nabla\phi$. But, of course, (8.2) is *purely heuristic.*

# 9 Ergodicity

We now want to prove the analogue of Theorem 3.18 for $\mu \in \mathcal{G}_c(m, \Phi) \cap \mathcal{P}_2$ first proven in [AKR97b, Sect. 6]. We start, however, with a general result on ergodicity which can be formulated for large classes of probability measures $\mu$ on $(\Gamma_X, \mathcal{B}(\Gamma_X))$. In this general situation also the proofs become more transparent, since they are quite universal. Therefore, in this section we include all proofs from [AKR97b, Sect. 6].

## 9.1 A general result on irreducibility resp. ergodicity

Consider the situation studied in Subsection 4.1 after Remark 4.2, i.e., $\mu$ is a probability measure on $(\Gamma_X, \mathcal{B}(\Gamma_X))$ ($X$ a Riemannian manifold) satisfying (IbP), in particular, hence (4.3) holds and $(\mathcal{E}_\mu^\Gamma, \mathcal{F}C_b^{\infty,\mu})$ (as defined in (4.2)) is closable on $L^2(\Gamma_X; \mu)$ and its closure $(\mathcal{E}_\mu^\Gamma, D(\mathcal{E}_\mu^\Gamma))$ is a quasi–regular local symmetric Dirichlet form on $L^2(\ddot{\Gamma}_X; \mu)$ (cf. Proposition 4.3, Theorem 4.12, and Corollary 4.13). We recall that according to (4.9) we set $H_0^{1,2}(\Gamma_X; \mu) := D(\mathcal{E}_\mu^\Gamma)$.

We start with the following simple lemma.

**Lemma 9.1** *Let $(\mathcal{E}, D(\mathcal{E}))$ be a symmetric Dirichlet form on $L^2(\Gamma_X;\mu)$. Then $(\mathcal{E}, D(\mathcal{E}))$ is irreducible (i.e., $F = const.$ provided $F \in D(\mathcal{E})$ such that $\mathcal{E}(F,F) = 0$) if and only if for all* bounded *$F \in D(\mathcal{E})$ with $\mathcal{E}(F,F) = 0$ it follows that $F =$* const.

**Proof.** Suppose $F \in D(\mathcal{E})$ with $\mathcal{E}(F,F) = 0$. Then (e.g. by [MR92, Chap. I, Proposition 4.17]) we have that for all $n \in \mathbb{N}$

$$F_n := (F \wedge n) \vee (-n) \in D(\mathcal{E})$$

and $F_n \underset{n\to\infty}{\longrightarrow} F$ w.r.t. $(\mathcal{E}(\cdot,\cdot) + (\cdot,\cdot)_{L^2(\Gamma_X;\mu)})^{1/2}$. Furthermore, $\mathcal{E}(F_n, F_n) \le \mathcal{E}(F,F)$ for all $n \in \mathbb{N}$ (cf. e.g. [MR 92, Chap. I, Theorem 4.12]). Now the assertion follows easily. $\square$

We recall the following known result characterizing the irreducibility of $(\mathcal{E}^\Gamma_\mu, H^{1,2}_0(\Gamma_X;\mu))$. Let $\mathbf{M} = (\mathbf{\Omega}, \mathbf{F}, (\mathbf{F}_t)_{t\ge 0}, (\mathbf{\Theta}_t)_{t\ge 0}, (\mathbf{X}_t)_{t\ge 0}, (\mathbf{P}_\gamma)_{\gamma\in\ddot{\Gamma}_X})$ be the conservative diffusion process associated with $(\mathcal{E}^\Gamma_\mu, H^{1,2}_0(\Gamma_X;\mu))$ (cf. Theorem 4.14). As usual we set

$$\mathbf{P}_\mu := \int \mathbf{P}_\gamma \; \mu(d\gamma) \; .$$

We have the analogue of Proposition 3.17 whose proof also works in this case.

**Proposition 9.2** *The following assertions are equivalent:*

*(i)* $\mathbf{P}_\mu$ *is (time) ergodic (i.e., every bounded $\mathbf{F}$–measurable function $G : \mathbf{\Omega} \to \mathbb{R}$, which is $\mathbf{\Theta}_t$–invariant for all $t \ge 0$, is constant $\mathbf{P}_\mu$–a.e.*

*(ii)* $(\mathcal{E}^\Gamma_\mu, H^{1,2}_0(\Gamma_X;\mu))$ *is irreducible.*

*(iii)* $(e^{-H^\Gamma_\mu t})_{t>0}$ *is irreducible (i.e., if $G \in L^2(\Gamma_X;\mu)$ such that $e^{-H^\Gamma_\mu t}(GF) = G\, e^{-H^\Gamma_\mu t} F$ for all $F \in L^\infty(\Gamma_X;\mu)$, $t > 0$, then $G = const.$).*

*(iv) If $F \in L^2(\Gamma_X;\mu)$ such that $e^{-H^\Gamma_\mu t}F = F$ for all $t > 0$, then $F = const.$.*

*(v)* $(e^{-H^\Gamma_\mu t})_{t>0}$ *is ergodic (i.e.,*

$$\int \left( e^{-H^\Gamma_\mu t} F - \int F \; d\mu \right)^2 d\mu \to 0 \text{ as } t \to \infty \text{ for all } F \in L^2(\Gamma_X;\mu))$$

*(vi) If $F \in D(H^\Gamma_\mu)$ with $H^\Gamma_\mu F = 0$, then $F = const.$ ("uniqueness of ground state").*

Now we consider the weak (1,2)–Sobolev space $W^{1,2}(\Gamma_X;\mu)$ introduced in (4.11) and the corresponding bilinear form $(\mathcal{E}^\Gamma_\mu, W^{1,2}(\Gamma_X;\mu))$.

In view of Lemma 9.1 we define $(\mathcal{E}^\Gamma_\mu, W^{1,2}(\Gamma_X;\mu))$ to be *irreducible* if for all bounded $F \in W^{1,2}(\Gamma_X;\mu)$ with $\mathcal{E}^\Gamma_\mu(F,F) = 0$, it follows that $F = const.$ However, unlike as for $(\mathcal{E}^\Gamma_\mu, H^{1,2}_0(\Gamma_X;\mu))$ we *cannot* drop the term "bounded" in this definition, since it is not clear whether the larger bilinear form $(\mathcal{E}^\Gamma_\mu, W^{1,2}(\Gamma_X;\mu))$ is also a Dirichlet form. We emphasize here that the Dirichlet property was heavily used in the proof of Lemma 9.1. In particular, we do not know whether an analogue of Proposition 9.2 holds for $(\mathcal{E}^\Gamma_\mu, W^{1,2}(\Gamma_X;\mu))$

**Remark 9.3** Obviously, the irreducibility of $(\mathcal{E}^\Gamma_\mu, W^{1,2}(\Gamma_X;\mu))$ implies that of $(\mathcal{E}^\Gamma_\mu, H^{1,2}_0(\Gamma_X;\mu))$.

Suppose for every $v \in V_0(X)$ we are given a sequence $(B_{v,n})_{n\in\mathbb{N}}$ of $\mathcal{B}(\Gamma_X)$–measurable $\mathbb{R}$–valued functions on $\Gamma_X$. Let $(\mathrm{IbP})^B$ denote the set of all measures $\mu$ on $(\Gamma_X, \mathcal{B}(\Gamma_X))$ such that

(IbP 1) $B^\mu_v := L^1(\Gamma_X;\mu) - \lim_{n\to\infty} B_{v,n}$ exists and $B^\mu_v \in L^2(\Gamma_X;\mu)$ for all $v \in V_0(X)$.

(IbP 2) $\int \gamma(K)^2\ \mu(d\gamma) < \infty$ for all compact $K \subset X$ and
$$\int \nabla^\Gamma_v F\ d\mu = -\int F\ B^\mu_v\ d\mu \text{ for all } F \in \mathcal{F}C^\infty_b \text{ and all } v \in V_0(X).$$

Obviously, $(\mathrm{IbP})^B$ is a convex set. Let $\mathrm{ex}(\mathrm{IbP})^B$ denote the set of its extreme points. Clearly, if $\mu_1, \mu_2 \in (\mathrm{IbP})^B$ and $M \in ]0,\infty[$ such that $M\mu_1 - \mu_2$ is a non–zero positive measure on $(\Gamma_X, \mathcal{B}(\Gamma_X))$, then $[(M\mu_1 - \mu_2)(\Gamma_X)]^{-1}(M\mu_1 - \mu_2) \in (\mathrm{IbP})^B$. For every $\mu \in (\mathrm{IbP})^B$ both $(\mathcal{E}^\Gamma_\mu, H^{1,2}_0(\Gamma_X;\mu))$ and $(\mathcal{E}, W^{1,2}(\Gamma_X;\mu))$ are defined according to Subsection 4.1.

**Theorem 9.4** *Let $\mu \in (\mathrm{IbP})^B$. Then the following assertions are equivalent:*

*(i)* $\mu \in ex(\mathrm{IbP})^B$.

*(ii)* $\{\nu \in (\mathrm{IbP})^B \mid \nu = \rho\cdot\mu$ *for some bounded, $\mathcal{B}(\Gamma_X)$–measurable function* $\rho : \Gamma_X \to \mathbb{R}_+\} = \{\mu\}$.

*(iii)* $(\mathcal{E}^\Gamma_\mu, W^{1,2}(\Gamma_X;\mu))$ *is irreducible.*

**Proof.** (i) $\Rightarrow$ (ii): Assume (i) holds. Let $\rho : \Gamma_X \to \mathbb{R}_+$, bounded, $\mathcal{B}(\Gamma_X)$–measurable such that $\nu := \rho \cdot \mu \in (\mathrm{IbP})^B$, and let $M := \sup_{\gamma \in \Gamma_X} \rho(\gamma)$. Define

$$\mu_1 := \frac{M - \rho}{M - 1} \cdot \mu \ .$$

Then $\mu_1 \in (\mathrm{IbP})^B$ and $\mu = \frac{M-1}{M} \mu_1 + \frac{1}{M} \nu$. By assumption (i) it follows that $\mu_1 = \nu$ which implies $\rho = 1$, and (ii) is proved.

(ii) $\Rightarrow$ (i): Assume (ii) holds. Let $\mu_1, \mu_2 \in (\mathrm{IbP})^B$ and $t \in ]0, 1[$ such that $\mu = t\mu_1 + (1 - t)\mu_2$. Then they are both absolutely continuous w.r.t. $\mu$ with bounded densities. By assumption (ii) it follows that $\mu_1 = \mu = \mu_2$. Consequently, $\mu \in \mathrm{ex}(\mathrm{IbP})^B$.

(ii) $\Rightarrow$ (iii): Assume (ii) holds. Let $G \in W^{1,2}(\Gamma_X; \mu) \cap L^\infty(\Gamma_X; \mu)$ such that $\mathcal{E}^\Gamma_\mu(G, G) = 0$. We have to show that $G = const$. Since $1 \in \mathcal{F}C_b^\infty \subset W^{1,2}(\Gamma_X; \mu)$ and $d^\mu 1 = \nabla^\Gamma_\mu 1 = 0$, replacing $G$ by $G - \operatorname{essinf} G$ we may assume that $G \geq 0$, and, in addition, that $\int G \, d\mu = 1$. Define $\nu := G \cdot \mu$. Then since $d^\mu G = 0$, Proposition 4.6 and (IbP2) imply that for all $v \in V_0(X)$ and all $F \in \mathcal{F}C_b^\infty$

$$\begin{aligned} \int \nabla^\Gamma_v F \, d\nu &= \int \langle d^\mu(FG), v \rangle_{T\Gamma_X} \, d\mu \\ &= -\int F \, B^\mu_v \, d\nu \end{aligned}$$

where we used (4.10) in the last step. Hence (IbP2) holds. Since $G$ is bounded, clearly $B_{v,n} \to B^\mu_v$ as $n \to \infty$ in $L^1(\Gamma_X; \nu)$ and $B^\mu_v \in L^2(\Gamma_X; \nu)$. So, $\nu \in (\mathrm{IbP})^B$ and assumption (ii) implies that $G = 1$, hence (iii) is proved.

(iii) $\Rightarrow$ (ii): Assume (iii) holds. Let $\rho : \Gamma_X \to \mathbb{R}_+$, bounded, $\mathcal{B}(\Gamma_X)$–measurable such that $\nu := \rho \cdot \mu \in (\mathrm{IbP})^B$. Then $B_{v,n} \to B^\mu_v$ as $n \to \infty$ in $L^1(\Gamma_X; \nu)$ and $B^\mu_v \in L^2(\Gamma_X; \nu)$ for all $v \in V_0(X)$, hence by (IbP2) applied to $\nu$

$$\int \mathrm{div}_\mu V \ \rho \ d\mu = 0 \quad \text{for all } V \in \mathcal{V}\mathcal{F}C_b^\infty.$$

Hence by Remark 4.4 and (4.11) it follows that

$$\rho \in W^{1,2}(\Gamma_X; \mu) \text{ and } \ d^\mu \rho = 0,$$

i.e., $\mathcal{E}^\Gamma_\mu(\rho, \rho) = 0$. By assumption (iii) it follows that $\rho \equiv 1$ and (ii) is proved. $\square$

**Corollary 9.5** *Let $\mu \in \mathrm{ex}(\mathrm{IbP})^B$. Then $(\mathcal{E}_\mu^\Gamma, H_0^{1,2}(\Gamma_X;\mu))$ is irreducible and all equivalent assertions in Proposition 9.2 hold.*

**Proof.** Theorem 9.4 and Remark 9.3. $\square$

## 9.2 Applications to canonical Gibbs and Ruelle measures

We start with the fee case.

### (a) Mixed Poisson measures

We consider the situation described at the end of Subsection 6.2, i.e., the general "manifold–case" with $\Phi \equiv 0$. So, by Theorem 6.6 our canonical Gibbs measures $\mathcal{G}_c(\sigma)$ are exactly the mixed Poisson measures $\mu_{\lambda,\sigma}$.

So let $\sigma = \rho \cdot m$ be as in Remark 2.3 (i), i.e., we assume

(a.1) $\rho \in L^1_{loc}(X;m)$ such that $\rho^{1/2} \in H^{1,2}_{loc}(X;m)$.

Furthermore, consider the following conditions

(a.2) $\sigma(X) = \infty$ and $\beta^\sigma \in L^1_{loc}(X;m)$ where $\beta^\sigma := \dfrac{\nabla^X \rho}{\rho}$ (as in (2.21)).

(a.3) $(H_\sigma^X, \mathcal{D}^\sigma)$ is essentially self–adjoint on $L^2(X;\sigma)$ (which is e.g. the case if $X$ is complete and $|\beta|_{TX} \in L^p_{loc}(X;m)$ for some $p > \dim X$; cf. Remark 3.4).

(a.4) $(H_\sigma^X, D(H_\sigma^X))$ is conservative (cf. condition (B) in Theorem 3.3).

In conditions (a.3), (a.4) the operator $(H_\sigma^X, D(H_\sigma^X))$ is the generator of the Dirichlet form $(\mathcal{E}_\sigma^X, D(\mathcal{E}_\sigma^X))$ defined as the closure of $(\mathcal{E}_\sigma^X, D^\sigma)$ on $L^2(X;m)$ (cf. Theorem 6.12).

Let $(\mathrm{IbP})^B$ be defined as in Subsection 9.1 with $B_{v,n} := \langle \mathrm{div}_\sigma^X v, \cdot\rangle$ , $n \in \mathbb{N}$, $v \in V_0(X)$ (cf. (2.22)). Let $\mathcal{M}_2$ denote the set of all probability measures $\mu$ on $(\Gamma_X, \mathcal{B}(\Gamma_X))$ such that $\int \gamma(K)^2\,\mu(d\gamma) < \infty$ for all compact $K \subset X$.

**Lemma 9.6** *Assume that (a.1) and (a.2) hold. Then*

$$\begin{aligned}(\mathrm{IbP})^B &= \{\mu_{\lambda,\sigma} \mid \lambda \textit{ satisfies } (3.1)\}\\ &= \mathcal{G}_c(\sigma) \cap \mathcal{M}_2(\Gamma_X)\ .\end{aligned}$$

**Proof.** The first equality immediately follows from Remark 2.3 (i), the second from Theorem 6.6. □

As before, for a convex set $\mathcal{K}$ of measures on $(\Gamma_X, \mathcal{B}(\Gamma_X))$ we denote the set of its extreme points by ex$\mathcal{K}$.

**Lemma 9.7** *Assume (a.1) and (a.2) hold. Then*

$$\begin{aligned}\mathrm{ex}(\mathrm{IbP})^B &= \mathrm{ex}\{\mu_{\lambda,\sigma} \mid \lambda \textit{ satisfies } (3.1)\}\\ &= \{\pi_{z\sigma} \mid z \in [0, \infty[\}.\end{aligned}$$

**Proof.** The first equality follows from Lemma 9.6. To prove the second we note that since

$$\begin{aligned}&\{\mu_{\lambda,\sigma} \mid \lambda \text{ satisfies } (3.1)\}\\ =\ &\{\mu_{\lambda,\sigma} \mid \lambda \text{ any probability measure on } (\mathbb{R}_+, \mathcal{B}(\mathbb{R}_+))\} \cap \mathcal{M}_2(\Gamma_X),\end{aligned}$$

it is straightforward to check that

$$\begin{aligned}&\mathrm{ex}\{\mu_{\lambda,\sigma} \mid \lambda \text{ satisfies } (3.1)\}\\ \subset\ &\mathrm{ex}\{\mu_{\lambda,\sigma} \mid \lambda \text{ any probability measure on } (\mathbb{R}_+, \mathcal{B}(\mathbb{R}_+))\}.\end{aligned}$$

Since $\pi_{z\sigma} \in \mathcal{M}_2(\Gamma_X)$ for all $z \in \mathbb{R}_+$, the second equality in the assertion follows by Theorem 6.6 (i) and (ii). □

**Theorem 9.8** *Assume (a.1) and (a.2) hold and let* $\mu \in \mathcal{G}_c(\sigma) \cap \mathcal{M}_2(\Gamma_X)$ *(*$= \{\mu_{\lambda,\sigma} \mid \lambda \textit{ satisfies } (3.1)\}$*, cf. Lemma 9.6). Then the following assertions are equivalent:*

*(i)* $\mu = \pi_{z\sigma}$ *for some* $z \in [0, \infty[$.

*(ii)* $(\mathcal{E}^\Gamma_\mu, W^{1,2}(\Gamma_X; \mu))$ *is irreducible.*

*Furthermore, the assertions in Corollary 9.5 hold for* $\mu := \pi_{z\sigma}$ *and any* $z \in \mathbb{R}_+$.

**Proof.** Lemmas 9.6, 9.7 and Theorem 9.4. □

Theorem 3.18, left unproved in Subsection 3.4, is now a special case of the following result:

**Theorem 9.9** *Assume that conditions (a.1) – (a.4) hold and let $\mu$ be as in Theorem 9.8. Then $H_0^{1,2}(\Gamma_X;\mu) = W^{1,2}(\Gamma_X;\mu)$ and any of the assertions (i) – (vi) of Proposition 9.2 is equivalent to:*

*(vii)* $\mu = \pi_{z\sigma}$ *for some* $z \in [0,\infty[$.

**Proof.** Proposition 4.6 and Theorem 9.8. □

### (b) Canonical Gibbs measures

We now consider the situation described in Section 7. In particular, $X := \mathbb{R}^d$, $\Gamma := \Gamma_{\mathbb{R}^d}$. We assume that the pair potential $\phi$ satisfies (SS), (LR), (C), and (D). (For cases where (C) does not necessarily hold see [AKR97b, Subsection 6.3 (b) and (c)]).

Obviously, if $(\mathrm{IbP})^B$ is defined as in Subsection 9.1 with $B_{v,n} := B_v^\phi = L_v^\phi + \langle \operatorname{div} v, \cdot\rangle$, $n \in \mathbb{N}$, $v \in V_0(\mathbb{R}^d)$ (and $L_v^\phi$ as in (7.5)), then by Theorem 7.11

$$(\mathrm{IbP})_2^B = \mathcal{G}_c(m,\Phi)_2 \ , \tag{9.1}$$

where for a set $\mathcal{K}$ of probability measures on $(\Gamma,\mathcal{B}(\Gamma))$ we set

$$\mathcal{K}_2 := \mathcal{K} \cap \mathcal{P}_2$$

with $\mathcal{P}_2$ as in Subsection 7.2.

**Lemma 9.10** *For any convex set $\mathcal{K}$ of probability measures on* $(\Gamma,\mathcal{B}(\Gamma))$

$$\mathrm{ex}(\mathcal{K} \cap \mathcal{P}_2) = \mathrm{ex}\mathcal{K} \cap \mathcal{P}_2 \ .$$

The proof of Lemma 9.10 is straightforward, hence omitted.

**Theorem 9.11** *Let $\mu \in \mathcal{G}_c(m,\Phi)_2$. Then the following assertions are equivalent:*

*(i)* $\mu \in \mathrm{ex}\mathcal{G}_c(m,\Phi)$.

*(ii)* $(\mathcal{E}^{\Gamma}_{\mu}, W^{1,2}(\Gamma;\mu))$ *is irreducible.*

*Furthermore, the assertion in Corollary 9.5 holds for* $\mu \in \mathrm{ex}\mathcal{G}_c(m,\Phi)_2$.

**Proof.** (9.1), Lemma 9.10 and Theorem 9.4. □.

**Remark 9.12** Note that since $\emptyset \neq \mathcal{G}^t_{gc}(z,\phi) \subset \mathcal{G}_c(m,\Phi)_2$ for all $z > 0$, it immediately follows from Remark 6.5 that $\mathrm{ex}\mathcal{G}_c(m,\Phi)_2 \neq \emptyset$.

### (c) Ruelle measures

We are still in the situation of Subsection 9.2 (b).

**Theorem 9.13** *Let* $z \in ]0,\infty[$ *and* $\mu \in \mathcal{G}^t_{gc}(z,\phi)$, *i.e.,* $\mu$ *is a Ruelle measure. Then the following assertions are equivalent:*

*(i)* $\mu \in \mathrm{ex}\mathcal{G}^t_{gc}(z,\phi)$.

*(ii)* $(\mathcal{E}^{\Gamma}_{\mu}, W^{1,2}(\Gamma;\mu))$ *is irreducible.*

*Furthermore, the assertion in Corollary 9.5 holds for* $\mu \in \mathrm{ex}\mathcal{G}^t_{gc}(z,\phi)$.

**Proof.** Since $\mathcal{G}^t_{gc}(z,\phi) \subset \mathcal{G}_c(m,\Phi)_2$ by Lemma 7.10 and since $\mathrm{ex}\mathcal{G}^t_{gc}(z,\phi) \subset \mathrm{ex}\mathcal{G}_c(m,\phi)$ for $z > 0$ by [Ge79, Theorem 6.14], all assertions follow from Theorem 9.11. □

By [Ru70] there exists $z_0 \in ]0,\infty[$ such that

$$\#\mathcal{G}^t_{gc}(z,\phi) = 1 \text{ for all } z \in ]0, z_0[.$$

As an immediate consequence of Theorem 9.13 we thus obtain:

**Corollary 9.14** *Let* $z \in ]0, z_0[$ *and* $\mu \in \mathcal{G}^t_{gc}(z,\phi)$. *Then* $(\mathcal{E}^{\Gamma}_{\mu}, W^{1,2}(\Gamma;\mu))$ *is irreducible and the assertion in Corollary 9.5 holds for* $\mu$. *In particular, the corresponding stochastic dynamics started with* $\mu$ *is (time) ergodic.*

**Acknowledgment** It is a pleasure to thank the organizers for very pleasant stays and very fruitful conferences in Anogia and at the MSRI in Berkeley. We also thank the EU and the MSRI for financial support. Financial support from DFG through Project Ro1195/2–1 and the SFB–343–Bielefeld is also gratefully acknowledged.